\documentclass[10pt, reqno]{amsart}

\usepackage{hyperref} 
\usepackage{amssymb, mathrsfs,bbold}
\usepackage{amsmath, amsthm}
\usepackage{enumerate}
\usepackage{dsfont}
\usepackage{color}
\usepackage[a4paper]{geometry}
\usepackage[utf8]{inputenc}   
\usepackage{stmaryrd}
\usepackage{titlesec, wasysym}
\usepackage{appendix}
\usepackage{todonotes, fancyhdr}
\usepackage{chngcntr}
\usepackage{cancel}
\usepackage{longtable}

\usepackage{tikz-cd}
\usetikzlibrary{calc}

\usepackage{setspace}
\renewcommand{\baselinestretch}{0.99}

\usepackage[all]{xy}
\numberwithin{subsection}{section}
\numberwithin{subsubsection}{subsection}
\numberwithin{equation}{section} 

\renewcommand{\theenumi}{\alph{enumi}}
\renewcommand{\labelenumi}{\textsf{(\theenumi)}}

\allowdisplaybreaks[4]

\definecolor{macouleur}{rgb}{0.2,0.5,0.8}

\newenvironment{Dem}[1][\unskip]{%
    \begin{list}{\hspace{1.15cm}{\color{macouleur} {\it Proof #1 --}}}{   
        \setlength{\topsep}{0pt}%
        \setlength{\leftmargin}{0pt}%
        \setlength{\rightmargin}{0pt}%
        \setlength{\listparindent}{0pt}%
        \setlength{\itemindent}{0pt}%
        \setlength{\parsep}{0pt}%
        \addtolength{\leftmargin}{0pt} 
        \addtolength{\rightmargin}{0pt}%
    } \item }{\hfill {\color{macouleur} $\rhd$}\end{list}\smallskip}

\newenvironment{Dem*}[1][\unskip]{%
    \begin{list}{\hspace{0cm}{\sf \textbf{{\small Proof #1 --}}}}{   %
        \setlength{\topsep}{0pt}%
        \setlength{\leftmargin}{0pt}%
        \setlength{\rightmargin}{0pt}%
        \setlength{\listparindent}{0pt}%
        \setlength{\itemindent}{0pt}%
        \setlength{\parsep}{0pt}%
        \addtolength{\leftmargin}{20pt}%
        \addtolength{\rightmargin}{0pt}%
    } \item }{\hfill $\rhd$\end{list}\smallskip}

\renewcommand\thesection       {\arabic{section}}
\renewcommand\thesubsection    {\thesection{\boldmath $.$}\arabic{subsection}}
\renewcommand\thesubsubsection    {\thesection{\boldmath $.$}\arabic{subsection}{\boldmath $.$}\arabic{subsubsection}} 

\titleformat{\section}[block] 
{\filcenter\normalfont\sffamily\bfseries\Large}  
{{\hspace{-0.87cm}}{\color{macouleur} \thesection} \hspace{0em} {\color{macouleur} --}\vspace{0cm}}{0.5em}{\color{macouleur} {}} 

\titleformat{\subsection}[runin]
{\filcenter\normalfont\sffamily\bfseries\large}  
{{\hspace{0cm}}{\color{macouleur} \thesubsection} \hspace{0em} {\color{macouleur} --} \vspace{0.1cm}}{.2em}{\color{macouleur} {}}   
\titlespacing{\subsection}{-0pc}{1.5ex plus .1ex minus .2ex}{0pc}   

\titleformat{\subsubsection}[runin]
{\filcenter\normalfont\sffamily\bfseries}   
{\filright\sffamily{\hspace{0cm}}{\color{macouleur} \thesubsubsection}\hspace{0em} {\color{macouleur} --}}{0.4em}{\color{macouleur} {}}\titlespacing{\subsection}{-0pc}{1.5ex plus .1ex minus .2ex}{0pc}

\newtheoremstyle{mystyle}
{3pt}               
{3pt}               
{\it }                      
{}                      
{\bfseries}      
{}                      
{0.5em}                 
{\hspace{0cm}{\color{macouleur} \textit{#2 --} {\hspace{-0.2cm}}\textit{ #1}}}   

\theoremstyle{mystyle}

\newtheorem{thm}{Theorem.}   
\newtheorem*{thm*}{Theorem}

\newtheorem{cor}[thm]{{Corollary.} }
\newtheorem{lem}[thm]{{Lemma}. }
\newtheorem{prop}[thm]{{Proposition.}}
\newtheorem{defn}[thm]{{Definition.}}
\newtheorem*{rem*}{Remark.}

\newtheoremstyle{mystyle3}
{3pt}               
{3pt}               
{\it }                      
{}                      
{\bfseries}      
{}                      
{0.5em}                 
{\hspace{-0.8cm}{\textbf{\textit{#2}} --} {\hspace{-0.02cm}}{\textbf{\textit{#1}}}}

\theoremstyle{mystyle3}
\newtheorem{asmp}{Assumption.}

\newtheoremstyle{mystyle2}
{3pt}               
{3pt}               
{\it }                      
{}                      
{\sffamily}    
{}                      
{0.5em}                 
{\llap{#2 }{\it #1{\hspace{0.2cm}--}}}

\theoremstyle{mystyle2}

\newtheorem*{definition*}{Definition}
\newtheorem*{theorem*}{Theorem}
\newtheorem*{Remark*}{Remark}
\newtheorem*{lem*} {Lemma}
\newtheorem*{defn*} {Definition}
\newtheorem*{prop*} {Proposition}
\newtheorem*{cor*} {Corollary}

\newcommand{\ssk}{\smallskip}

\renewcommand{\epsilon}{\varepsilon}

\newcommand\bbC{\mathbb{C}}

\newcommand\bbE{\mathbb{E}}
\newcommand\bbF{\mathbb{F}}

\newcommand\bbN{\textbf{\textsf{N}}}
\newcommand\bbQ{\textbf{\textsf{Q}}}
\newcommand\bbR{\textbf{\textsf{R}}}
\newcommand\bbS{\mathbb{S}}
\newcommand{\bbT}{\textbf{\textsf{T}}}
\newcommand{\bbX}{\mathbb{X}}
\newcommand\bbZ{\textbf{\textsf{Z}}}

\newcommand{\mcC}{\mathcal{C}} 
\newcommand{\mcD}{\mathcal{D}}
\newcommand{\mcE}{\mathcal{E}}

\newcommand{\mcI}{\mathcal{I}}
\newcommand{\mcJ}{\mathcal{J}}

\newcommand{\mcL}{\mathcal{L}}

\newcommand{\mcN}{\mathcal{N}}

\newcommand{\mcP}{\mathcal{P}}
\newcommand{\mcQ}{\mathcal{Q}}

\newcommand\mcZ{\mathcal{Z}}

\newcommand{\bbB}{\mathbb{B}}
\newcommand{\bbU}{\mathbb{U}}
\newcommand{\bbV}{\mathbb{V}}

\newcommand{\bsp}{\mathfrak{p}}
\newcommand{\bsq}{\mathfrak{q}}
\newcommand{\bsr}{\mathfrak{r}}
\newcommand{\bsu}{\boldsymbol{u}}
\newcommand{\bsv}{\boldsymbol{v}}
\newcommand{\bsw}{\boldsymbol{w}}
\newcommand{\bsD}{\boldsymbol{D}}

\newcommand{\mfs}{\mathfrak{s}}

\newcommand{\id}{\text{\rm id}}

\makeatletter
\newcommand*{\defeq}{\mathrel{\rlap{%
                     \raisebox{0.3ex}{$\m@th\cdot$}}%
                     \raisebox{-0.3ex}{$\m@th\cdot$}}%
                     =}
\makeatother

\makeatletter
\newcommand*{\eqdef}{=\mathrel{\rlap{%
                     \raisebox{0.3ex}{$\m@th\cdot$}}%
                     \raisebox{-0.3ex}{$\m@th\cdot$}}%
                     }
\makeatother

\newcommand{\trees}{\overline{\mathbb{G}}}
\newcommand{\Trees}{\mathbb{G}}

\newcommand{\treesp}{\overline{\mathsf{T}}}
\newcommand{\Treesp}{\mathsf{T}}

\newcommand{\basis}{{\overline{\mathbb{B}}}}
\newcommand{\Basis}{\mathbb{B}}

\newcommand{\icprod}{\overline{\triangle}}
\newcommand{\iCprod}{\triangle}

\newcommand{\Cprod}{\Delta}

\newcommand{\Models}{\mathscr{M}_{\text{\rm ad}}}
\newcommand{\mModels}{\overline{\mathscr{M}_{\text{\rm ad}}^m}}

\newcommand{\spa}{\text{\rm span}}

\newcommand{\balpha}{\boldsymbol{\alpha}}

\newcommand{\tri}{|\!|\!|}

\newcommand{\rtdb}[2]
{\begin{tikzpicture}[baseline=#1pt]
\coordinate (A1) at (0,0);
#2
\end{tikzpicture}
}

\newcommand{\pol}[3]
{\coordinate (A#1) at ($(A#2)+0.4*({cos(#3)},{sin(#3)})$);}

\newcommand{\drc}[1]
{\foreach \n in {#1} \filldraw[white] (A\n) circle (1.7pt);
\foreach \n in {#1} \draw (A\n) circle (1.7pt);}
\newcommand{\drcc}[1]
{\foreach \n in {#1} \filldraw[white] (A\n) circle (2pt);
\foreach \n in {#1} \draw (A\n) circle (2pt);
\foreach \n in {#1} \fill (A\n) circle (0.8pt);}

\newcommand{\drb}[1]
{\foreach \n in {#1} \fill (A\n) circle (1.7pt);}

\newcommand{\drl}[2]
{\foreach \n in {#2} \draw (A#1)--(A\n);}
\newcommand{\drbl}[2]
{\foreach \n in {#2} \draw[line width=1.5pt] (A#1)--(A\n);}
\newcommand{\drll}[2]
{\foreach \n in {#2} \draw[double distance=0.7pt] (A#1)--(A\n);}

\newcommand{\fcst}{l}
\newcommand{\vcst}{\ell}

\begin{document}

\begin{center}
{\LARGE\sffamily{\textbf{Regularity structures for quasilinear singular SPDEs}   \vspace{0.5cm}}}
\end{center}

\begin{center}
{\sf I. BAILLEUL\footnote{I.B. acknowledges support from the CNRS and PIMS, and the ANR via the ANR-16-CE40-0020-01 grant.}, 
M. HOSHINO\footnote{M.H. acknowledges support from JSPS KAKENHI Grant Numbers 19K14556 {\color{black}and 23K12987}.} 
and S. KUSUOKA\footnote{S.K. acknowledges support from JSPS KAKENHI Grant Number 21H00988.}}
\end{center}

\vspace{1cm}

\begin{center}
\begin{minipage}{0.8\textwidth}
\renewcommand\baselinestretch{0.7} \scriptsize \textbf{\textsf{\noindent Abstract.}} We prove the well-posed character of a regularity structure formulation of the quasilinear generalized (KPZ) equation and give an explicit form for a renormalized equation in the full subcritical regime. Under the assumption that the BPHZ models associated with a non-translation invariant operator converge, we obtain a convergence result for the solutions of the regularized renormalized equations.  This conditional result covers the spacetime white noise case. 
\end{minipage}
\end{center}

\bigskip

{\color{macouleur} {\sf 
\begin{center}
\begin{minipage}[t]{11cm}
\baselineskip =0.35cm
{\scriptsize 

\center{\textbf{Contents}}

\vspace{0.1cm}

\textbf{ 1.~Introduction\dotfill 
\pageref{SectionIntro}}

\textbf{ 2.~Setting\dotfill 
\pageref{SectionSetting}}

\hspace{0.2cm} \textbf{ 2.1~Function spaces\dotfill \pageref{SubsectionFunctionSpaces}}

\hspace{0.2cm} \textbf{ 2.2~The regularity structure\dotfill \pageref{SubsectionRS}}

\hspace{0.2cm} \textbf{ 2.3~Models and modelled distributions\dotfill \pageref{SubsectionModels}}

\textbf{ 3.~Local well-posedness\dotfill  \pageref{SectionLocalWellPosed}}

\textbf{ 4.~Renormalization matters\dotfill  \pageref{SectionRenormalization}}

\hspace{0.2cm} \textbf{ 4.1~Notations\dotfill  \pageref{SectionNotations}}

\hspace{0.2cm} \textbf{ 4.2~Coherence and morphism property for the $\star$ product\dotfill \pageref{SubsectionCoherenceStarProduct}}

\hspace{0.2cm} \textbf{ 4.3~Strong preparation maps and their associated models\dotfill \pageref{SubsectionPreparationMaps}}

\hspace{0.2cm} \textbf{ 4.4~Renormalized equation\dotfill \pageref{SubsectionRenormalizedEquation}}

\hspace{0.2cm} \textbf{ 4.5~Examples satisfying Assumption \ref{asmp1}\dotfill \pageref{SubsectionAsmp1}}

\hspace{0.2cm} \textbf{ 4.6~Examples satisfying Assumption \ref{asmp2}\dotfill \pageref{SubsectionAsmp2}}

\textbf{ A.~Appendix\dotfill   \pageref{SectionAppendix}}

\hspace{0.2cm} \textbf{ A.1.~Gaussian kernels\dotfill   \pageref{SubsectionGaussian}}

\hspace{0.2cm} \textbf{ A.2.~Existence of the fundamental solution\dotfill   \pageref{SubsectionFundamental}}

\hspace{0.2cm} \textbf{ A.3.~Uniqueness of the fundamental solution\dotfill   \pageref{AppendixUniqueness}}

\hspace{0.2cm} \textbf{ A.4.~Temporally homogeneous operator\dotfill   \pageref{AppendixProperties}}

\hspace{0.2cm} \textbf{ A.5.~Anisotropic Taylor formula\dotfill   \pageref{AppendixTaylor}}


}\end{minipage}
\end{center}
}}   \vspace{1cm}

\section{Introduction}
\label{SectionIntro}

Denote by $\bbT$ the one dimensional torus.
We consider the one dimensional space-periodic quasilinear generalized (KPZ) equation 
\begin{equation} \label{EqQgKPZ}
\big(\partial_t - a(u)\partial_x^2 \big) u = f(u)\xi + g(u)(\partial_x u)^2,
\end{equation}
for regular enough functions $a,f,g$, where $a$ takes values in a compact interval of $(0,\infty)$, the initial condition $u_0\in C^{0+}(\bbT)\defeq\bigcup_{\alpha>0}C^\alpha(\bbT)$ is given, and $\xi$ is a random spacetime distribution.
Especially, we are interested in the spacetime white noise.
This equation falls within the class of singular stochastic partial differential equations (SPDEs) of parabolic type. All equations of this class share the same defect: The low regularity of some terms in a singular SPDE prevents the expected regularizing effect of the dynamics from giving sense to a number of products in the equations. In the case at hand, Equation \eqref{EqQgKPZ}, one expects a parabolic type dynamics to have a resolvent that improves regularity by $2$. The `subcritical' nature of the dynamics is here encoded in the fact that the spacetime distribution $\xi$ is (almost surely) assumed to have regularity $\alpha_0-2$, for $0<\alpha_0<2$. It is then formally consistent to expect a solution $u$ of Equation \eqref{EqQgKPZ} to have parabolic regularity $\alpha_0$, as $(\partial_xu)^2$ will then have regularity $2(\alpha_0-1)$, bigger than $\alpha_0-2$, the expected regularity of the term $f(u)\xi$. With the right hand side of regularity $\alpha_0-2$ a Schauder type continuity estimate satisfied by the resolvent of the evolution gives indeed $u$ a regularity $\alpha_0$. The problem with that regularity analysis is that for $u$ of regularity $\alpha_0$ none of the products $f(u)\xi$ and $\vert\partial_xu\vert^2$ make sense, even less $g(u)\vert\partial_xu\vert^2$, when $0<\alpha_0<1$, the case of interest.

The development of the study of semilinear subcritial singular SPDEs was launched by the two groundbreaking works \cite{Hai14} of M. Hairer, on regularity structures, and \cite{GIP} of M. Gubinelli, P. Imkeller \& N. Perkowski, on paracontrolled calculus. Both of them introduced new settings and new tools to make sense of such equations and solve them uniquely under some small parameter condition. Despite the difference of languages and tools used in regularity structures and paracontrolled calculus both settings provide a similar understanding of a subcritical singular parabolic SPDE. The mantra of their common approach to the product problem is that if one can make sense of a number of analytically ill-defined `reference products' that only involve the noise $\xi$, not in an $\omega$-wise sense but as random variables, then one can make sense of the ill-defined products in the equation for all functions $u$ that locally look like linear combinations of the reference random variables. Regularity structures and paracontrolled calculus differ in the tools used to make sense of that comparison with reference random variables. In both settings, working with a random noise turns out to be crucial to construct these reference random variables by probabilistic means.

We refer the reader to the overviews \cite{ChandraWeber, CorwinShen} of Chandra \& Weber and Corwin \& Shen for non-technical introductions to the domain of semilinear singular SPDEs, to the books \cite{FrizHairerBook, Berglund} of Friz \& Hairer and Berglund for a mildly technical introduction to regularity structures, and to Bailleul \& Hoshino's Tourist's Guide \cite{RSGuide} for a thorough tour of the analytic and algebraic sides of the theory. Readers interested in paracontrolled calculus will find a nice account of the fundamentals in Gubinelli's panorama \cite{GubinelliPanorama}.

\medskip

The first works on quasilinear singular SPDEs by Otto \& Weber \cite{OttoWeber}, Furlan \& Gubinelli \cite{FurlanGubinelli} and Bailleul, Debussche, \& Hofmanov\'a \cite{BDH} all three investigated the generalized (PAM) equation in the regime where the noise is $(\alpha_0-2)$ regular and $\alpha_0>2/3$. Interestingly each of these works used a different method: A variant of regularity structures in \cite{OttoWeber}, a variant of paracontrolled calculus based on the use of the paracomposition operator for \cite{FurlanGubinelli}, and the initial form of paracontrolled calculus in \cite{BDH}. On the paracontrolled side Bailleul \& Mouzard \cite{BailleulMouzard} extended the high order paracontrolled calculus toolbox to deal with the paracontrolled equivalent of Equation \eqref{EqModifiedQgKPZ} in the spacetime white noise regime $\alpha_0>2/5$. On the regularity structures side Otto \& Weber deepened their framework in their works \cite{OSSW2} with Sauer \& Smith, dedicated to the study of the equation with linear additive forcing
\begin{equation} \label{EqQgKPZLinear}
\partial_t u - a(u)\partial_x^2u = \xi.
\end{equation}
They obtained in particular in \cite{OSSW2} an explicit form of a renormalized equation for \eqref{EqQgKPZLinear} backed up by the general convergence result proved by Linares, Otto, Tempelmayr \& Tsatsoulis in \cite{LOTT} that holds for a large class of random noises in the full subcritical regime. Our general formula for the counterterm in the renormalized equation generalizes theirs. The algebraic machinery behind their approach was further analysed by Linares, Otto \& Tempelmayr in \cite{LOT}. Meanwhile Gerencs\'er \& Hairer provided in \cite{GerencserHairer} an analysis of a regularity structure counterpart of Equation \eqref{EqQgKPZ}, in the full subcritical regime. Their method allowed for an analysis of the renormalized equation only in the regime $\alpha_0>1/2$. By implementing some tricky integration by parts-type formulas Gerencs\'er was able in \cite{Gerencser} to obtain the renormalized equation for the special case of Equation \eqref{EqQgKPZLinear} from the analysis of \cite{GerencserHairer} in the spacetime white noise regime $\alpha_0>2/5$. We prove in the present work the well-posed character of a regularity structure formulation of the quasilinear generalized (KPZ) equation and give an explicit form for a renormalized equation in the full subcritical regime, with a simple expression in a number of cases. 
Our result is about the algebraic part and the analytic part of the theory of regularity structures.
While we strongly expect that the probabilistic convergence result (\textbf{\textit{Assumption \ref{asmp1}}} below) holds, we leave the complete proof as a future work.

\medskip

Following \cite{BDH, BailleulMouzard} we set 
$$
L^{v} \defeq  a(v)\partial_x^2
$$
for a sufficiently regular function $v$ on $[0,\infty)\times\bbT$ and rewrite Equation \eqref{EqQgKPZ} in the form
\begin{equation} \label{EqQgKPZBis}
\big(\partial_t - L^{v}+c\big) u = f(u)\xi + g(u)(\partial_x u)^2 +cu+\big(a(u)-a(v)\big)\partial_x^2u
\end{equation}
for a large positive constant $c$. We consider \eqref{EqQgKPZBis} as a `perturbation' of the non-translation invariant generalized (KPZ) equation
$$
\big(\partial_t - L^{v}+c\big) u = f(u)\xi + g(u)(\partial_x u)^2+cu.
$$
In Section {\sf \ref{SectionSetting}} below we will reformulate Equation \eqref{EqQgKPZBis} to an equation of modelled distributions, the detailed form of which is \eqref{EqModifiedQgKPZzeta}. We call this reformulated equation the modelled form of Equation \eqref{EqQgKPZBis}. We will see in Theorem \textit{\ref{ThmMainAnalytic}} that given any admissible model $\sf M$ on our regularity structure, the modelled form of Equation \eqref{EqQgKPZBis} has a unique solution over a model-dependent time interval $(0,t_0({\sf M}))$, in an appropriate class of modelled distributions. This analytical statement holds in the full subcritical range provided the model is part of the data. Such a statement was already proved by Gerencs\'er \& Hairer in \cite{GerencserHairer} in a different setting. However their choice of formulation for \eqref{EqQgKPZ} did not allow them to write down in the full subcritical range the renormalized equation. The spacetime white noise regime is in particular out of range of their result. The regularity structure in which we formulate the modelled form of Equation \eqref{EqQgKPZBis} is different from the regularity structure used in \cite{GerencserHairer}. Working with an appropriate choice of model $\sf M$ that is the natural analogue in our setting of the Bruned-Hairer-Zambotti (BHZ) renormalized model from \cite{BHZ} we are able to give in Theorem  \textit{\ref{ThmMainRenormalized}} below a renormalized equation in the full subcritical regime. 

Denote by $\varepsilon\in(0,1]$ a regularization parameter and by $\xi^\varepsilon\in C^\infty(\textbf{\textsf{R}}\times\textbf{\textsf{T}})$ an $\varepsilon$-regularized noise $\xi$. 
Our main results take a conditional form involving two `assumptions'. \textbf{\textit{Assumption \ref{asmp1}}} is stated in Section {\sf\ref{SectionAdmissibleModelsPrepMaps}} and assumes the convergence of the natural BPHZ model ${\sf M}^\varepsilon$ associated with $\xi^\epsilon$ and the non-translation invariant operator $(\partial_t - L^v+c)$ as $\varepsilon$ goes to $0$.
We expect that under some assumptions similar to previous studies, \textbf{\textit{Assumption \ref{asmp1}}} is satisfied -- see Section {\sf\ref{SubsectionAsmp1}} for some examples. But we need a lot of technical calculations to extend a graphical approach \cite{CH} or an inductive approach \cite{HS, RandomModel} to our non translation-invariant setting. 

\medskip

\begin{thm} \label{ThmMainRenormalized}
Fix $\alpha_0\in(0,1)$, $f,g\in C_b^n(\bbR)$ with an integer $n\ge 2/\alpha_0+1$, and $a\in C_b^1(\bbR)$ such that $\inf a>0$.
Let $\alpha\in(0,\alpha_0)$ and $u_0\in C^\alpha(\bbT)$.
Choose any function $v$ on $[0,\infty)\times\bbT$ sufficiently close to $u_0$ -- see Condition \eqref{eq:neighborhoodu_0} for the precise meaning. 
Under \textbf{Assumption \ref{asmp1}}, there exist an infinite set of rooted decorated trees $\Basis_-$, a family of positive integers $\{S(\tau)\}_{\tau\in\Basis_-}$, a family of continuous functions $\{\frak{F}_v[\tau]\in C(\bbR^2)\}_{\tau\in\Basis_-}$, and a family of continuous functions $\{\ell_v^\varepsilon(\cdot,\tau)\in C(\bbR\times\bbT)\}_{\tau\in\Basis_-}$, such that the solution $u^\varepsilon$ to
\begin{equation} \label{EqRenormalizedEquationInfinite} \begin{split}
\big(\partial_t - a(u^\varepsilon)\partial_x^2\big) u^\varepsilon = f(u^\varepsilon)\xi^\varepsilon &+ g(u^\varepsilon)(\partial_x u^\varepsilon)^2    
+ \sum_{\tau\in\bbB_-} \frac{\vcst^\varepsilon_{v}(\cdot,\tau)}{S(\tau)} \, \frak{F}_v[\tau]\big(u^\varepsilon,\partial_xu^\varepsilon\big)
\end{split} \end{equation}
starting from $u_0$ converges in $C\big([0,t_0)\times\bbT\big)$ for a random time $t_0>0$ in probability as $\varepsilon>0$ goes to $0$. The last infinite series of \eqref{EqRenormalizedEquationInfinite} absolutely converges in $[0,t_0)\times\bbT$ for each $\varepsilon$. 
\end{thm}

\medskip

We will see along the way that $\frak{F}_v[\tau]$ is at most linear with respect to $\partial_xu^\varepsilon$ and that $\ell_v^\varepsilon(\cdot,\tau)$ is a non-local functional of the function $a(v(\cdot))$ independent of $u^\varepsilon$. Therefore the sum in \eqref{EqRenormalizedEquationInfinite} is not a trivial term like $-g(u^\varepsilon)(\partial_xu^\varepsilon)^2$. This sum is named {\it counterterm}. Theorem {\it \ref{ThmMainRenormalized}} extends the results of \cite{OttoWeber, FurlanGubinelli, BDH, BailleulMouzard, OSSW2, GerencserHairer} and deals with the quasilinear generalized (KPZ) in the full subcritical regime. The reader familiar with regularity structures will see that our arguments extend immediately to coupled systems of generalized (KPZ) equations \cite{BGHZ, FH17}. Such a generalization is left to the reader and we concentrate here on the renormalized equation for \eqref{EqQgKPZ}. 

In Section {\sf\ref{SubsectionAsmp1}} we discuss \textbf{\textit{Assumption 1}} by some examples. Such a convergence result was proved in several works \cite{CH,HS,RandomModel} in semilinear settings. Whether one way or another -- using a direct graphical approach as in \cite{CH} or an inductive proof as in \cite{HS, RandomModel}, we expect that the convergence result holds true in the region
$$
\alpha_0>\frac14.
$$
Otherwise the counterterm like \eqref{EqRenormalizedEquationInfinite} is not valid because for the solution $Z^\varepsilon$ of the linear stochastic heat equation $(\partial_t-\partial_x^2)Z^\varepsilon=\xi^\varepsilon$, the variance of $(\partial_xZ^\varepsilon)^2$ diverges. See \cite{Hos16} about such restriction in a semilinear setting and \cite{KPZrougher} about a different type renormalization when $\alpha_0<1/4$.


\ssk

The reader may feel uncomfortable about the fact that the functions $\ell^\epsilon_v(\cdot,\tau^{\bsp})$ depend on the somewhat arbitrary choice of function $v$ satisfying the condition \eqref{eq:neighborhoodu_0} below. One can give a simpler representation of the counterterm when the functions $\vcst^\varepsilon_v(\cdot,\tau^{\bsp})$ can be traded off for a local functional of $a(v(\cdot))$ -- meaning that $\vcst^\varepsilon_v(z,\tau^{\bsp})$ can be replaced by a function of $a(v(z))$. This is the content of {\it \textbf{Assumption \ref{asmp2}}} stated in Section {\sf\ref{SubsectionRenormalizedEquation}}.

\medskip

\begin{thm} \label{ThmMainRenormalizedSimplified}
Under \textbf{Assumption \ref{asmp1}} and \textbf{Assumption\textbf{\ref{asmp2}}} there exist a finite subset $\bbF_-\subset\bbB_-$ and a family of continuous functions $\{\fcst^\varepsilon_{\boldmath{.}}(\tau)\in C(\bbR)\}_{\tau\in\bbF_-}$ independent of $u^\varepsilon$ and $v$ such that the counterterm in \eqref{EqRenormalizedEquationInfinite} is of the form
\begin{equation} \label{EqRenormalizedEquation}
\sum_{\tau\in\bbF_-} \frac{\fcst^\varepsilon_{a(u^\varepsilon)}(\tau)}{S(\tau)} \, \frak{F}[\tau](u^\varepsilon,\partial_xu^\varepsilon) + O(1)
\end{equation}
with $O(1)$ independent of $\varepsilon$, and some functions $\frak{F}[\tau]$ that do not depend on $v$.
\end{thm}

\medskip

Note that apart from the $O(1)$ term in \eqref{EqRenormalizedEquation}, which we can discard in the renormalized equation, the counterterm is independent of $v$. We give in Section {\sf\ref{SubsectionAsmp2}} some examples of quasilinear equations where \textbf{\textit{Assumption 2}} holds true.

\medskip

Dealing with quasilinear singular SPDEs rather than with semilinear equations requires a twist that appears in the form of an infinite dimensional ingredient. It is related in our formulation \eqref{EqModifiedQgKPZ} to the fact that our structure needs to be stable under the operator $\mcI_{(0,2)}$. In the previous works using regularity structures \cite{OttoWeber, GerencserHairer, OSSW2} this infinite dimensional feature appeared under the form of a one parameter family of heat kernels or abstract integration operators. Our regularity structure is different from the regularity structures used in these works. Its model space $T=\bigoplus_{\beta\in A}T_\beta$ has infinite dimensional homogeneous spaces $T_\beta$ whose basis elements are the usual trees associated with the generalized (KPZ) equation with an additional integer decoration $p$ on each edge accounting for how many times the operator $\mcI_{(0,2)}$ is applied to this edge. The same infinite dimensional ingredient appeared in Bailleul \& Mouzard's work \cite{BailleulMouzard} in a paracontrolled setting.

\medskip

\noindent \textbf{\textit{Organization of the work} --} We set the scene in Section {\sf\ref{SectionSetting}}, where the function spaces we work with are introduced together with our regularity structure. We introduced in particular a non-classical spacetime elliptic operator to define our parabolic spaces. For the reader's convenience some properties of its heat kernel are proved in full detail in Appendix {\sf \ref{SectionAppendix}}. Section {\sf\ref{SectionLocalWellPosed}} is dedicated to proving the local well-posedness of the modelled form of Equation \eqref{EqQgKPZBis} in the full subcritical regime. The analysis of the renormalized equation is done in Section {\sf\ref{SectionRenormalization}}, where we give in particular an explicit description of the renormalization counterterm.

\bigskip

\noindent \textbf{\textit{Notations} --} {\it 
We use the sign $a\defeq b$ to define $a$ as being equal to $b$. We denote by $\bbN,\bbZ,\bbQ,\bbR$ the sets of all nonnegative integers, integers, rational numbers, real numbers, respectively. We define the the parabolic scaling as the vector 
$$
\mfs = (\mfs_1,\mfs_2) = (2,1) \in\bbN^2.
$$
We represent by $z=(t,x)\in\bbN^2$ a generic spacetime variable, for which we set 
$$
\|z\|_{\mfs} \defeq  |t|^{\frac{1}{2}}+|x|.
$$ 
We use bold letters ${\bf k},{\bf l},{\bf m},{\bf n}$ to represent multiindices in $\bbN^2$, and write ${\bf e}_1=(1,0)$ and ${\bf e}_2=(0,1)$.
For any multiindex ${\bf k}=(k_1,k_2)\in\bbN^2$ we set
$$
{\bf k}!\defeq k_1!k_2!,\qquad
|{\bf k}|_{\mfs} \defeq 2k_1+k_2,\qquad
\partial_z^{\bf k} \defeq \partial_t^{k_1}\partial_x^{k_2}.
$$  
For ${\bf k}=(k_1,k_2),{\bf l}=(l_1,l_2)\in\bbN^2$, we write ${\bf l}\le{\bf k}$ if $l_1\le k_1$ and $l_2\le k_2$.
We set 
$$
{{\bf k}\choose{\bf l}}\defeq{\bf1}_{{\bf l}\le{\bf k}}\, \frac{{\bf k}!}{{\bf l}!({\bf k}-{\bf l})!},
$$
where and in what follows, for any condition $C(x_1,\dots,x_n)$ of variables $x_1,\dots,x_n$, we denote by ${\bf1}_C$ the indicator function of the set of all $(x_1,\dots,x_n)$ satisfying $C(x_1,\dots,x_n)$.

For any $\alpha>0$, we define $C^\alpha(\bbT)$ as the collection of functions $f$ on $\bbT$ which are $n$-th differentiable and such that $f^{(n)}$ is $(\alpha-n)$-H\"older continuous, where $n$ is the unique integer such that $n<\alpha\le n+1$. We also define $\mcC_{\mfs}^\alpha(\bbR\times\bbT)$ as the set of functions $f$ on $\bbR\times\bbT$ such that $\partial_z^{\bf k}f$ exists and is bounded for any ${\bf k}\in\textbf{\textsf{N}}^2$ with $|{\bf k}|_{\mfs}<\alpha$, and for any ${\bf k}\in\bbN^2$ and $i\in\{1,2\}$ such that $|{\bf k}|_\mfs<\alpha\le|{\bf k}+{\bf e}_i|$, and for any $z\in\bbR^2$ and $h\in\bbR$, one has
$$
\big|\partial_z^{\bf k}f(z+h{\bf e}_i)-\partial_z^{\bf k}f(z)\big| \lesssim |h|^{\frac{\alpha-|{\bf k}|_\mfs}{\mfs_i}}.
$$
}

\medskip

\section{The setting}
\label{SectionSetting}

In this section we introduce the functional setting and the regularity structure in which we study the modelled form of Equation \eqref{EqQgKPZBis}.

\subsection{Function spaces{\boldmath $.$} \hspace{0.1cm}}
\label{SubsectionFunctionSpaces}

In what follows we sometimes regard a function on $\bbR\times\bbT$ as a function on $\bbR^2$ which is periodic with respect to the second variable.
The following basic facts are proved in Appendix {\sf \ref{SectionAppendix}} -- see Theorem {\it \ref{thm:GaussGamma}}, Theorem {\it \ref{thm:GammaUnique}}, and Corollary {\it \ref{cor:TaylorFundamentalSolution}}. Pick $\alpha\in(0,1]$ and an arbitrary function $v\in \mcC_{\mfs}^\alpha(\bbR\times\bbT)$. Recall the definition of the operator $L^v = a(v)\partial_x^2$.

\medskip

\begin{prop}\label{2.1QGauss}
The fundamental solution $Q_{ts}^{v,0}(x,y)$ of the operator $\partial_t-L^{v}$ satisfies the estimate
\begin{equation} \label{EqHeatKernelEstimate}
\big|\partial_t^n\partial_x^kQ_{ts}^{v,0}(x,y)\big| \le \frac{c_0e^{c_0 (t-s)}}{(t-s)^{(1+k+2n)/2}} \, \exp\bigg(\hspace{-0.1cm}-c_1\frac{|x-y|^2}{t-s} \bigg)
\end{equation}
for any $k,n\in\bbN$ such that $k+2n\le 2$, $x,y\in\bbR$, and $s<t\in\bbR$, for some positive constants $c_0,c_1$ depending only on $\inf a>0$, $\|a\|_{C_b^1(\bbR)}$, and $\|v\|_{\mcC_{\mfs}^\alpha(\bbR\times\bbT)}$.
Moreover one has $\int_{\bbR}Q_{ts}^{v,0}(x,y)dy=1$.
\end{prop}

\medskip

Here and below a number of objects are indexed by $v$. A better choice would have been to index them by $a(v)$, which would have resulted in longer expressions. We choose for that purpose to use the index $v$. We also define the spacetime elliptic operator
\begin{align}\label{EqSpacetimeOperator}
\mcL^{v} \defeq \big(\partial_t-L^{v}\big)(\partial_t+\partial_x^2) = \partial_t^2-a(v)\partial_x^4 - \big(a(v)-1\big)\partial_t\partial_x^2
\end{align}
on $\bbR^2$ with spatially periodic coefficients, and introduce the parabolic operator 
$$
\partial_\theta-\mcL^v
$$ 
with the additional variable $\theta>0$. For $C_1>0$ and $\theta>0$, we define the function
\begin{equation} \label{EqDefnGaussianLike}
{\sf G}_\theta^{(C_1)}(z)\defeq
\frac1{\theta^{3/4}}\exp\bigg\{\hspace{-0.1cm} - C_1\bigg(\frac{t^2}{\theta}+\frac{|x|^{4/3}}{\theta^{1/3}}\bigg)\bigg\}
\end{equation}
of $z=(t,x)\in\bbR^2$.

\medskip

\begin{prop}
The fundamental solution $\mcQ_\theta^{v,0}(z,w)$ of $\partial_\theta-\mcL^{v}$ on $(0,+\infty)\times\bbR^2$ satisfies the estimates
\begin{equation} \label{EqSpacetimeHeatKernelEstimate}
\big|\partial_z^{\bf k}\mcQ_\theta^{v,0}(z,w)\big| \le \frac{C_0e^{C_0\theta}}{\theta^{|{\bf k}|_{\mfs}/4}} {\sf G}_\theta^{(C_1)}(z-w)
\end{equation}
and
\begin{align} \label{EqSpacetimeHeatKernelEstimate2}
\begin{aligned}
\Big|\partial_z^{\bf k}\mcQ_\theta^{v,0}(z',w) - \sum_{|{\bf k}+{\bf l}|_{\mfs}\le 4} &\frac{(z'-z)^{{\bf k}+{\bf l}}}{{\bf l}!}\partial_z^{{\bf k}+{\bf l}}\mcQ_\theta^{v,0}(z,w)\Big|   \\
&\le \frac{C_0e^{C_0\theta}\|z'-z\|_{\mfs}^{4-|{\bf k}|_{\mfs}+\delta}}{\theta^{(4+\delta)/4}} \, \Big\{{\sf G}_\theta^{(C_1)}(z'-w)+{\sf G}_\theta^{(C_1)}(z-w)\Big\}
\end{aligned}
\end{align}
for any $|{\bf k}|_{\mfs}\le4$ and $\delta\in(0,\alpha)$, with some positive constants $C_0,C_1$ depending only on $\inf a>0$, $\|a\|_{C_b^1(\bbR)}$ and $\|v\|_{\mcC_\mfs^\alpha(\bbR\times\bbT)}$. Moreover one has $\int_{\bbR^2}\mcQ_\theta^{v,0}(z,w)dw=1$.
\end{prop}

\medskip

Note that $Q_{ts}^{v,0}$ and $\mcQ_\theta^{v,0}$ are spatially periodic in the sense that
$$
Q_{ts}^{v,0}(x+k,y)=Q_{ts}^{v,0}(x,y-k),\qquad
\mcQ_\theta^{v,0}(z+k{\bf e}_2,w)=\mcQ_\theta^{v,0}(z,w-k{\bf e}_2)
$$
for each $k\in\bbZ$. We can see these properties from the proof of Theorem {\it \ref{thm:GaussGamma}} by using the periodicity of $a(v)$.

Next we describe a class of possible choices for $v$. Recall that $\alpha_0-2$ is the spacetime H\"older regularity of the noise $\xi$ in Equation \eqref{EqQgKPZ}. We consider some initial condition $u_0\in C^\alpha(\bbT)$ with $\alpha\in(0,\alpha_0)$.

\medskip

\begin{defn} \label{def:Valpha}
For any $\alpha\in(0,1)$ and $T>0$, define $V^\alpha[0,T]$ as the set of $f\in\mcC_\mfs^\alpha([0,T]\times\bbT)$ such that the following quantity is finite.
\begin{align*}
\|f\|_{V^\alpha[0,T]}
\defeq
&\sup_{z,z'\in[0,T]\times\bbT,\,z\neq z'}\frac{|f(z')-f(z)|}{\|z'-z\|_{\mfs}^\alpha}   
+ \sum_{1\le|{\bf k}|_\mfs\le2}\sup_{t\in(0,T]}t^{(|{\bf k}|_\mfs-\alpha)/2} \, \|\partial_z^{\bf k}f(t,\cdot)\|_{L^\infty(\bbT)}   \\
&\quad+ \sup_{0<t<t'\le T}t^{(2-\alpha)/2} \, \frac{\|\partial_xf(t',\cdot)-\partial_xf(t,\cdot)\|_{L^\infty(\bbT)}}{|t'-t|^{\frac12}}.
\end{align*}
\end{defn}

\medskip

We will choose later a function $v\in V^\alpha[0,T]$ satisfying
\begin{align}\label{eq:neighborhoodu_0}
\left\{
\begin{aligned}
&\|v\|_{V^\alpha[0,T]}\le C\|u_0\|_{C^\alpha(\bbT)},\\
&\big\| e^{t\partial_x^2}u_0 - v \big\|_{L^\infty([0,T]\times\bbT)}\le \delta\|u_0\|_{C^\alpha(\bbT)}
\end{aligned}
\right.
\end{align}
for some constant $C>0$ and a sufficiently small positive constant $\delta$, which will be chosen later depending only on $\|u_0\|_{C^\alpha(\bbT)}$.
Here $\{e^{t\partial_x^2}\}_{t\ge0}$ denotes the heat semigroup associated with $\partial_x^2$.
The most natural choice is $v(t,x)=(e^{t\partial_x^2}u_0)(x)$ -- see Proposition {\it \ref{smoothingofQt1}} for the proof that $e^{(\cdot)\partial_x^2}u_0\in V^\alpha[0,T]$.
Other than that, we choose in Section \ref{subsubsection:qgKPZ} a $t$-independent smooth function $v(x)=(e^{\kappa\partial_x^2}u_0)(x)$ for a sufficiently small $\kappa>0$.
We then extend the domain of $v$ to $\bbR\times\bbT$ setting
$$
v(t,x)=
\begin{cases}
v(0,x), & \textrm{for } t\le 0,   \\
v(T,x), & \textrm{for } t\ge T.
\end{cases}
$$
and consider the spacetime operator \eqref{EqSpacetimeOperator}. Since $\|v\|_{\mcC_{\mfs}^\alpha(\bbR\times\bbT)}\lesssim\|u_0\|_{C^\alpha(\bbT)}$, the constants $c_0,c_1, C_0,C_1$ in \eqref{EqHeatKernelEstimate}, \eqref{EqDefnGaussianLike}, \eqref{EqSpacetimeHeatKernelEstimate} can be chosen to depend only on $\inf a>0$, $\|a\|_{C_b^1(\bbR)}$, and $\|u_0\|_{C^\alpha(\bbT)}$. Therefore, {\bf \textit{all the mutliplicative constants appearing sometimes implicitly in some inequalities below are independent of the choice of $v$.}}

\ssk

For any bounded continuous function $f$ on $\bbR^2$ and ${\bf k}\in\bbN^2$ with $|{\bf k}|_{\mfs}\le4$, we set
$$
\big(\partial^{\bf k}\mcQ_\theta^{v,0}f \big) (z) \defeq \int_{\bbR^2} \partial_z^{\bf k}\mcQ_\theta^{v,0}(z,w) f(w)dw.
$$
We use the operators $\mcQ_\theta^{v,0}$ to define the full scale of anisotropic parabolic H\"older spaces.

\medskip

\begin{defn*}
For $\beta<0$, we define $\mcC_{\mfs}^{\beta,v}(\bbR^2)$ as the completion of the set of bounded continuous functions $f$ on $\bbR^2$ under the norm
$$
\|f\|_{\mcC_{\mfs}^{\beta,v}(\bbR^2)}\defeq \sup_{0<\theta\le1}\theta^{-\beta/4} \big\|\mcQ_\theta^{v,0}f\big\|_{L^\infty(\bbR^2)}.
$$
We also define the closed subspace $\mcC_\mfs^{\beta,v}(\bbR\times\bbT)$ as the completion of the set of bounded continuous functions $f$ on $\bbR\times\bbT$, which are regarded as spatially periodic functions on $\bbR^2$, under the above norm.
\end{defn*}

%
%
%

\medskip

Next we rewrite the resolvent of the operator $\partial_t-L^{v}$ in terms of the operators $\mcQ_\theta^{v,0}$. With an eye on the heat kernel estimates \eqref{EqHeatKernelEstimate} and \eqref{EqSpacetimeHeatKernelEstimate} pick a positive constant $c>c_0\vee C_0$ and write 
$$
Q_{ts}^{v,c} \defeq e^{-c(t-s)}Q_{ts}^{v,0}
$$ 
and 
$$
\mcQ_\theta^{v,c} \defeq e^{-c\theta}\mcQ_\theta^{v,0}.
$$
Then the operators $c-\mcL^{v}$ and $\partial_t-L^{v}+c$ have inverses of the forms
$$
\big(c-\mcL^{v}\big)^{-1} f = \int_0^\infty \mcQ_\theta^{v,c}f\,\,d\theta = \int_0^1 \mcQ_\theta^{v,c}f \, \,d\theta + \mcQ_1^{v,c}\circ(c-\mcL^{v})^{-1} f
$$
and
$$
(P^{v,c}g)(t)\defeq
\big((\partial_t-L^{v}+c)^{-1}g\big)(t) = \int_{-\infty}^t Q_{ts}^{v,c}g(s)\,ds.
$$
See Theorems {\it\ref{thm:InverseParabolic}} and {\it\ref{thm:SchauderGeneral}} for the proofs and details about the domains of $f$ and $g$.
Then we can write $P^{v,c}$ in terms of $(c-\mcL^{v})^{-1}$.
Indeed, for any suitable functions $f$ on $\bbR\times\bbT$, by setting $g=(c-\mcL^{v})^{-1}f$ and $h=(\partial_t+\partial_x^2)g$ we have
\begin{equation*}
(\partial_t-L^{v}+c)h = \mcL^{v}g+ch=-f+c(g+h),
\end{equation*}
hence
\begin{equation*} \begin{split}
P^{v,c}f &= -h+cP^{v,c}(g+h)   \\
&= -(\partial_t+\partial_x^2)(c-\mcL^{v})^{-1}f   
+ cP^{v,c}(1+\partial_t+\partial_x^2) (c-\mcL^{v})^{-1}f.
\end{split} \end{equation*}
Thus, by setting
$$
K^{v,c}f \defeq -\int_0^1 (\partial_t+\partial_x^2)\mcQ_\theta^{v,c}f\,\,d\theta \eqdef \int_0^1 K_\theta^{v,c}f\,\,d\theta
$$
and
$$
S^{v,c}f \defeq K_1^{v,c}(c-\mcL^{v})^{-1}f + cP^{v,c}(1+\partial_t+\partial_x^2)
(c-\mcL^{v})^{-1}f,
$$
one has the decomposition
\begin{equation} \label{K+R}
P^{v,c}f = K^{v,c}f + S^{v,c}f.
\end{equation}
The letter `$S$' in $S^{v,c}$ is chosen for the `smooth' part. This choice is justified by the regularizing properties of this operator stated in the next statement.

\medskip

\begin{thm}\label{thm:SchauderK+R}
Let $\beta\in(-2,0)\setminus\{-1\}$. Then $K^{v,c}$ is a continuous operator from $\mcC_{\mfs}^{\beta,v}(\bbR\times\bbT)$ into $\mcC_{\mfs}^{\beta+2}(\bbR\times\bbT)$ and $S^{v,c}$ is a continuous operator from $\mcC_{\mfs}^{\beta,v}(\bbR\times\bbT)$ into $\mcC_{\mfs}^{\alpha\wedge(\beta+2)+2-}(\bbR\times\bbT)\defeq\bigcap_{\varepsilon>0}\mcC_{\mfs}^{\alpha\wedge(\beta+2)+2-\varepsilon}(\bbR\times\bbT)$.
\end{thm}

\medskip

\begin{Dem}
The former is obtained from a similar argument to the proof of Theorem {\it \ref{thm:SchauderGeneral}}. The latter is obtained from a combination of Theorem {\it \ref{thm:InverseParabolic}} and Theorem {\it \ref{thm:SchauderGeneral}}. Note that $(c-\mcL^{v})^{-1}$ sends $\mcC_{\mfs}^{\beta,v}(\bbR\times\bbT)$ into $\mcC_{\mfs}^{\beta+4}(\bbR\times\bbT)$ continuously by Theorem {\it \ref{thm:SchauderGeneral}}, and the operator $P^{v,c}$ sends $(1+\partial_t+\partial_x^2)\big(c-\mcL^{v}\big)^{-1}f\in \mcC_{\mfs}^{\beta+2}(\bbR\times\bbT)$ into $\mcC_{\mfs}^{\alpha\wedge(\beta+2)+2-}(\bbR\times\bbT)$ by Theorem {\it \ref{thm:InverseParabolic}}.
\end{Dem}

\medskip

{\bf \textit{We fix a constant $c>c_0\vee C_0$ throughout this paper.}}
We omit the letter `$c$' in $Q^{v},\mcQ^{v},P^v,K^{v}$, and $S^{v}$ unless it needs to be emphasized.

\medskip

\subsection{The regularity structure{\boldmath $.$} \hspace{0.1cm}}
\label{SubsectionRS}

In this section we construct the regularity structure associated with Equation \eqref{EqQgKPZBis}. 
A desired reformulation of \eqref{EqQgKPZBis} is given by the equation
\begin{equation}\label{EqModifiedQgKPZ}
\bsu={\sf L}\big(Q^{v}u_0\big)+{\sf P}^{v,{\sf M}}
\Big\{F(\bsu) \, \Xi+G(\bsu)(\bsD\bsu)^2+c\bsu
+\big\{A(\bsu)-A({\sf L}v)\big\}\bsD^2\bsu\Big\}
\end{equation}
of modelled distributions $\bsu$ defined on an appropriate regularity structure, where $\Xi$ is a placeholder of the noise $\xi$, the operator ${\sf L}$ stands for the canonical lift operator of a spacetime/spatial function to the polynomial part of the regularity structure, the operator ${\sf P}^{v,\sf M}$ is the model dependent integration operator on modelled distributions intertwined to $P^v$ via the reconstruction operator, and the operator $\bsD$ is a natural derivative operator on the space of modelled functions. All these operators are precisely defined in Section {\sf\ref{SubsectionModels}}. In precise, Equation \eqref{EqModifiedQgKPZ} is not a correct form since the parameters $\gamma$ and $\eta$ that define the space $\mcD^{\gamma,\eta}$ of modelled distributions are not considered here. The more precise form \eqref{EqModifiedQgKPZzeta} will be given in Section {\sf\ref{SectionLocalWellPosed}}.

\ssk

We first define a `preparatory' collection of rooted decorated trees as considered in Section {\it 8} of \cite{Hai14} and Section {\it 5} of \cite{BHZ}.
We consider two node types $\Xi$ and ${\bf 1}$, which represent the noise $\xi$ and the constant function $1$, respectively.

\medskip

\begin{defn}\label{defn:pretree}
A finite connected non-planar graph without loops is called a \textbf{tree}.
For any tree $\tau$, its node set, edge set, and root are denoted by $N_\tau$, $E_\tau$, and $\varrho_\tau$, respectively.
The node set $N_\tau$ has a partial order $\le$ defined by $u\le v$ if $u$ lies on the unique path connecting $v$ to the root.
A maximal node with respect to $\le$ is called a leaf.

A tree $\tau$ with maps $\frak{t}:N_\tau\to\{\Xi,{\bf 1}\}$, $\frak{n}:N_\tau\to \bbN^2$, and $\frak{e}:E_\tau\to\bbN^2$ is called a \textbf{preparatory tree}. We say that a preparatory tree $\tau$ is \textbf{good} if $\frak{t}(v)=\Xi$ for any $v\in N_\tau\setminus\{\varrho_\tau\}$ which is a leaf. We denote by $\trees$ the set of all good preparatory trees and by $\trees_\bullet$ the set of all $\tau\in\trees$ such that $\frak{t}(\varrho_\tau)={\bf1}$. Any $\tau\in\trees$ is factorized in the sense of the tree product at the root as
\begin{equation}\label{Eq:FactorizationPreparatory}
\tau=X^{\bf k}\zeta\prod_{i=1}^N\mcI_{{\bf n}_i}(\tau_i)
\end{equation}
with the convention that $\Pi_{i=1}^0=1$, where $N\in\bbN$, ${\bf k},{\bf n}_i\in\bbN^2$, $\zeta\in\{\Xi,{\bf1}\}$, $\tau_i\in\trees\setminus\{X^{\bf k}\}_{{\bf k}\in\bbN^2}$, the notation $\mcI_{{\bf n}_i}(\tau_i)$ denotes a planted tree given by the grafting of the tree $\tau_i$ onto a new root via an edge decorated by ${\bf n}_i$, and $X^{\bf k}$ is a tree with single node decorated by $(\frak{t},\frak{n})=({\bf 1},{\bf k})$. The node type ${\bf 1}$ is identified with $X^{\bf 0}$. Since $\tau$ is good, none of the $\tau_i$ is of the form $X^{\bf k}$. We say that $\tau\in\trees$ with the factorization \eqref{Eq:FactorizationPreparatory} \textbf{strongly conforms to the rule {\rm\bf P}} if all $\tau_i$ strongly conform to the rule {\bf P}, and any of the following conditions hold.
\begin{itemize}
\item[({\bf P}1)]
${\bf n}_i=0$ for any $i$.
\item[({\bf P}2)]
$\zeta={\bf1}$ and ${\bf n}_i=0$ except at most two ${\bf n}_i=(0,1)$.
\item[({\bf P}3)]
$\zeta={\bf1}$ and ${\bf n}_i=0$ except at most one ${\bf n}_i=(0,2)$.
\end{itemize}
We denote by $\basis$ the set of all $\tau\in\trees$ which strongly conform to the rule {\bf P}. The symbol $\basis$ is chosen with the word `basis' in mind.

Finally we say that $\tau\in\trees_\bullet$ \textbf{conforms to the rule {\rm\bf P}} if it is factorized as \eqref{Eq:FactorizationPreparatory} with $\zeta={\bf1}$ and all $\tau_i$ belonging to $\basis$.
We denote by $\basis_\bullet$ the set of all $\tau\in\trees_\bullet$ which conform to the rule {\bf P}.
\end{defn}

\medskip

 We make remarks about the above definition. First, by the definition of $\trees$, we do not consider the trees that contain a leaf of the form $\mcI_{\bf n}(X^{\bf k})$. Such a restriction is usual in the literature of regularity structures as in Definition {\it5.7} of \cite{Hai14}. For ${\bf k}={\bf 0}$ this reflects that the operator $K^{v}$ sends any constant into zero. Indeed,
$$
K^{v}(1) = -\int_0^1 \big(\partial_t+\partial_x^2\big)\mcQ_\theta^{v}(1) \, d\theta = - \int_0^1\big(\partial_t+\partial_x^2\big) e^{-c\theta}\,d\theta = 0.
$$
We do not claim that $K^{v}(x^{\bf k}) = 0$ for ${\bf k}\neq{\bf 0}$ in general but it does not matter because we do not use the symbols $\mcI_{\bf n}(X^{\bf k})$ with ${\bf k}\neq{\bf0}$ to solve Equation \eqref{EqModifiedQgKPZ} in the space $\mcD^{\gamma,\eta}$ of modelled distributions with $\gamma<2$ -- see Theorem {\it\ref{ThmMainAnalytic}} below for details.

As a second remark, note that the above definition of the rule {\bf P} is recursive. Equivalently, we can define $\basis$ as the sum of the sets $\{\basis^{(n)}\}_{n=0}^\infty\subset\trees$ defined by $\basis^{(0)}\defeq\emptyset$ and
\begin{align*}
\basis^{(n+1)}\defeq\basis^{(n)}\cup\big\{F\big((\tau_i)_{i=1}^N\big)\, ;\, F\in\frak{P},\ \tau_i\in\overline{\bbB}^{(n)}\setminus\{X^{\bf k}\}_{{\bf k}\in\bbN^2}\big\},
\end{align*}
where $\frak{P}$ denotes the set of all functionals $F\big((\tau_i)_{i=1}^N\big)=X^{\bf k}\zeta\prod_{i=1}^N\mcI_{{\bf n}_i}(\tau_i)$ satisfying any of ({\bf P}1), ({\bf P}2) or ({\bf P}3).

\medskip

For any fixed $\alpha_0\in(0,1)\setminus\bbQ$, we define the \textbf{\textit{homogeneity map}} $|\cdot|:\trees\to\bbR$ by setting
\begin{align*}
&|X^{\bf k}| \defeq |{\bf k}|_{\mfs},\qquad |\Xi| \defeq \alpha_0-2,   \qquad
|\mcI_{\bf n}(\tau)| \defeq |\tau| + 2 - |{\bf n}|_{\mfs},\qquad |\tau\sigma| \defeq |\tau|+|\sigma|.
\end{align*}
We choose irrational $\alpha_0$ to ensure that
$$
\big|\basis\setminus\{X^{\bf k}\}_{{\bf k}\in\bbN^2}\big| \cap\bbZ = \emptyset.
$$

Since the operator $\mcI_{(0,2)}$ does not change the homogeneity, an infinite number of preparatory trees in $\overline{\bbB}$ may have the same homogeneity. 
Hence we have to treat the model space involving infinite linear combinations of preparatory trees. To deal with such infinite sums it will be convenient to introduce a new set of symbols $\mcI_{\bf n}^p$, with $p\in\bbN$ and ${\bf n}\in\textbf{\textsf{N}}^2$.

\medskip

\begin{defn}
A tree $\tau$ with maps $\frak{t}:N_\tau\to\{\Xi,{\bf 1}\}$, $\frak{n}:N_\tau\to \bbN^2$, $\frak{e}:E_\tau\to\bbN^2$ and $\bsp:E_\tau\to\bbN$ is called a \textbf{contracted tree}. We say that a contracted tree $\tau$ is \textbf{good} if $\frak{t}(v)=\Xi$ for any $v\in N_\tau\setminus\{\varrho_\tau\}$ which is a leaf, and we denote by $\Trees$ the set of all good contracted trees. Similarly to Definition \ref{defn:pretree}, we use the notation $\mcI_{\bf n}^p(\tau)$ to denote a planted tree given by the grafting of the contracted tree $\tau$ onto a new root via an edge decorated by $({\bf n},p)$.

We define the map $\pi:\trees\to\Trees$ recursively by the properties 
\begin{equation*} \begin{split}
&\pi(X^{\bf k}\zeta) = X^{\bf k}\zeta\qquad(\zeta\in\{\Xi,{\bf1}\})   \\
&\pi(\tau\sigma) = \pi(\tau)\pi(\sigma)
\end{split} \end{equation*}
and
$$
\pi\big(\mcI_{\bf n}\big((\mcI_{(0,2)})^p(\tau)\big)\big)=\mcI_{\bf n}^p\big(\pi(\tau)\big)
$$
for any planted trees $\tau\in\trees$ which is not of the form $\tau=\mcI_{(0,2)}(\sigma)$, where $(\mcI_{(0,2)})^p$ is the $p$-fold iteration of the operator $\mcI_{(0,2)}$.
Finally we define the two sets of good contracted trees
$$
\Basis\defeq\pi(\basis),\qquad
\Basis_\bullet\defeq\pi(\basis_\bullet).
$$
\end{defn}

\medskip

The map $\pi$ is injective and its left inverse map $\mcE:\Trees\to\trees$ is defined by
\begin{align*}
\mcE(X^{\bf k}\zeta) \defeq X^{\bf k}\zeta,\qquad 
\mcE(\tau\sigma)=\mcE(\tau)\mcE(\sigma),\qquad
\mcE(\mcI_{\bf n}^p(\tau)) \defeq \mcI_{\bf n}\big((\mcI_{(0,2)})^p\mcE(\tau)\big).
\end{align*}
However $\pi$ is not surjective because the contracted tree $\mcI_{\bf n}^p(\mcI_{(0,2)}^q(\Xi))$ is not in the image of $\pi$. We provide alternative recursive constructions of $\Basis$ and $\Basis_\bullet$ that involves an intermediate set $\bbC$.

\medskip

\begin{defn}\label{defn:contractedtrees}
Any $\tau\in\Trees$ is factorized at the root as
\begin{equation}\label{Eq:FactorizationContracted}
\tau=X^{\bf k}\zeta\prod_{i=1}^N\mcI_{{\bf n}_i}^{p_i}(\tau_i),
\end{equation}
where $N\in\bbN$, ${\bf k},{\bf n}_i\in\bbN^2$, $p_i\in\bbN$, $\zeta\in\{\Xi,{\bf 1}\}$, and $\tau_i\in\Trees\setminus\{X^{\bf k}\}_{{\bf k}\in\bbN^2}$.
We say that a contracted tree with the factorization \eqref{Eq:FactorizationContracted} \textbf{strongly conforms to the rule {\rm\bf C}} if all the $\tau_i$ strongly conform to the rule {\bf C} and any of the following conditions hold.
\begin{itemize}
\item[({\bf C}1)]
${\bf n}_i=0$ for any $i$.
\item[({\bf C}2)]
$\zeta={\bf1}$ and ${\bf n}_i=0$ except at most two ${\bf n}_i=(0,1)$.
\item[({\bf C}3)]
$\zeta={\bf1}$, ${\bf n}_i=0$ except at most one ${\bf n}_i=(0,2)$, and either of $N\ge 2$ or ${\bf k}\neq{\bf 0}$ holds.
\end{itemize}
We denote by $\bbC$ the set of all $\tau\in\Trees$ which strongly conform to the rule {\bf C}. Equivalently we can define $\bbC$ as the union of the sets $\{\bbC^{(n)}\}_{n=0}^\infty\subset\Trees$ by setting $\bbC^{(0)}\defeq\emptyset$ and
\begin{align*}
\bbC^{(n+1)}\defeq\bbC^{(n)}\cup\big\{F\big((\tau_i)_{i=1}^N\big)\, ;\, F\in\frak{C},\ \tau_i\in\bbC^{(n)}\setminus\{X^{\bf k}\}_{{\bf k}\in\bbN^2}\big\},
\end{align*}
where $\frak{C}$ denotes the set of all functionals $F\big((\tau_i)_{i=1}^N\big)=X^{\bf k}\zeta\prod_{i=1}^N\mcI_{{\bf n}_i}^{p_i}(\tau_i)$ satisfying any of \emph{({\bf C}1), ({\bf C}2)} or \emph{({\bf C}3)}.
\end{defn}

\medskip

\begin{prop}\label{prop:extensioncontraction}
The sets $\bbB$ and $\bbB_\bullet$ are given by
\begin{align*}
\bbB &= \bbC\cup\Big\{\mcI_{(0,2)}^p(\tau)\,;\,\tau\in\bbC\setminus\{X^{\bf k}\}_{{\bf k}\in\bbN^2},\ p\in\bbN\Big\},   \\
\bbB_\bullet &= \Bigg\{X^{\bf k}\prod_{i=1}^N\mcI_{{\bf n}_i}^{p_i}(\tau_i)\ ;\, 
N\in\bbN,\ {\bf k},{\bf n}_i\in\bbN^2,\ p_i\in\bbN,\ \tau_i\in\bbC\setminus\{X^{\bf k}\}_{{\bf k}\in\bbN^2}
\Bigg\}.
\end{align*}
\end{prop}

\medskip

\begin{Dem}
We prove only the first identity since the proof of the second identity is similar.

Set 
$$
\bbB' \defeq \bbC\cup\Big\{\mcI_{(0,2)}^p(\tau)\,;\,\tau\in\bbC\setminus\{X^{\bf k}\}_{{\bf k}\in\bbN^2},\ p\in\bbN\Big\}.
$$ 
Pick generic $\tau=X^{\bf k}\zeta\prod_{i=1}^N\mcI_{{\bf n}_i}(\tau_i)\in\basis$ with $\tau_i=(\mcI_{(0,2)})^{p_i}(\sigma_i)$ and $\sigma_i\notin\mcI_{(0,2)}(\basis)$ and assume that $\pi(\sigma_i)\in\bbB'$. Then by definition,
$$
\pi(\tau)=X^{\bf k}\zeta\prod_{i=1}^N\mcI_{{\bf n}_i}^{p_i}\big(\pi(\sigma_i)\big).
$$
Since $\sigma_i$ is not of the form $\mcI_{(0,2)}(\eta)$, we actually have $\pi(\sigma_i)\in\bbC$. Since any of ({\bf P}1), ({\bf P}2) or ({\bf P}3) holds at the root of $\tau$, we see that $\pi(\tau)\in\bbC$ or otherwise $\pi(\tau)=\mcI_{(0,2)}^{p_1}\big(\pi(\sigma_1)\big)$. Hence we have $\pi(\tau)\in\bbB'$. Conversely, pick generic $\tau=X^{\bf k}\zeta\prod_{i=1}^N\mcI_{{\bf n}_i}^{p_i}(\tau_i)\in\bbB'$ and assume that $\tau_i=\pi(\sigma_i)$ for some $\sigma_i\in\basis$.
Since $\sigma_i\in\bbC$, it is not of the form $\mcI_{(0,2)}(\eta)$, and then 
$$
\tau=\pi\Bigg(X^{\bf k}\zeta\prod_{I=1}^N\mcI_{{\bf n}_i}\big((\mcI_{(0,2)})^{p_i}(\sigma_i)\big)\Bigg)\in\pi(\basis).
$$
Thus we have $\bbB=\pi(\basis)=\bbB'$.
\end{Dem}

\medskip

We denote by 
$$
\mcZ : \Trees\to\Trees
$$ 
the idempotent map setting the $\bsp$-decoration of each tree to zero and define
$$
\bbF\defeq\mcZ(\bbB).
$$
We write $\tau^{\bsp}$ for a generic contracted tree when we want to emphasize its $\bsp$-decoration.
Recalling that there is no restriction on the choice of $\bsp$ in Definition {\it\ref{defn:contractedtrees}}, we have the parametrization
\begin{equation}\label{Eq:parametrizationofB}
\bbB = \big\{\tau^{\bsp}\, ;\, \tau\in\bbF,\ \bsp:E_\tau\to\bbN\big\}.
\end{equation}
An infinite number of contracted trees in $\bbB$ may have the same homogeneity, but the set $\bbF$ satisfies the locally finite condition as the cases of locally subcritical semilinear SPDEs.

\medskip

\begin{prop}
The set $A \defeq \big\{\tau\in\bbF\, ;\, |\tau|<\gamma\big\}$ is finite for any $\gamma\in\bbR$.
\end{prop}

\medskip

\begin{Dem}
By modifying the recursive definition of $\bbC$, we define the sets $\{\bbC_0^{(n)}\}_{n=0}^\infty\subset\Trees$ by setting $\bbC_0^{(0)}\defeq\emptyset$ and
\begin{align*}
\bbC_0^{(n+1)}\defeq\bbC_0^{(n)}\cup\big\{F\big((\tau_i)_{i=1}^N\big)\, ;\, F\in\frak{C}_0,\ \tau_i\in\bbC_0^{(n)}\setminus\{X^{\bf k}\}_{{\bf k}\in\bbN^2}\big\},
\end{align*}
where $\frak{C}_0$ denotes the set of all functionals of the form $F\big((\tau_i)_{i=1}^N\big)=X^{\bf k}\zeta\prod_{i=1}^N\mcI_{{\bf n}_i}^{0}(\tau_i)$ satisfying any of ({\bf C}1), ({\bf C}2) and ({\bf C}3). Then by setting $\bbC_0\defeq\bigcup_{n=0}^\infty\bbC_0^{(n)}$ we have
$$
\bbF=\bbC_0\cup\big\{\mcI_{(0,2)}^0(\tau)\,;\,\tau\in\bbC_0\setminus\{X^{\bf k}\}_{{\bf k}\in\bbN^2}\big\}.
$$
Hence it is sufficient to show the following three claims for each $n$ by an induction.

\ssk

\begin{itemize}
\item[(1)]
$|\tau|\ge\alpha_0-2$ for any $\tau\in\bbC_0^{(n)}$.
\item[(2)]
The set $\{\tau\in\bbC_0^{(n)}\, ;\, |\tau|<\gamma\}$ is finite for any $\gamma\in\bbR$ and $n\in\bbN$.
\item[(3)]
Once we set $\alpha^{(n)}=\inf\big\{|\tau|\, ;\, \tau\in\bbC_0^{(n)}\setminus\bbC_0^{(n-1)}\big\}$, we have $\alpha^{(n+1)}\ge\alpha^{(n)}+\alpha_0$.
\end{itemize}

\ssk 

Since $\bbC_0^{(1)}=\big\{X^{\bf k}\zeta\, ;\, {\bf k}\in\bbN^2,\ \zeta\in\{\Xi,{\bf1}\}\big\}$, both (1) and (2) hold for $n=1$.
Suppose that (1) and (2) hold for some $\bbC_0^{(n)}$.
Then since $|\mcI_{\bf 0}^0(\tau_i)|\ge\alpha_0>0$ for any $\tau_i\in\bbC_0^{(n)}$, we have $\big|F\big((\tau_i)_{i=1}^N\big)\big|\ge\alpha_0-2$ for any $F\in\frak{C}_0$. Moreover, if $\tau=F\big((\tau_i)_{i=1}^N\big)\in\bbC_0^{(n)}$ satisfies $|\tau|<\gamma$ for some $\gamma\in\bbR$, then the number $N$ cannot exceed some constant $C$ which depends only on $\alpha_0$ and $\gamma$.
Thus (1) and (2) hold for $\bbC_0^{(n+1)}$.

To show (3) we pick $\tau=F\big((\tau_i)_{i=1}^N\big)\in\bbC_0^{(n+1)}\setminus\bbC_0^{(n)}$ such that $|\tau|=\alpha^{(n+1)}$. If $F$ satisfies (C1), $\tau$ should be of the form $\Xi\mcI_{\bf 0}^0(\sigma)$ with $\sigma\in\bbC_0^{(n)}\setminus\bbC_0^{(n-1)}$. If $F$ satisfies (C2), $\tau$ should be of the forms $\mcI_{(0,1)}^0(\Xi)^2\mcI_{\bf 0}^0(\sigma)$ or $\mcI_{(0,1)}^0(\Xi)\mcI_{(0,1)}^0(\sigma)$ with $\sigma\in\bbC_0^{(n)}\setminus\bbC_0^{(n-1)}$.
If $F$ satisfies (C3), $\tau$ should be of the form $\tau=\mcI_{(0,2)}^0(\Xi)\mcI_{\bf 0}^0(\sigma)$ with $\sigma\in\bbC_0^{(n)}\setminus\bbC_0^{(n-1)}$. In any case $|\tau|\ge\alpha_0+\alpha^{(n)}$.
\end{Dem}

\medskip

We use the parametrization \eqref{Eq:parametrizationofB} to define the norms on some linear spaces spanned by some sets of contracted trees.

\medskip

\begin{defn}\label{defn:Banachspacespannedbycontractedtrees}
For any $\tau\in\mcZ(\Trees)$, we denote by ${\sf T}_\tau$ the linear space spanned by 
$$
\mcZ^{-1}(\tau)=\{\tau^{\bsp}\}_{\bsp:E_\tau\to\bbN}.
$$
Picking parameters $m>0$ and $r\in[1,\infty)$ we define $T_\tau^{(m,r)}$ as the completion of ${\sf T}_\tau$ under the norm defined by
$$
\Bigg\| \sum_{\tau^{\bsp}\in\mcZ^{-1}(\tau)} \lambda_{\bsp}\tau^{\bsp} \Bigg\|_{m,r} \defeq 
\Bigg(\sum_{\tau^{\bsp}\in\mcZ^{-1}(\tau)} |\lambda_{\bsp}|^r\,m^{r|\bsp|}\Bigg)^{\frac{1}{r}},
$$
where 
$$
|\bsp| \defeq  \sum_{e\in E_\tau} \bsp(e).
$$
We also define
$$
T^{(m,r)} \defeq \bigoplus_{\tau\in\bbF}T_\tau^{(m,r)}
$$
as the algebraic sum. 
\end{defn}

\medskip

We consider only $r\in\{1,2\}$ in this paper.
We use $\|\cdot\|_{m,1}$ norms to prove the continuity property of a number of operators. In Section {\sf \ref{SectionRenormalization}} we adopt the $\|\cdot\|_{m,2}$ norms in order to exploit the duality between $T^{(m,2)}$ and $T^{(1/m,2)}$ -- see Proposition {\it\ref{prop:dualityTmT1/m}}. As shown in following statement, the difference between $T^{(m,1)}$ and $T^{(m,2)}$ does not matter if we ignore an infinitesimal small change of $m$.

\medskip

\begin{prop}\label{prop:T1T2embedding}
For any $\tau\in\mcZ(\Trees)$ and $m'<m$, one has the continuous embeddings 
$$
T_\tau^{(m,1)}\hookrightarrow T_\tau^{(m,2)}\hookrightarrow T_\tau^{(m',1)}.
$$
\end{prop}

\medskip

\begin{Dem}
Let $\lambda=\sum_{\bsp}\lambda_{\bsp}\tau^{\bsp}\in{\sf T}_\tau$ be a finite sum. If $\lambda\in T_\tau^{(m,1)}$, since $|\lambda_{\bsp}|\,m^{|\bsp|}\le\|\lambda\|_{m,1}<\infty$, we have
$$
\|\lambda\|_{m,2}^2\le \|\lambda\|_{m,1}\sum_{\bsp}|\lambda_{\bsp}|\,m^{|\bsp|}
\le \|\lambda\|_{m,1}^2.
$$
If $\lambda\in T_\tau^{(m,2)}$ and $m'<m$, by Cauchy--Schwarz inequality we have
$$
\|\lambda\|_{m',1}
\le\Bigg(\sum_{\bsp}|\lambda_{\bsp}|^2m^{2|\bsp|}\Bigg)^{\frac12}\Bigg(\sum_{\bsp}\bigg(\frac{m'}{m}\bigg)^{2|\bsp|}\Bigg)^{\frac12}
\lesssim\|\lambda\|_{m,2}.
$$
\end{Dem}

\medskip

We see that the tree product is indeed continuous in the above norm.

\medskip

\begin{prop}\label{prop:treeproductcontinuous}
Let $\tau,\sigma\in\mcZ(\Trees)$.
For any $m>0$, the bilinear extension ${\sf T}_\tau\times{\sf T}_\sigma\to{\sf T}_{\tau\sigma}$ of the tree product is continuously extended into $T_\tau^{(m,1)}\times T_\sigma^{(m,1)}\to T_{\tau\sigma}^{(m,1)}$.
\end{prop}

\medskip

\begin{Dem}
Set $\eta=\tau\sigma$.
For any $\lambda=\sum_{\tau^{\bsp}\in\mcZ^{-1}(\tau)}\lambda_{\bsp}\tau^{\bsp}\in{\sf T}_\tau$ and $\mu=\sum_{\sigma^{\bsq}\in\mcZ^{-1}(\sigma)}\mu_{\bsq}\sigma^{\bsq}\in{\sf T}_\sigma$ we have
\begin{align*}
\|\lambda\mu\|_{m,1}
&=\sum_{\eta^{\bsr}\in\mcZ^{-1}(\eta)}\bigg|\sum_{\tau^{\bsp},\sigma^{\bsq};\tau^{\bsp}\sigma^{\bsq}=\eta^{\bsr}}
\lambda_{\bsp}\mu_{\bsq}\bigg|\,m^{|\bsr|}\\
&\le\sum_{\eta^{\bsr}\in\mcZ^{-1}(\eta)}\sum_{\tau^{\bsp},\sigma^{\bsq};\tau^{\bsp}\sigma^{\bsq}=\eta^{\bsr}}
|\lambda_{\bsp}\mu_{\bsq}|\,m^{|\bsr|}
=\sum_{\tau^{\bsp},\sigma^{\bsq}}|\lambda_{\bsp}\mu_{\bsq}|\,m^{|\bsp|+|\bsq|}\\
&=\bigg(\sum_{\tau^{\bsp}\in\mcZ^{-1}(\tau)}|\lambda_{\bsp}|\,m^{|\bsp|}\bigg)
\bigg(\sum_{\sigma^{\bsq}\in\mcZ^{-1}(\sigma)}|\mu_{\bsq}|\,m^{|\bsq|}\bigg)
=\|\lambda\|_{m,1}\|\mu\|_{m,1}.
\end{align*}
\end{Dem}

\medskip

To construct a regularity structure for Equation \eqref{EqQgKPZ} we equip the linear spaces spanned by the preparatory or contracted trees with some coproducts, similarly to \cite{Hai14} and \cite{BHZ}.
We write
$$
\treesp\defeq\spa(\trees),\qquad
\treesp_\bullet\defeq\spa(\trees_\bullet)
$$
and
$$
\Treesp\defeq\spa(\Trees),\qquad
\Trees_\bullet\defeq\pi(\trees_\bullet),\qquad
\Treesp_\bullet\defeq\spa(\Trees_\bullet).
$$
Moreover we denote by $\treesp\mathbin{\widehat{\otimes}}\treesp_\bullet$ and $\Treesp\mathbin{\widehat{\otimes}}\Treesp_\bullet$ the tensor products of `bigraded spaces', which is defined in Section {\it2.3} of \cite{BHZ} to permit infinite linear combinations. Since we usually consider the truncation of $\treesp\mathbin{\widehat{\otimes}}\treesp_\bullet$ into finite sums, we omit the details of bigraded spaces here.
Infinite sums in bigraded spaces are purely algebraic objects and different from our topological infinite sums in $T^{(m,r)}$.

\medskip

\begin{defn}
We define the linear map $\icprod:\treesp\to\treesp\mathbin{\widehat{\otimes}}\treesp_\bullet$ by the identities
\begin{align}\label{eq:recursiveDelta}
\begin{aligned}
\icprod\zeta &= \zeta\otimes {\bf1}\qquad(\zeta\in\{\Xi,{\bf1}\}),\qquad
\icprod X^{\bf k} = \sum_{{\bf l}\leq {\bf k}} {{\bf k} \choose {\bf l}} X^{{\bf l}}\otimes X^{{\bf k}-{\bf l}},\\
\icprod(\tau\sigma)&=(\icprod\tau)(\icprod\sigma),\qquad
\icprod(\mcI_{\bf n}\tau) =\big(\mcI_{\bf n}\otimes\id\big)\icprod\tau
+\sum_{{\bf k}\in\bbN^2}\frac{X^{\bf k}}{{\bf k}!}\otimes\mcI_{{\bf n}+{\bf k}}\tau.
\end{aligned}
\end{align}
In the last formula, the grafting map $\mcI_{\bf n}$ is linearly extended by
$$
\mcI_{\bf n}(\tau)=\begin{cases}
\mcI_{\bf n}(\tau)&(\tau\in\basis\setminus\{X^{\bf k}\}_{{\bf k}\in\bbN^2}),\\
0&(\tau\in\{X^{\bf k}\}_{{\bf k}\in\bbN^2}).
\end{cases}
$$
Furthermore, we define the linear map
$$
\iCprod\defeq(\pi\otimes\pi)\circ\icprod\circ\mcE:\Treesp\to\Treesp\mathbin{\widehat{\otimes}}\Treesp_\bullet.
$$
\end{defn}

\medskip

We recall from Section {\it3} of \cite{BHZ} some important properties of $\icprod$. First, for the preparatory tree $\tau$ with decorations $\frak{n},\frak{e}$, the expansion of $\icprod\tau_{\frak{e}}^{\frak{n}}$ is explicitly given by
\begin{equation}\label{eq:expliciticprod}
\icprod\tau_{\frak{e}}^{\frak{n}}
=\sum_{\sigma}\sum_{\frak{n}_\sigma,\frak{e}_{\partial\sigma}}\frac1{\frak{e}_{\partial\sigma}!}
{\frak{n}\choose\frak{n}_\sigma}
\sigma_{\frak{e}}^{\frak{n}_\sigma+\pi\frak{e}_{\partial\sigma}}\otimes
(\tau/\sigma)_{\frak{e}+\frak{e}_{\partial\sigma}}^{[\frak{n}-\frak{n}_\sigma]_\sigma},
\end{equation}
where the first sum is over the set of subtrees $\sigma$ of $\tau$ which contains the root $\varrho_\tau$ of $\tau$, the second sum is over the set of functions $\frak{n}_\sigma:N_\sigma\to\bbN^2$ and
$$
\frak{e}_{\partial\sigma} : \partial\sigma\to\bbN^2
$$ 
where
$$
\partial\sigma\defeq\{e=\{u,v\}\in E_\tau\,;\,u\in N_\sigma, v\in N_\tau\setminus N_\sigma\big\}.
$$
The other notations are defined as follows.
\begin{itemize}
\item
$\frak{e}_{\partial\sigma}!\defeq\prod_{e\in\partial\sigma}\frak{e}(e)!$, ${\frak{n}\choose\frak{n}_\sigma}\defeq\prod_{v\in N_\sigma}{{\frak{n}(v)}\choose{\frak{n}_\sigma(v)}}$
\item
$\pi\frak{e}_{\partial\sigma}:N_\sigma\to\bbN^2$ is defined by $\pi\frak{e}_{\partial\sigma}(u)=\sum_{e=\{u,v\}\in\partial\sigma}\frak{e}_{\partial\sigma}(e)$.
\item
$\tau/\sigma$ is a tree obtained by contracting $\sigma$ in $\tau$ into one single node.
\item
$[\frak{n}-\frak{n}_\sigma]_\sigma:N_{\tau/\sigma}\to\bbN^2$ is defined by $[\frak{n}-\frak{n}_\sigma]_\sigma(\varrho_{\tau/\sigma})=\sum_{u\in N_\sigma}\big(\frak{n}(u)-\frak{n}_\sigma(u)\big)$ and $[\frak{n}-\frak{n}_\sigma]_\sigma(v)=\frak{n}(v)$ for $v\in N_{\tau/\sigma}\setminus\{\varrho_{\tau/\sigma}\}$.
\end{itemize}
Another important property of $\icprod$ is the co-associativity
$$
(\icprod\otimes\id)\icprod=(\id\otimes\icprod)\icprod.
$$
Since $\mcE\circ\pi=\id$, the co-associativity is also inherited by $\iCprod$. Moreover since the rule {\bf P} is `complete' in the sense of Definition {\it5.22} of \cite{BHZ}, these operators are stable on the set of strongly conforming or conforming trees. That is we have
\begin{align*}
\icprod:\spa(\basis)\to\spa(\basis)\mathbin{\widehat{\otimes}}\spa(\basis_\bullet),\\
\iCprod:\spa(\Basis)\to\spa(\Basis)\mathbin{\widehat{\otimes}}\spa(\Basis_\bullet).
\end{align*}
In addition to these algebraic properties, the coproduct is continuous in the following sense. For any $\tau\in\mcZ(\Trees)$ we denote by 
$$
\iota_\tau:{\sf T}_\tau\to\Treesp
$$ 
the canonical injection and by 
$$
{\sf Q}_\tau:\Treesp\to{\sf T}_\tau
$$ 
the canonical projection.

\medskip

\begin{prop}\label{prop:coproductcontinuous}
For any $\tau,\sigma\in\bbF$ and $\eta\in\mcZ(\bbB_\bullet)$, we define the linear map
$$
\iCprod_{\tau\to\sigma,\eta}\defeq({\sf Q}_\sigma\otimes{\sf Q}_\eta)\circ\iCprod\circ\iota_\tau:{\sf T}_\tau\to{\sf T}_\sigma\otimes{\sf T}_\eta.
$$
For any $m>0$ this map is continuously extended into
$$
\iCprod^{(m,1)}_{\tau\to\sigma,\eta}:T_\tau^{(m,1)}\to T_\sigma^{(m,1)}\otimes T_\eta^{(m,1)},
$$
where $T_\sigma^{(m,1)}\otimes T_\eta^{(m,1)}$ is defined as a projective tensor product of Banach spaces.
\end{prop}

\medskip

\begin{Dem}
We first consider the explicit action of $\iCprod$. Since $\mcI_{\bf n}(X^{\bf k})=0$, for any planted tree $\mcI_{\bf n}^p(\tau)$ with $\tau\in\bbC$ we have
\begin{align*}
\iCprod\mcI_{\bf n}^p(\tau)
&=(\pi\otimes\pi)\icprod\mcI_{\bf n}\big((\mcI_{(0,2)})^p\mcE(\tau)\big)\\
&=(\pi\otimes\pi)\bigg\{\big(\mcI_{\bf n}(\mcI_{(0,2)})^p\otimes\id\big)\icprod\mcE(\tau)
+\sum_{\bf k}\frac{X^{\bf k}}{{\bf k}!}\otimes\mcI_{{\bf n}+{\bf k}}\big((\mcI_{(0,2)})^p\mcE(\tau)\big)\bigg\}\\
&=(\mcI_{\bf n}^p\otimes\id)\iCprod\tau+\sum_{\bf k}\frac{X^{\bf k}}{{\bf k}!}\otimes\mcI_{{\bf n}+{\bf k}}^p(\tau)
\end{align*}
unless $\iCprod\tau$ does not produce the term of the form $\mcI_{(0,2)}^q(\sigma)\otimes\eta$ -- the contracted tree $\mcI_{\bf n}^p(\mcI_{(0,2)}^q(\sigma))$ does not belong to $\Basis$ by the rule ({\bf C}3). However, even for such cases, if we regard $\mcI_{\bf n}^p$ as the linear operator from $\Trees$ to $\Trees$ which sends planted trees of the form $\mcI_{(0,2)}^q(\sigma)$ into $\mcI_{\bf n}^{p+q+1}(\sigma)$, the above representation of $\iCprod\mcI_{\bf n}^p(\tau)$ still holds.

This investigation implies that we can almost treat the parameter $p$ of an edge $\mcI_{\bf n}^p$ simply as an `edge type' of typed trees in the sense of \cite{BHZ} except the possibility of contractions. The precise meaning is as follows. For any $\tau\in\bbF$, we assume that the (contracted version of) expansion \eqref{eq:expliciticprod} is given by
$$
\iCprod\tau=\sum_{\sigma^{\bsq},\eta}c_{\sigma^{\bsq}\eta}^\tau\, \sigma^{\bsq}\otimes\eta,
$$
where $\sigma^{\bsq}$ and $\eta$ are some contracted trees. Then for any $\bsp:E_\tau\to\bbN$ we have the expansion with the same coefficients
$$
\iCprod\tau^{\bsp}=\sum_{\sigma^{\bsq},\eta}c_{\sigma^{\bsq}\eta}^\tau\, \sigma^{\bsq+[\bsp]_\sigma}\otimes\eta^{\bsp\vert_\eta},
$$
where $\bsp\vert_\eta$ is the restriction of $\bsp$ to the edge set $E_\eta$ of $\eta$. The map $[\bsp]_\sigma:E_\sigma\to\bbN$ is defined for $e\in E_\sigma$ by $[\bsp]_\sigma(e)=\sum_{e'\in E_\tau,\, \pi(e')=e}\bsp(e')$, where $\pi(e')=e$ means that an edge $e'$ of $\tau$ is contracted to an edge $e$ of $\sigma$. Since 
$$
\big\|\sigma^{\bsq+[\bsp]_\sigma}\big\|_{T_\sigma^{(m,1)}}
\|\eta^{\bsp\vert_\eta}\|_{T_\eta^{(m,1)}}= m^{|\bsp|+|\bsq|}
\le m^{|\bsp|} \big(1\vee m^{|E_\tau|}\big),
$$
we have the boundedness of $\iCprod_{\tau\to\sigma,\eta}$ for any $\lambda=\sum_{\bsp}\lambda_{\bsp}\tau^{\bsp}\in {\sf T}_\tau$ as follows.
\begin{align*}
\|\lambda\|_{T_\sigma^{(m,1)}\otimes T_\eta^{(m,1)}}
&\le|c_{\sigma^{\bsq}\eta}^\tau|\sum_{\bsp}|\lambda_{\bsp}| \big\|\sigma^{\bsq+[\bsp]_\sigma}\big\|_{T_\sigma^{(m,1)}} \|\eta^{\bsp\vert_\eta}\|_{T_\eta^{(m,1)}}\\
&\le|c_{\sigma^{\bsq}\eta}^\tau| \, \big(1\vee m^{|E_\tau|}\big) \|\lambda\|_{T_\tau^{(m,1)}}.
\end{align*}
\end{Dem}

\medskip

We complete the construction of the regularity structure. We define $\Basis_+$ as the subset of $\Basis_\bullet$ consisting of contracted trees $X^{\bf k}\prod_{i=1}^N\mcI_{{\bf n}_i}^{p_i}(\tau_i)$ such that $|\tau_i|+2-|{\bf n}_i|_\mfs>0$ for each $i$ and write

$$
{\sf T}^{\bbB}\defeq\spa(\bbB),\qquad
{\sf T}_+^{\bbB}\defeq\spa(\bbB_+).
$$
Denoting by $p_+:\Treesp_\bullet\to\Treesp_+^{\bbB}$ the canonical projection, we can define the two linear maps
$\Cprod:\Treesp^{\bbB}\to\Treesp^{\bbB}\otimes\Treesp_+^{\bbB}$ and $\Cprod_+:\Treesp_+^{\bbB}\to\Treesp_+^{\bbB}\otimes\Treesp_+^{\bbB}$ by
$$
\Cprod \defeq \big(\id\otimes p_+\big) \iCprod \vert_{{\sf T}^{\bbB}},   \qquad
\Cprod_+\defeq \big(p_+\otimes p_+\big) \iCprod \vert_{\Treesp_+^{\bbB}}.
$$
Similarly to \cite{BHZ} we can conclude that ${\sf T}_+^{\bbB}$ is a Hopf algebra with the tree product and the coproduct $\Delta_+$, and ${\sf T}^{\bbB}$ is a right comudule over ${\sf T}_+^{\bbB}$ with the coproduct $\Delta$.
Hence the pair
$$
\mathscr{T} \defeq \big(({\sf T}_+^{\Basis},\Delta_+),({\sf T}^{\Basis},\Delta)\big)
$$
is a \textbf{\textit{concrete regularity structure}} in the sense of \cite{RSGuide}. The set ${\sf G}_+$ of all algebra morphisms $g:{\sf T}_+\to\bbR$ becomes a group with respect to the convolution product 
$$
g*h\defeq(g\otimes h)\Delta_+.
$$
By Proposition {\it\ref{prop:coproductcontinuous}}, these operators are continuous with respect to $\|\cdot\|_{m,1}$ norms.

\medskip

\begin{cor}
We set 
$$
\bbF_+\defeq\mcZ(\bbB_+)
$$ 
and define 
$$
T_+^{(m,1)} \defeq \bigoplus_{\tau\in\bbF_+}T_\tau^{(m,1)}
$$ 
for any $m>0$ as the algebraic sum of Banach spaces. For any $m>0$ the linear map 
$$
\Cprod:\Treesp^{\bbB}\to\Treesp^{\bbB}\otimes\Treesp_+^{\bbB}
$$
is continuously extended into 
$$
\Delta^{(m,1)}:T^{(m,1)}\to T^{(m,1)}\otimes T_+^{(m,1)}.
$$ 
Similarly then map 
$$
\Cprod_+:\Treesp_+^{\bbB}\to\Treesp_+^{\bbB}\otimes\Treesp_+^{\bbB}
$$
is continuously extended into 
$$
\Delta_+^{(m,1)}:T_+^{(m,1)}\to T_+^{(m,1)}\otimes T_+^{(m,1)}.
$$
\end{cor}

\medskip

We sometimes consider some subsets of $\Basis$ closed with respect to the coproduct.

\medskip

\begin{defn}
We say that a subset $\bbS\subset\Trees$ of the form $\bbS=\mcZ^{-1}(\bbS_0)=\{\tau^{\bsp}\}_{\tau\in\bbS_0,\, \bsp:E_\tau\to\bbN}$ for some $\bbS_0\subset\mcZ(\Trees)$ is \textbf{complete}. A complete subset $\bbS\subset\Basis$ is called a \textbf{\textit{sector}} if 
$$
\Delta(\bbS)\subset {\sf T}^{\bbS}\otimes {\sf T}_+^{\Basis},\qquad
{\sf T}^{\bbS}\defeq\spa(\bbS).
$$
The quantity $\min\big\{|\tau|\, ;\, \tau\in\bbS\big\}$ is called the \textbf{\textit{regularity}} of the sector $\bbS$. 
\end{defn}

\medskip

Especially we consider the sectors
\begin{align*}
&\bbX \defeq \{X^{\bf k}\}_{{\bf k}\in\bbN^2},   \\
&\bbB_- \defeq \big\{\tau\in\bbB\,;\,|\tau|<0\big\},   \\
&\bbU_{\bf n}\defeq\bbX\cup\Big\{\mcI_{\bf n}^p(\tau)\,;\,\tau\in\bbC,\ |\tau|<0,\ p\in\bbN\Big\}   \qquad   ({\bf n}\in\{{\bf0},(0,1),(0,2)\})
\end{align*}
in the next section. Moreover for any sector $\bbS$ and $\gamma\in\bbR$ the set 
$$
\bbS_{<\gamma}\defeq\{\tau\in\bbS\,;\,|\tau|<\gamma\}
$$ 
is also a sector. We close this section with a table of notations used throughout this paper.

\medskip

\begin{center}
\begin{longtable}{ll} \hline
Symbol & Meaning \\\hline
&\\[-8pt]
$\trees$ & The set of all {\it good} preparatory trees \\
$\trees_\bullet$ & The set of all $\tau\in\trees$ such that $\frak{t}(\varrho_\tau)={\bf1}$ \\
$\basis$ & The set of all $\tau\in\trees$ which strongly conform to the rule {\rm P}\\
$\basis_\bullet$ & The set of all $\tau\in\trees$ which conform to the rule {\rm P} \\
$\Trees$ & The set of all good contracted trees \\
$\Trees_\bullet$ & The set of all $\tau\in\Trees$ such that $\frak{t}(\varrho_\tau)={\bf1}$ \\
$\pi$ & The contraction map from $\trees$ to $\Trees$ \\
$\mcE$ & The expansion map from $\Trees$ to $\trees$ \\
$\mcZ$ & The map setting the $\bsp$ decoration of each contracted tree to zero \\
$\Basis$, $\Basis_\bullet$ & $\pi(\basis)$, $\pi(\basis_\bullet)$ \\
$\bbC$ & The set of all $\tau\in\Trees$ which strongly conform to the rule {\rm C}\\
$\Basis_+$ & The set of all $X^{\bf k}\prod_{i=1}^N\mcI_{{\bf n}_i}(\tau_i)\in\Basis_\bullet$ such that $|\tau_i|+2-|{\bf n}_i|_\mfs>0$ for each $i$ \\
$\bbB_-$ & $\{\tau\in\bbB\,;\,|\tau|<0\}$\\
$\bbX$ & $\{X^{\bf k}\}_{{\bf k}\in\bbN^2}$\\
$\bbU_{\bf n}$ & $\bbX\cup\big\{\mcI_{\bf n}^p(\tau)\,;\,\tau\in\bbC,\ |\tau|<0,\ p\in\bbN\big\}$\\
$\treesp$, $\treesp_\bullet$, $\Treesp$, $\Treesp_\bullet$ & Linear spaces spanned by $\trees$, $\trees_\bullet$, $\Trees$, $\Trees_\bullet$, respectively\\
$\Treesp^\Basis$, $\Treesp_+^\Basis$ & Linear spaces spanned by $\basis$, $\basis_+$, respectively\\
${\sf T}^{\bbS}$ & The linear space spanned by $\bbS\subset\bbB$\\
${\sf T}_\tau$ & The linear space spanned by $\mcZ^{-1}(\tau)=\{\tau^{\bsp}\}_{\bsp:E_\tau\to\bbN}$ \\
${\sf Q}_\tau$ & The canonical projection ${\sf T}\to{\sf T}_\tau$\\
$T_\tau^{(m,r)}$ & The completion of ${\sf T}_\tau$ under the $\|\cdot\|_{m,r}$ norm\\
$\icprod$, $\iCprod$ & The coproducts $\treesp\to\treesp\mathbin{\widehat{\otimes}}\treesp_\bullet$ and $\Treesp\to\Treesp\mathbin{\widehat{\otimes}}\Treesp_\bullet$ allowing infinite sums\\
$\Delta$, $\Delta_+$ & The truncated coproducts $\Treesp^{\bbB}\to\Treesp^{\bbB}\otimes\Treesp_+^{\bbB}$ and $\Treesp_+^{\bbB}\to\Treesp_+^{\bbB}\otimes\Treesp_+^{\bbB}$ \\
\hline
\end{longtable}
\end{center}

\medskip

\subsection{Models and modelled distributions{\boldmath $.$} \hspace{0.15cm}}
\label{SubsectionModels}

In this section, we recall from \cite{Hai14, Semi} some analytic notions and theorems to formulate Equation \eqref{EqModifiedQgKPZ}.
First, we define the notions of models and modelled distributions and state the reconstruction theorem.
Next, we introduce a class of admissible models and construct the lift of the operator $(\partial_t-L^v+c)^{-1}$ to the space of modelled distributions.

\ssk

For any sector $\bbS\subset\bbB$ we denote by $\bbS_+\subset\Trees$ the minimum complete subset such that the subalgebra ${\sf T}_+^{\bbS}$ of ${\sf T}_+^{\bbB}$ generated by
$$
P(\bbS_+) \defeq \bbX \cup \Big\{\mcI_{\bf n}^p(\tau)\,;\,{\bf n}\in\bbN^2,\ p\in\bbN,\ \tau\in\bbS_+\setminus\bbX,\ |\tau|+2>|{\bf n}|_\mfs\Big\}
$$
satisfies
$$
\Delta{\sf T}^{\bbS}\subset {\sf T}^{\bbS}\otimes {\sf T}_+^{\bbS},\qquad
\Delta_+{\sf T}_+^{\bbS}\subset{\sf T}_+^{\bbS}\otimes{\sf T}_+^{\bbS},\qquad
{\sf A}_+({\sf T}_+^{\bbS})\subset{\sf T}_+^{\bbS},
$$
where ${\sf A}_+$ is the antipode of Hopf algebra $\Treesp_+^{\Basis}$. Then 
$$
\mathscr{T}^{\bbS} \defeq \big(({\sf T}_+^{\bbS},\Delta_+),({\sf T}^{\bbS},\Delta)\big)
$$ 
is a concrete sub-regularity structure of $\mathscr{T}$.

\medskip

\begin{defn}\label{defn:growthfactor}
Let $\bbS$ be a sector. A pair $\sf M=(g,\Pi)$ made up of a family of algebra morphisms
$$
{\sf g}_z:{\sf T}_+^{\bbS}\to\bbR\qquad(z\in\bbR^2)
$$
and a linear map
$$
{\sf\Pi}:{\sf T}^\bbS\to \mcC_{\mfs}^{\alpha_0-2,v}(\bbR^2)
$$
is called \textbf{a model on $\bbS$ of growth factor $m>0$} if, setting
$$
{\sf g}_{z'z}\defeq{\sf g}_{z'}*{\sf g}_z^{-1},  \qquad  {\sf \Pi}_z^{\sf g}\defeq({\sf \Pi}\otimes{\sf g}_z^{-1})\Delta
$$
one has
$$
\|{\sf g}\|_\tau\defeq\sup_{\tau^{\bsp}\in\mcZ^{-1}(\tau)}\sup_{z,z'\in\bbR^2,\,z\neq z'}\frac{\big|{\sf g}_{z'z}(\tau^{\bsp})\big|}{m^{|\bsp|}\|z'-z\|_{\mfs}^{|\tau|}}<\infty, 
$$
for all $\tau\in\mcZ(P(\bbS_+))$ and
$$
\|{\sf\Pi}\|_\sigma \defeq \sup_{\sigma^{\bsq}\in\mcZ^{-1}(\sigma)}\sup_{\theta\in(0,1]}\sup_{z\in\bbR^2} \frac{\big|\mcQ_\theta^{v}\big({\sf\Pi}_z^{\sf g}\sigma^{\bsq}\big)(z)\big|} {m^{|\bsq|} \, \theta^{|\sigma|/4}} < \infty
$$
for all $\sigma\in\mcZ(\bbS)$. For a sector $\bbS$ such that $\mcZ(\bbS)$ and $\mcZ(P(\bbS_+))$ are finite sets we define
$$
\|{\sf M}\|_{\bbS}\defeq\max_{\tau\in \mcZ(P(\bbS_+))}\|{\sf g}\|_{\tau}+\max_{\sigma\in \mcZ(\bbS)}\|{\sf\Pi}\|_\sigma.
$$
For two models ${\sf M}^{(i)}=({\sf g}^{(i)},{\sf \Pi}^{(i)})$, with $i\in\{1,2\}$, we define $\|{\sf g}^{(1)}:{\sf g}^{(2)}\|_\tau$, $\|{\sf\Pi}^{(1)}:{\sf\Pi}^{(2)}\|_\sigma$ by replacing ${\sf g}_{z'z}$ and ${\sf\Pi}_z^{\sf g}$ above with ${\sf g}_{z'z}^{(1)}-{\sf g}_{z'z}^{(2)}$ and $({\sf\Pi}^{(1)})_z^{{\sf g}^{(1)}}-({\sf\Pi}^{(2)})_z^{{\sf g}^{(2)}}$ respectively. Moreover the model $\sf M$ is said to be \textbf{smooth} if $({\sf \Pi}\tau)(\cdot)$ is a continuous function for any $\tau\in\bbS$.
Finally the model $\sf M$ is said to be \textbf{spatially periodic} if
$$
{\sf g}_{(z'+{\bf e}_2)\, (z+{\bf e}_2)}={\sf g}_{z'z},\qquad
\mcQ_\theta^{v}\big({\sf\Pi}_{z+{\bf e}_2}^{\sf g}(\cdot)\big) (z+{\bf e}_2) = \mcQ_\theta^{v}\big({\sf\Pi}_z^{\sf g}(\cdot)\big)(z)
$$
for any $z,z'\in\bbR^2$.
\end{defn}

\medskip

These conditions ensure that the functions ${\sf g}_{z'z}$ and ${\sf \Pi}_z^{\sf g}$ are continuous with respect to the $\|\cdot\|_{m,1}$ norm. Then by Proposition {\it\ref{prop:coproductcontinuous}} the linear operator $\widehat{{\sf g}_{z'z}}$ on ${\sf T}^{\bbS}$ defined for $\tau\in {\sf T}^{\bbS}$ by
$$
\widehat{{\sf g}_{z'z}}(\tau)\defeq(\id\otimes{\sf g})\Delta\tau 
$$
is continuously extended into the linear operator on 
$$
S^{(m,1)}\defeq\bigoplus_{\tau\in\mcZ(\bbS)}T_\tau^{(m,1)}.
$$ 
Therefore the same analytical arguments as in \cite{RSGuide, Semi, Singular} work well for the proofs of the results stated in this section. 

\medskip

\begin{defn} \label{DefnModelledDistributions}
Let $\bbS$ be a sector and let $\sf M=(g,\Pi)$ be a spatially periodic model on $\bbS$ of a growth factor $m>0$.
For any interval $I\subset\bbR$ and $\eta\leq\gamma$, we denote by $\mcD^{\gamma,\eta}(I)=\mcD_m^{\gamma,\eta}(I,\bbS;{\sf g})$ the set of functions 
$$
\bsu:(I\setminus\{0\})\times\bbT\to S_{<\gamma}^{(m,1)}\defeq\bigoplus_{\tau\in\mcZ(\bbS),\,|\tau|<\gamma}T_\tau^{(m,1)}
$$ 
such that
\begin{align*}
&\llparenthesis\, \bsu\, \rrparenthesis_{\mcD^{\gamma,\eta}(I)}
\defeq \max_{\tau\in\mcZ(\bbS),\,|\tau|<\gamma} 
\sup_{z\in (I\setminus\{0\})\times\bbT}\frac{\|{\sf Q}_\tau\bsu(z)\|_{m,1}}{\big(|t|^{1/2}\wedge1\big)^{(\eta-|\tau|)\wedge0}} < \infty,   \\
&\|\bsu\|_{\mcD^{\gamma,\eta}(I)} \defeq \max_{\tau\in\mcZ(\bbS),\,|\tau|<\gamma} 
\sup_{\substack{z,z'\in (I\setminus\{0\})\times\bbT \\ \|z'-z\|_\mfs\le|t|^{1/2}\wedge|t'|^{1/2}\wedge1}}
\frac{\big\|{\sf Q}_\tau\big\{\bsu(z')-\widehat{{\sf g}_{z'z}}\bsu(z)\big\}\big\|_{m,1}}{\big(|t|^{1/2}\wedge|t'|^{1/2}\wedge1\big)^{\eta-\gamma}\|z'-z\|_\mfs^{\gamma-|\tau|}} < \infty,
\end{align*}
where $t$ and $t'$ represent the first variables of $z$ and $z'$ respectively.
Equipped with the norm $|\!|\!| \bsu |\!|\!|_{\mcD^{\gamma,\eta}(I)} \defeq \llparenthesis\, \bsu\, \rrparenthesis_{\mcD^{\gamma,\eta}(I)} + \|\bsu\|_{\mcD^{\gamma,\eta}(I)}$, the space $\mcD_m^{\gamma,\eta}(I,\bbS;{\sf g})$ is a Banach space. 
In addition, for two models ${\sf M}^{(i)}=({\sf g}^{(i)},{\sf \Pi}^{(i)})$ and $\bsu^{(i)}\in\mcD_m^{\gamma,\eta}(I,\bbS;{\sf g}^{(i)})$ with $i\in\{1,2\}$, we define
$$
|\!|\!| \bsu^{(1)}:\bsu^{(2)} |\!|\!|_{\mcD^{\gamma,\eta}(I)} 
\defeq \llparenthesis\, \bsu^{(1)}:\bsu^{(2)}\, \rrparenthesis_{\mcD^{\gamma,\eta}(I)} + \|\bsu^{(1)}:\bsu^{(2)}\|_{\mcD^{\gamma,\eta}(I)},
$$ 
by replacing ${\sf Q}_\tau\bsu(z)$ and ${\sf Q}_\tau\{\bsu(z')-\widehat{{\sf g}_{z'z}}\bsu(z)\}$ above with ${\sf Q}_\tau\{\bsu^{(1)}(z)-\bsu^{(2)}(z)\}$ and ${\sf Q}_\tau\{\bsu^{(1)}(z')-\widehat{{\sf g}_{z'z}^{(1)}}\bsu^{(1)}(z)\}-{\sf Q}_\tau\{\bsu^{(2)}(z')-\widehat{{\sf g}_{z'z}^{(2)}}\bsu^{(2)}(z)\}$ respectively.
\end{defn}

\medskip

We sometimes omit $m,\bbS$, and ${\sf g}$ from $\mcD_m^{\gamma,\eta}(I,\bbS;{\sf g})$ to lighten notations. On the other hand, the dependence on $\gamma,\eta$, and $I$ is more important. We can prove the following reconstruction theorem similarly to Theorem {\it20} of \cite{RSGuide} and Theorem {\it4.1} of \cite{Semi} -- see Proposition {\it3.8} and Corollary {\it3.9} of \cite{Singular} for details. We define ${\sf Q}_{<\gamma}$ as the canonical projection from $T^{(m,1)}$ to the subspace $\bigoplus_{\tau\in\bbF,\,|\tau|<\gamma}T_\tau^{(m,1)}$.

\medskip

\begin{thm} \label{ThmReconstruction}
Let $\bbS$ be a sector of regularity $\beta\in(-2,0]$ and let $\sf M=(g,\Pi)$ be a spatially periodic model on $\bbS$ of a growth factor $m>0$. If $\gamma\in(0,\beta+2)$ and $\eta\in(\gamma-2,\gamma]$, there exists a unique continuous linear operator
$$
{\sf R}^{\sf M} : \mcD_m^{\gamma,\eta}(\bbR,\bbS;{\sf g})\to\mcC_{\mfs}^{\eta\wedge\beta,v}(\bbR\times\bbT)
$$
such that the bound
$$
\Big| \mcQ_\theta^v\big({\sf R}^{\sf M}\bsu - {\sf \Pi}_z^{\sf g} \bsu(z)\big)(z) \Big| \lesssim \big(\vert t\vert^{1/2}\wedge1\big)^{\eta-\gamma} \, \theta^{\gamma/4}\|{\sf M}\|_{\bbS_{<\gamma}}\|\bsu\|_{\mcD^{\gamma,\eta}(\bbR)}
$$
holds for any $\bsu\in\mcD_m^{\gamma,\eta}(\bbR,\bbS;{\sf g})$ and $z=(t,x)\in(\bbR\setminus\{0\})\times\bbT$. 
Moreover, the mapping $({\sf M},\bsu)\mapsto{\sf R}^{\sf M}\bsu$ is locally Lipschitz continuous in terms of the metrics defined above.
Finally, if $\sf M$ is smooth then ${\sf R}^{\sf M}\bsu$ is realized as a measurable function and
$$
{\sf R}^{\sf M}\bsu(z)=\big({\sf\Pi}_z^{\sf g}\bsu(z)\big)(z)
$$
for any $z\in(\bbR\setminus\{0\})\times\bbT$. 
\end{thm}

\medskip

We define some operations on modelled distributions. Recall from Proposition {\it\ref{prop:treeproductcontinuous}} that the tree product is continuous with respect to the $\|\cdot\|_{m,1}$ norm. Therefore the following statements follow from the same proofs as Propositions \textit{6.12} and \textit{6.13} of \cite{Hai14}.

\medskip

\begin{prop}\label{prop:productMD}
Let $\bbS$ be a sector and let $\sf M=(g,\Pi)$ be a spatially periodic model on $\bbS$ of a growth factor $m>0$. Let $I\subset\bbR$ be an interval.
\begin{itemize}
\setlength{\itemsep}{10pt}
\item[(1)]
Let $\bbS_1,\bbS_2\subset\bbS$ be sectors of regularities $\alpha_1$ and $\alpha_2$ respectively, and such that the tree product $\bbS_1\times \bbS_2\to \bbS$ is defined. Then for any $\bsu_i\in\mcD_m^{\gamma_i,\eta_i}(I,\bbS_i;{\sf g})$ with $i\in\{1,2\}$ and $\eta_i\le\gamma_i$, we have
$$
{\sf Q}_{<\gamma}(\bsu_1\cdot\bsu_2)\in\mcD_m^{\gamma,\eta}(I,\bbS;{\sf g})
$$
with $\gamma=(\gamma_1+\alpha_2)\wedge(\gamma_2+\alpha_1)$ and $\eta=(\eta_1+\alpha_2)\wedge(\eta_2+\alpha_1)\wedge(\eta_1+\eta_2)$. Moreover, the mapping $(\bsu_1,\bsu_2)\mapsto{\sf Q}_{<\gamma}(\bsu_1\cdot\bsu_2)$ is locally Lipschitz continuous.
\item[(2)]
Let $\bbV\subset\bbS$ be a sector containing $\bf1$ and that $\beta=\min\{|\tau|\,;\,\tau\in\bbV\setminus\{{\bf1}\}\}>0$. (Such a sector is called \emph{function-like} in the sense of \cite{Hai14}.) Assume that the product of $n$ trees $\bbV^n\to\bbS$ is defined for any $n\ge1$. For any $\bsv\in\mcD_m^{\gamma,\eta}(I,\bbV;{\sf g})$ with $0\le\eta\le\gamma$ and a function $h\in C^\kappa(\bbR)$ with $\kappa\ge\max\{\gamma/\beta,1\}+1$, we define
$$
H(\bsv) \defeq {\sf Q}_{<\gamma}\Bigg(\sum_{n=0}^\infty\frac{h^{(n)}(v_0)}{n!} \, \big(\bsv-v_0{\bf1}\big)^n\Bigg),
$$
where $v_0$ denotes the ${\bf1}$-component of $\bsv$. Then $H(\bsv)\in\mcD_m^{\gamma,\eta}(I,\bbS;{\sf g})$, and the mapping $\bsv\mapsto H(\bsv)$ is locally Lipschitz continuous.
\end{itemize}
\end{prop}

\medskip

Hereafter, we consider particular sectors $\bbB_-$ and $\bbU_{\bf n}$ for ${\bf n}\in\{{\bf0},(0,1),(0,2)\}$.
We define the linear map $\mcI_{\bf n}:{\sf T}^{\bbB_-\cup\{{\bf1}\}}\to{\sf T}^{\bbU_{\bf n}}$ (called an \emph{abstract integration} in the sense of \cite{Hai14}) by
\begin{align*}
\mcI_{\bf n}(\tau)=\begin{cases}
\mcI_{\bf n}^0(\tau)&\textrm{ if } \tau\in\bbB_-\ \text{is not of the form $\mcI_{(0,2)}^p(\sigma)$},   \\
\mcI_{\bf n}^{p+1}(\sigma)&\textrm{ if } \text{$\tau=\mcI_{(0,2)}^p(\sigma)$ and $\sigma$ is not of the form $\mcI_{(0,2)}^q(\eta)$},   \\
0&\textrm{ if } \tau={\bf1}.
\end{cases}
\end{align*}
Since this definition ensures that
$$
\|\mcI_{\bf n}(\tau)\|_{m,1}\le(1\vee m)\|\tau\|_{m,1}
$$
for any $\tau\in\bbB_-\cup\{{\bf1}\}$, the map $\mcI_{\bf n}$ is continuously extended into the map from $\bigoplus_{\tau\in\mcZ(\bbB_-)\cup\{{\bf1}\}}T_\tau^{(m,1)}$ to $\bigoplus_{\tau\in\mcZ(\bbU_{\bf n})}T_\tau^{(m,1)}$.

\ssk

The notion of admissible model in the present setting is defined as follows. Recall from \eqref{K+R} the decomposition 
$$
P^v = K^{v}+S^{v}
$$
of $P^v=(\partial_t-L^v+c)^{-1}$.

\medskip

\begin{defn}\label{defn:smoothadmissiblemodels}
We say that a smooth model $\sf M=(g,\Pi)$ on $\bbB_-\cup\bbU_{\bf0}\cup\bbU_{(0,1)}$ is \textbf{admissible} if one has
\begin{align*}
&{\sf g}_z(X^{\bf k})=z^{\bf k}, \qquad
({\sf\Pi}X^{\bf k})(z)=z^{\bf k},\\
&({\sf \Pi}X^{\bf k}\tau)(z)=z^{\bf k}({\sf\Pi}\tau)(z)\qquad(\tau\in\bbB_-\ \text{such that}\ X^{\bf k}\tau\in\bbB_-),\\
&{\sf \Pi}(\mcI_{\bf n}\tau) = \partial^{\bf n}K^{v}({\sf \Pi}\tau)\qquad(\tau\in\bbB_-,\ {\bf n}\in\{{\bf0},(0,1),(0,2)\}).
\end{align*}
We denote by $\Models^m$ the set of all smooth admissible models of growth factor $m>0$ and define $\mModels$ as the closure of $\Models^m$ with respect to the metric of the space of models -- see Definition \ref{defn:growthfactor}. 
We simply write $\|{\sf M}\|$ instead of $\|{\sf M}\|_{\bbB_-\cup\bbU_{\bf0}\cup\bbU_{(0,1)}}$ for each ${\sf M}\in\mModels$.
\end{defn}

\medskip

Any smooth admissible model is determined by the ${\sf\Pi}$ part. Indeed, the ${\sf g}$ part is given by
$$
{\sf g}_z^{-1}\big(\mcI_{\bf n}\tau\big) = -\sum_{|{\bf k}|_{\mfs}<|\tau|+2-|{\bf n}|_{\mfs}} \frac{(-z)^{\bf k}}{{\bf k}!}\, \Big(\big(\partial_z^{{\bf n}+{\bf k}}K^{v} \big)\big({\sf \Pi}_z^{\sf g}\tau\big)\Big)(z)
$$
for any $\tau\in\bbB_-\cup\{{\bf 1}\}$ and ${\bf n}\in\bbN^2$ such that $|\tau|+2>|{\bf n}|_\mfs$.
Then ${\sf\Pi}_z^{\sf g}$ satisfies the identity
$$
{\sf\Pi}_z^{\sf g}\mcI_{\bf n}\tau
=(\partial^{\bf n}K^v)({\sf\Pi}_z^{\sf g}\tau)-\sum_{|{\bf k}|_\mfs<|\tau|+2-|{\bf n}|_\mfs}\frac{(\cdot-z)^{\bf k}}{{\bf k}!} \, \big(\partial^{{\bf n}+{\bf k}}K^v\big)({\sf\Pi}_z^{\sf g}\tau)(z)
$$
for any $\tau\in\bbB_-\cup\{{\bf 1}\}$ -- see e.g. Proposition {\it 15} of \cite{RSGuide}. The proof of the multi-level Schauder estimates can be done along the same lines as in Hairer's original statement, Theorem \textit{5.12} of \cite{Hai14}. But we need a slight modification because the kernel $K_\theta(z,w)$ is only twice differentiable with respect to the first variable rather than smooth. See Theorem \textit{5.12} of \cite{Semi} or Lemma {\it4.2} and Corollary {\it4.6} of \cite{Singular} for the proofs of Lemma {\it\ref{LemPointwiseEstimate}} and Theorem {\it\ref{thm:multilevelSchauderK}} below. 

\medskip

\begin{lem}   \label{LemPointwiseEstimate}
Let ${\sf M=(g,\Pi)}\in\mModels$.
For any $\sigma^{\bsp}\in\bbB_-$, $z\in\bbR\times\bbT$, $\theta\in(0,1]$, and ${\bf k}\in\bbN^2$ such that $|{\bf k}|_{\mfs}\le2$, one has
\begin{align*}
\Big| \big(\partial_z^{\bf k}K_\theta^{v}\big) \big({\sf\Pi}_z^{\sf g}\sigma^{\bsp}\big)(z) \Big|
\lesssim m^{|\boldsymbol{p}|} \, \theta^{(|\sigma|-|{\bf k}|_{\mfs}-2)/4}.
\end{align*}
Therefore, if $|{\bf k}|_{\mfs}<|\sigma|+2$ then the integral
$$
\big(\partial_z^{\bf k}K^{v}\big) \big({\sf\Pi}_z^{\sf g}\sigma^{\bsp}\big)(z) = \int_0^1 \big(\partial_z^{\bf k}K_\theta^{v}\big) \big({\sf\Pi}_z^{\sf g}\sigma^{\bsp}\big)(z) \, \,d\theta
$$
converges.
Similarly, for any $\bsu\in\mcD_m^{\gamma,\eta}(\bbR,\bbB_-\cup\{{\bf1}\};{\sf g})$ with $\gamma\in(0,\alpha_0)$ and $\eta\in(\gamma-2,\gamma]$, $z\in(\bbR\setminus\{0\})\times\bbT$, $\theta\in(0,1]$, and ${\bf k}\in\bbN^2$ such that $|{\bf k}|_\mfs\le2$, one has
\begin{align*}
\Big| \big(\partial_z^{\bf k}K_\theta^{v}\big) \big({\sf R}^{\sf M}\bsu - {\sf \Pi}_z^{\sf g} \bsu(z)\big)(z) \Big|
\lesssim \big(\vert t\vert^{1/2}\wedge1\big)^{\eta-\gamma} \, \theta^{(\gamma-|{\bf k}|_\mfs-2)/4}\|{\sf M}\|\|\bsu\|_{\mcD^{\gamma,\eta}(\bbR)},
\end{align*}
hence the integral
$$
\big(\partial_z^{\bf k}K^{v}\big) \big({\sf R}^{\sf M}\bsu - {\sf \Pi}_z^{\sf g} \bsu(z)\big)(z) = \int_0^1 \big(\partial_z^{\bf k}K_\theta^{v}\big) \big({\sf R}^{\sf M}\bsu - {\sf \Pi}_z^{\sf g} \bsu(z)\big)(z) \, \,d\theta
$$
converges.
\end{lem}

\medskip

Recall that $\alpha\in(0,\alpha_0)$ is the regularity of the initial value $u_0$ -- see Section {\sf \ref{SubsectionFunctionSpaces}}.

\medskip

\begin{thm}\label{thm:multilevelSchauderK}
Let ${\sf M=(g,\Pi)}\in\mModels$.
For any $z\in\bbR\times\bbT$, we define the continuous linear map $\mcJ^v(z):{\sf T}^{\bbB_-\cup\{{\bf1}\}}\to\bbR$ by
\begin{align*}
\mcJ^{v}(z)\tau &\defeq \sum_{|{\bf k}|_{\mathfrak{s}}<|\tau|+2} (\partial^{\bf k}K^{v})({\sf\Pi}_z^{\sf g}\tau)(z)\,\frac{X^{\bf k}}{{\bf k}!},\qquad\tau\in\bbB_-\cup{\bf1}.
\end{align*}
Moreover, for any $\bsu\in\mcD_m^{\gamma,\eta}(\bbR,\bbB_-\cup\{{\bf1}\};{\sf g})$ with $\gamma\in(0,\alpha_0)$ and $\eta\in(\gamma-2,\gamma]$, we set
\begin{align*}
\big(\mcN^{v}\bsu\big)(z) &\defeq \sum_{|{\bf k}|_{\mathfrak{s}}<\gamma+2} (\partial^{\bf k}K^{v})\big({\sf R}^{\sf M}\bsu-{\sf\Pi}_z^{\sf g}\bsu(z)\big)(z) \, \frac{X^{\bf k}}{{\bf k}!},
\end{align*}
and for $\rho\le\gamma+2$, we set
\begin{align*}
\big({\sf K}_\rho^{v,{\sf M}}\bsu\big)(z) &\defeq 
{\sf Q}_{<\rho}\Big\{\big(\mcI_{\bf0}+\mcJ^{v}(z)\big)\bsu(z) + \big(\mcN^{v}\bsu\big)(z)\Big\}.
\end{align*}
If $\gamma\in(0,\alpha)$ and $\eta\in(\gamma-2,\gamma]$, the map ${\sf K}_\rho^{v,{\sf M}}$ sends continuously $\mcD_m^{\gamma,\eta}(\bbR,\bbB\cup\{{\bf1}\};{\sf g})$ into $\mcD_m^{\rho,(\eta+2)\wedge\alpha_0}(\bbR,\bbU_{\bf0};{\sf g})$ and the mapping $({\sf M},\bsu)\mapsto{\sf K}_\rho^{v,{\sf M}}\bsu$ is locally Lipschitz continuous.
Moreover, it holds that 
$$
{\sf R}^{\sf M}\big({\sf K}_\rho^{v,{\sf M}}\bsu\big) = K^v\big({\sf R}^{\sf M}\bsu\big)
$$ 
for any $\rho\in[(\eta+2)\wedge\alpha_0,2)$.
\end{thm}

\medskip

We can also realize the spatial differentiation on the sector $\bbU_{\bf0}\cup\bbU_{(0,1)}$.

\medskip

\begin{prop}
Let ${\sf M=(g,\Pi)}\in\mModels$. For ${\bf n}\in\{{\bf0},(0,1)\}$ we define the linear map 
$$
\bsD:{\sf T}^{\bbU_{\bf n}}\to{\sf T}^{\bbU_{{\bf n}+{\bf e}_2}}
$$ 
by
$$
\bsD X^{\bf k} \defeq k_2 X^{{\bf k}-{\bf e}_2}{\bf 1}_{k_2>0},\qquad \bsD\, \mcI_{\bf n}^p(\tau)\defeq\mcI_{{\bf n}+{\bf e}_2}^p(\tau).
$$
For any interval $I\subset\bbR$, $\gamma>1$, $\eta\le\gamma$, and $\bsu\in\mcD_m^{\gamma,\eta}(I,\bbU_{\bf n};{\sf g})$, the $\bbU_{{\bf n}+{\bf e}_2}$-valued function 
$$
(\bsD\bsu)(z)\defeq\bsD(\bsu(z))
$$ 
belongs to $\mcD_m^{\gamma-1,\eta-1}(I,\bbU_{{\bf n}+{\bf e}_2};{\sf g})$. Moreover the mapping $\bsu\mapsto\bsD\bsu$ is Lipschitz continuous.

If $\sf M$ is smooth then for any each $n\in\{1,2\}$, any $\gamma\in(0,\alpha)$, $\eta\in(\gamma-2,\gamma]$, $\rho\in(n,\gamma+2]$, and any $\bsu\in\mcD_m^{\gamma,\eta}(\bbR,\bbB\cup\{{\bf1}\};{\sf g})$, it holds that 
$$
{\sf R}^{\sf M}\big(\bsD^n{\sf K}_\rho^{v,{\sf M}}\bsu\big) = (\partial^{(0,n)}K^v)\big({\sf R}^{\sf M}\bsu\big).
$$ 
\end{prop}

\medskip

\begin{Dem}
Since $\bsD$ preserves the $\bsp$ decoration of each tree, it is continuous with respect to the $\|\cdot\|_{m,1}$ norm for any $m>0$.
Thus the first statement follows similarly to Propositions \textit{5.28} and \textit{6.15} of \cite{Hai14}.

To show the second statement, we set $\bsu(z)=\sum_{\tau\in\bbB_-}u_\tau(z)\tau+u_{\bf0}{\bf1}$. One has by definition
\begin{align*}
\bsD^n{\sf K}_\rho^{v,M}\bsu(z) = \sum_{\tau\in\bbB_-}u_\tau(z)\mcI_{(0,n)}\tau &+ \sum_{\tau\in\bbB_-}u_\tau(z)\sum_{\substack{|{\bf k}|_\mfs<(|\tau|+2)\wedge\rho \\ {\bf k}\ge n{\bf e}_2}} \partial_z^{\bf k}K^v({\sf\Pi}_z^{\sf g}\tau)(z)\frac{X^{{\bf k}-n{\bf e}_2}}{({\bf k}-n{\bf e}_2)!}   \\
&+\sum_{|{\bf k}|_\mfs<\rho,\,{\bf k}\ge n{\bf e}_2} \partial_z^{\bf k}K^{v}\big({\sf R}^{\sf M}\bsu-{\sf\Pi}_z^{\sf g}\bsu(z)\big)(z) \, \frac{X^{{\bf k}-n{\bf e}_2}}{({\bf k}-n{\bf e}_2)!}.
\end{align*}
By using the smoothness of $\sf M$ and the identities
$$
({\sf\Pi}_z\mcI_{(0,n)}\tau)(z)={\bf1}_{|\mcI_{(0,n)}\tau|\le0}(\partial^{(0,n)}K^v)({\sf\Pi}_z^{\sf g}\tau)(z),
\qquad
({\sf\Pi}_zX^{\bf k})(z)={\bf 1}_{{\bf k}={\bf0}},
$$
we have
\begin{align*}
{\sf R}^{\sf M}\big(\bsD^n{\sf K}_\rho^{v,M}\bsu\big)(z) &= \big({\sf\Pi}_z^{\sf g}\bsD^n{\sf K}_\rho^{v,M}\bsu(z)\big)(z)   \\
&= \sum_{\tau\in\bbB_-,\,|\tau|+2\le n}u_\tau(z)(\partial^{(0,n)}K^v)({\sf\Pi}_z^{\sf g}\tau)(z)   \\
&\qquad+ \sum_{\tau\in\bbB_-,\,|\tau|+2>n}u_\tau(z)(\partial^{(0,n)}K^v)({\sf\Pi}_z^{\sf g}\tau)(z)   \\
&\qquad+(\partial^{(0,n)}K^v)\big({\sf R}^{\sf M}\bsu-{\sf\Pi}_z^{\sf g}\bsu(z)\big)(z)   \\
&= (\partial^{(0,n)}K^v)\big({\sf R}^{\sf M}\bsu\big).
\end{align*}
This concludes the proof.
\end{Dem}

\medskip

We prove the following statements in order to apply Theorem {\it \ref{thm:multilevelSchauderK}} to the space of modelled distributions defined on $(0,T)\times\bbT$ for some $T>0$. We write $\tri\cdot\tri_{\mcD^{\gamma,\eta}(0,T)}$ instead of $\tri\cdot\tri_{\mcD^{\gamma,\eta}((0,T))}$ for the norm of that space. For any $\gamma>0$ and any smooth function $f$ on $(0,\infty)\times\bbT$ we define the $\bbX$-valued function ${\sf L}_{<\gamma}f$ by
$$
\big({\sf L}_{<\gamma}f\big)(z) \defeq \sum_{\vert {\bf k}\vert_{\frak{s}}<\gamma} \big(\partial_z^{\bf k}f\big)(z)\,\frac{X^{\bf k}}{{\bf k}!},  
$$
for $z\in (0,\infty)\times\bbT$. We choose the letter $\sf L$ for the word `lift'.

\ssk

\begin{prop}\label{prop:extensionMD}
Let ${\sf M=(g,\Pi)}\in\mModels$.
We fix a smooth non-increasing function $\chi:(0,\infty)\to[0,1]$ such that
$$
\chi(t)=\begin{cases}
1&(0<t\le1),\\
0&(t\ge2),
\end{cases}
$$
and set $\chi_T(t,x)\defeq\chi(t/T)$ for $T>0$.
Let $\gamma>0$ and $\eta\le\gamma$, and let $\bbS$ be a sector of a regularity $\beta\le0$ such that the tree product $X^{\bf k}\times\bbS\to\bbS$ is defined.
Then for any $\bsu\in\mcD_m^{\gamma,\eta}((0,2T),\bbS;{\sf g})$, the function defined by
\begin{align*}
\widetilde{\bsu}_T(z)\defeq
\begin{cases}
0&(t\le0,\ t\ge 2T),\\
{\sf Q}_{<\gamma}\big({\sf L}_{<\gamma-\beta}(\chi_T)(z)\cdot\bsu(z)\big)&(0<t<2T)
\end{cases}
\end{align*}
is an element of $\mcD_m^{\gamma,\eta\wedge\beta}(\bbR,\bbS;{\sf g})$ and satisfies
$$
|\!|\!| \widetilde{\bsu}_T |\!|\!|_{\mcD^{\gamma,\eta\wedge\beta}(\bbR)}\le C|\!|\!| \bsu |\!|\!|_{\mcD^{\gamma,\eta\wedge\beta}(0,2T)}
$$
for some $T$-independent constant $C>0$.
Moreover, we have $\widetilde{\bsu}_T\vert_{(0,T]\times\bbT}=\bsu\vert_{(0,T]\times\bbT}$.
\end{prop}

\medskip

\begin{Dem}
We can check that $|\!|\!| {\sf L}_{<\gamma-\beta}(\chi_T) |\!|\!|_{\mcD^{\gamma-\beta,0}(0,\infty)}\lesssim1$ by definition, so we have from Proposition {\it\ref{prop:productMD}} (1) that
$$
|\!|\!| \widetilde{\bsu}_T |\!|\!|_{\mcD^{\gamma,\eta\wedge\beta}(0,2T)}\lesssim |\!|\!| \bsu |\!|\!|_{\mcD^{\gamma,\eta}(0,2T)}.
$$
We can extend it into $|\!|\!| \widetilde{\bsu}_T |\!|\!|_{\mcD^{\gamma,\eta\wedge\beta}((0,2T])}\lesssim |\!|\!| \bsu |\!|\!|_{\mcD^{\gamma,\eta}(0,2T)}$ by the uniform continuity.
To show that $\widetilde{\bsu}_T\in\mcD^{\gamma,\eta\wedge\beta}(0,\infty)$, we pick $\tau\in\mcZ(\bbS)$, $z=(t,x)\in[2T,\infty)\times\bbT$ and $z'=(t',x')\in(0,2T)\times\bbT$.
By setting $z''=(2T,x')$ we have
\begin{align*}
\big\|{\sf Q}_\tau\big\{\widetilde{\bsu}_T(z') - &\widehat{{\sf g}_{z'z}}(\widetilde{\bsu}_T(z))\big\}\big\|_{m,1}   \\
&\leq \big\|{\sf Q}_\tau\big\{\widetilde{\bsu}_T(z')-\widehat{{\sf g}_{z'z''}}(\widetilde{\bsu}_T(z''))\big\}\big\|_{m,1}
+\big\|{\sf Q}_\tau\big\{\widehat{{\sf g}_{z'z''}}(\widetilde{\bsu}_T(z'')) - \widehat{{\sf g}_{z'z}}(\widetilde{\bsu}_T(z))\big\}\big\|_{m,1}   \\
&\lesssim|\!|\!| \widetilde{\bsu}_T |\!|\!|_{\mcD^{\gamma,\eta\wedge\beta}((0,2T])}(1\wedge t')^{(\eta\wedge\beta-\gamma)/2}\|z'-z''\|_\mfs^{\gamma-|\tau|}   \\
&\lesssim|\!|\!| \bsu |\!|\!|_{\mcD^{\gamma,\eta}(0,2T)}(1\wedge t')^{(\eta\wedge\beta-\gamma)/2}\|z'-z\|_\mfs^{\gamma-|\tau|}.
\end{align*}
For the case that $z\in(0,2T)\times\bbT$ and $z'\in[2T,\infty)\times\bbT$, by the properties of models we have
\begin{align*}
\big\|{\sf Q}_\tau\big\{\widetilde{\bsu}_T(z')-\widehat{{\sf g}_{z'z}}&(\widetilde{\bsu}_T(z))\big\}\big\|_{m,1}   \\
&= \big\|{\sf Q}_\tau\widehat{{\sf g}_{z'z}}\big\{\widehat{{\sf g}_{zz'}} (\widetilde{\bsu}_T(z')) - \widetilde{\bsu}_T(z)\big\}\big\|_{m,1}   \\
&\le\sum_{\sigma\in\mcZ(\bbS),\, |\tau|\le|\sigma|} \|z'-z\|_\mfs^{|\sigma|-|\tau|}\big\|{\sf Q}_\sigma\big\{\widehat{{\sf g}_{zz'}} (\widetilde{\bsu}_T(z')) - \widetilde{\bsu}_T(z)\big\}\big\|_{m,1}   \\
&\lesssim|\!|\!| \bsu |\!|\!|_{\mcD^{\gamma,\eta}(0,2T)}(1\wedge t)^{(\eta\wedge\beta-\gamma)/2} \sum_{\sigma\in\mcZ(\bbS),\, |\tau|\le|\sigma|} \|z'-z\|_\mfs^{|\sigma|-|\tau|} \|z'-z\|_\mfs^{\gamma-|\tau|}   \\
&\lesssim|\!|\!| \bsu |\!|\!|_{\mcD^{\gamma,\eta}(0,2T)}(1\wedge t)^{(\eta\wedge\beta-\gamma)/2}\|z'-z\|_\mfs^{\gamma-|\tau|}.
\end{align*}
On the other hand, $\widetilde{\bsu}_T\in\mcD^{\gamma,\eta\wedge\beta}(-\infty,0)$ is obvious from the definition.
Since $\|z'-z\|_\mfs\le|t'|^{1/2}\wedge|t|^{1/2}\wedge1$ implies that $t$ and $t'$ have the same sign, we obtain that $\widetilde{\bsu}_T\in\mcD^{\gamma,\eta\wedge\beta}(\bbR)$.

We obtain the last statement because ${\sf L}_{<\gamma-\beta}(\chi_T)(z)={\bf 1}$ for any $z\in (0,T]\times\bbT$.
\end{Dem}

\medskip

Hereafter we fix a function $\chi$ as in Proposition {\it\ref{prop:extensionMD}}.

\medskip

\begin{thm}\label{thm:multilevelSchauder}
Let ${\sf M=(g,\Pi)}\in\Models^m$, $\gamma\in(0,\alpha)$, $\eta\in(\gamma-2,\gamma]$, and $\rho\in[(\eta+2)\wedge\alpha_0,\gamma+2]$.
For any $\bsu\in\mcD_m^{\gamma,\eta}\big((0,2T),\bbB_-\cup\{{\bf1}\};{\sf g}\big)$ we define 
$$
{\sf S}_{\rho,T}^{v,\sf M}\bsu \defeq {\sf L}_{<\rho}\big\{S^{v}({\sf R}^{\sf M}\widetilde{\bsu}_T)\big\}
$$ 
and 
$$
{\sf P}_{\rho,T}^{v,\sf M}\bsu \defeq {\sf K}_\rho^{v,\sf M}\widetilde{\bsu}_T + {\sf S}_{\rho,T}^{v,\sf M}\bsu.
$$ 
Then ${\sf P}_{\rho,T}^{v,\sf M}\bsu\in\mcD_m^{\rho,(\eta+2)\wedge\alpha_0}((0,2T),\bbU_{\bf 0};{\sf g})$ and there exists a quadratic function $C(R)$ of $R>0$ such that, for any $\kappa>0$,
\begin{align*}
\big|\!\big|\!\big|{\sf P}_{\rho,T}^{v, \sf M}\bsu \big|\!\big|\!\big|_{\mcD^{\rho,(\eta+2)\wedge\alpha_0-\kappa}(0,2T)} 
\le C(\|{\sf M}\|) T^{\kappa/2} \, |\!|\!| \bsu |\!|\!|_{\mcD^{\gamma,\eta}(0,2T)}.
\end{align*}
Moreover, we have
$$
\left\{
\begin{aligned}
{\sf R}^{\sf M}\big({\sf P}_{\rho,T}^{v,\sf M}\bsu\big) &= P^v \big({\sf R}^{\sf M}\widetilde{\bsu}_T\big),
&&(\rho<2)\\
{\sf R}^{\sf M}\big(\bsD^n{\sf P}_{\rho,T}^{v,\sf M}\bsu\big) &= (\partial^{(0,n)}P^v) \big({\sf R}^{\sf M}\widetilde{\bsu}_T\big).
&&(n\in\{1,2\},\ \rho\in(n,\gamma+2])
\end{aligned}
\right.
$$ 
\end{thm}

\medskip

\begin{Dem}
We know that ${\sf K}_\rho^{v,{\sf M}}\widetilde{\bsu}_T\in\mcD^{\rho,(\eta+2)\wedge\alpha_0}(\bbU_{\bf 0})$ from Theorem {\it \ref{thm:multilevelSchauderK}} and Proposition {\it \ref{prop:extensionMD}}. 
Since $S^{v}$ sends ${\sf R}^{\sf M}\widetilde{\bsu}_T\in\mcC_{\mfs}^{\eta\wedge(\alpha_0-2),v}(\bbR\times\bbT)$ into $\mcC_{\mfs}^{\gamma+2}(\bbR\times\bbT)$ by Theorem {\it\ref{thm:SchauderK+R}}, we have ${\sf S}_{\rho,T}^{v,{\sf M}}\bsu \in \mcD^{\rho,\rho}(0,2T)$. Hence ${\sf P}_{\rho,T}^{v,\sf M}\bsu\in\mcD^{\rho,(\eta+2)\wedge\alpha_0}(0,2T)$. 

To show the estimate of ${\sf P}_{\rho,T}^{v,\sf M}\bsu$ it is convenient to consider the seminorms
$$
\llparenthesis\, \bsw\, \rrparenthesis_{\mcD^{\gamma,\eta}(0,2T)}'
\defeq \max_{\tau\in\mcZ(\bbS),\,|\tau|<\gamma} 
\sup_{z\in (0,2T)\times\bbT}\frac{\big\|{\sf Q}_\tau\bsw(z)\big\|_{m,1}}{\big(|t|^{1/2}\wedge1\big)^{\eta-|\tau|}}
$$
and 
$$
|\!|\!| \bsw |\!|\!|_{\mcD^{\gamma,\eta}(0,2T)}' \defeq \llparenthesis\, \bsw\, \rrparenthesis_{\mcD^{\gamma,\eta}(0,2T)}' + \|\bsw\|_{\mcD^{\gamma,\eta}(0,2T)},
$$
instead of $\llparenthesis\, \bsw\, \rrparenthesis_{\mcD^{\gamma,\eta}(0,2T)}$ and $|\!|\!| \bsw |\!|\!|_{\mcD^{\gamma,\eta}(0,2T)}$. It is obvious that $|\!|\!| \bsw |\!|\!|_{\mcD^{\gamma,\eta}(0,2T)}\le|\!|\!| \bsw |\!|\!|_{\mcD^{\gamma,\eta}(0,2T)}'$ while the reverse inequality fails. However if
$$
\lim_{t\to0}{\sf Q}_{<\eta} \bsw(t,x)=0,
$$
the reverse inequality $|\!|\!| \bsw |\!|\!|_{\mcD^{\gamma,\eta}(0,2T)}'\lesssim|\!|\!| \bsw |\!|\!|_{\mcD^{\gamma,\eta}(0,2T)}$ holds -- see Lemma {\it6.5} of \cite{Hai14}. The tree $\tau={\bf1}$ is the only tree $\tau$ in $\bbB_-\cup\{{\bf1}\}$ satisfying $|\tau|<(\eta+2)\wedge\alpha_0$. By the smoothness of $\sf M$ we have
$$
{\sf Q}_{\bf1}{\sf P}_{\rho,T}^{v, \sf M}\bsu (z)
=\big({\sf\Pi}_z^{\sf g}{\sf P}_{\rho,T}^{v, \sf M}\bsu (z)\big)(z)
= P^v \big({\sf R}^{\sf M}\widetilde{\bsu}_T\big)(z).
$$
Since ${\sf R}^{\sf M}\widetilde{\bsu}_T(z)=\big({\sf\Pi}_{z}^{\sf g}\widetilde{\bsu}_T(z)\big)(z)=0$ for $z\in(-\infty,0)\times\bbT$, we also have 
$$
P^v \big({\sf R}^{\sf M}\widetilde{\bsu}_T\big)(t,x)=\int_{-\infty}^t\int_{\bbT}Q_{ts}^v(x,y)({\sf R}^{\sf M}\widetilde{\bsu}_T)(s,y)dyds=0
$$
for any $t<0$. Since $P^v \big({\sf R}^{\sf M}\widetilde{\bsu}_T\big)$ is continuous by Theorem {\it\ref{thm:SchauderK+R}} we have
$$
\lim_{t\downarrow0}
{\sf Q}_{\bf1}{\sf P}_{\rho,T}^{v, \sf M}\bsu (t,x)
= \lim_{t\downarrow0} P^v \big({\sf R}^{\sf M}\widetilde{\bsu}_T\big)(t,x)
= 0,
$$
so the norms $|\!|\!|\cdot|\!|\!|$ and $|\!|\!|\cdot|\!|\!|'$ are equivalent for ${\sf P}_{\rho,T}^{v, \sf M}\bsu$ and we have
\begin{align*}
\big|\!\big|\!\big|{\sf P}_{\rho,T}^{v, \sf M}\bsu \big|\!\big|\!\big|_{\mcD^{\rho,(\eta+2)\wedge\alpha_0-\kappa}(0,2T)}
&\lesssim
\big|\!\big|\!\big|{\sf P}_{\rho,T}^{v, \sf M}\bsu \big|\!\big|\!\big|_{\mcD^{\rho,(\eta+2)\wedge\alpha_0-\kappa}(0,2T)}'   \\
&\lesssim
T^{\kappa/2}
\big|\!\big|\!\big|{\sf P}_{\rho,T}^{v, \sf M}\bsu \big|\!\big|\!\big|_{\mcD^{\rho,(\eta+2)\wedge\alpha_0}(0,2T)}'   \\
&\lesssim
T^{\kappa/2} \, |\!|\!| \bsu |\!|\!|_{\mcD^{\gamma,\eta}(0,2T)}.
\end{align*}
\end{Dem}

\vfill \pagebreak

\section{Local well-posedness}
\label{SectionLocalWellPosed}

In this section, we prove the local-in-time well-posedness of the modelled form \eqref{EqModifiedQgKPZ} of Equation \eqref{EqQgKPZBis}. More precisely, for a fixed $\beta\in(2-\alpha_0,2-\alpha_0+\alpha)$, we split the terms inside ${\sf P}^{v,\sf M}$ into two parts and consider the system of equations
\begin{equation}\label{EqModifiedQgKPZzeta}
\left\{
\begin{aligned}
\bsu_\gamma&= {\sf L}_{<\gamma}\big(Q^{v}u_0\big) + {\sf P}_{\gamma,T}^{v, \sf M}(\bsw_1+\bsw_2)
\qquad(\gamma\in\{\beta,\beta+\alpha_0\}),   \\
\bsw_1&= {\sf Q}_{<\beta+\alpha_0-2}\big\{F(\bsu_\beta) \, \Xi + G(\bsu_\beta)(\bsD\bsu_\beta)^2 + c\bsu_\beta \big\},   \\
\bsw_2&={\sf Q}_{<\beta+\alpha_0-2}\big\{\big( A(\bsu_\beta) - A({\sf L}_{<\beta}v)\big) \bsD^2\bsu_{\beta+\alpha_0}\big\}.
\end{aligned}
\right.
\end{equation}
We regard $\bsu_{\beta+\alpha_0}$ as a function of $\bsw_1$ and $\bsw_2$ and solve Equation \eqref{EqModifiedQgKPZzeta} with respect to $(\bsu_\beta,\bsw_1,\bsw_2)$.

\ssk

We emphasize some elementary facts before stating and proving the well-posedness result in Theorem \textit{\ref{ThmMainAnalytic}}. 

\medskip

\begin{lem}\label{Liftofpropagation}
Let ${\sf M}=({\sf g},{\sf \Pi})\in\mModels$.

\begin{enumerate}
\renewcommand{\theenumi}{(\roman{enumi})}
\renewcommand{\labelenumi}{(\roman{enumi})}
\setlength{\itemsep}{5pt}

	\item \label{Liftofpropagation1}
For $Q_t$ either equal to $e^{t\partial_x^2}$ or $Q_t^{v,c}$ with $c>0$, we set here
$$
(Qu_0)(t,x)\defeq (Q_tu_0)(x).
$$
For any $\beta\in(0,2+\alpha)$ and $\eta\le\alpha$ we have
$$
|\!|\!|{\sf L}_{<\beta}(Qu_0)|\!|\!|_{\mcD^{\beta,\eta}(0,T)}\lesssim\|u_0\|_{C^\alpha(\bbT)}.
$$
	\item \label{Liftofpropagation2}
For any $\beta\in(1,2)$ and $\eta\le\alpha$ we have
\begin{align*}
\big|\!\big|\!\big| {\sf L}_{<\beta}\big\{\big(Q^{v,c}-e^{(\cdot)\partial_x^2}\big)(u_0)\big\} \big|\!\big|\!\big|_{\mcD^{\beta,\eta}(0,T)} 
\lesssim T^{(\alpha-\eta)/2} \, \|u_0\|_{C^\alpha(\bbT)}.
\end{align*}

	\item \label{Liftofpropagation3}
Let $v\in V^\alpha[0,T]$ satisfy \eqref{eq:neighborhoodu_0}. For any $\beta\in(1,2)$ and $\eta\le\alpha$,we have
\begin{align*}
\tri {\sf L}_{<\beta}(e^{(\cdot)\partial_x^2}u_0-v) \tri_{\mcD^{\beta,\eta}(0,T)} 
\lesssim (\delta+T^{(\alpha-\eta)/2}) \, \|u_0\|_{C^\alpha(\bbT)}.
\end{align*}
\end{enumerate}
\end{lem}

\medskip

\begin{Dem}
We can show \ref{Liftofpropagation1} directly from the definition of norm by using Propositions {\it\ref{smoothingofQt1}} and {\it\ref{PropAnisotropicTaylor}}. For \ref{Liftofpropagation2}, it is sufficient to show that
\begin{align*}
\|Q^{v,c}u_0-e^{(\cdot)\partial_x^2}u_0\|_{L^\infty((0,T)\times\bbT)}\lesssim T^{\alpha/2}\|u_0\|_{C^\alpha(\bbT)}.
\end{align*}
This inequality follows from \eqref{smoothingofQt1-3-1}. Finally \ref{Liftofpropagation3} is obtained from the assumption \eqref{eq:neighborhoodu_0} and the inequality
$$
|\!|\!|{\sf L}_{<\beta}v|\!|\!|_{\mcD^{\beta,\eta}(0,T)}
\lesssim
\|v\|_{L^\infty((0,T)\times\bbT)}
+T^{(\alpha-\eta)/2}\|v\|_{V^\alpha[0,T]}.
$$
\end{Dem}

\medskip

\begin{thm} \label{ThmMainAnalytic}
Fix $\alpha\in(0,\alpha_0)$, $\beta\in(2-\alpha_0,2-\alpha_0+\alpha)$, and $m>0$.
Let $u_0\in C^\alpha(\bbT)$ and choose $v\in V^\alpha[0,T]$ satisfying \eqref{eq:neighborhoodu_0} for a sufficiently small constant $\delta>0$ depending only on $\|u_0\|_{C^\alpha(\bbT)}$. Then there exists a function $t_0:(0,\infty)\to(0,1]$ such that for any $L>0$ and any ${\sf M}\in\Models^m$ with $\|{\sf M}\|\le L$, Equation \eqref{EqModifiedQgKPZzeta} with $T=t_0(L)$ has a unique solution $(\bsu_\beta,\bsw_1,\bsw_2)$ in the class
\begin{align}\label{spaceuvw}
\mcD_m^{\beta,\alpha}\big((0,2T),\bbU_{\bf 0}\big) \times \mcD_m^{\beta+\alpha_0-2,2\alpha-2}\big((0,2T),\bbB_-\cup\{{\bf1}\}\big) \times \mcD_m^{\beta+\alpha_0-2,\alpha-2}\big((0,2T),\bbB_-\cup\{{\bf1}\}\big),
\end{align}
Moreover the solution map
$$
\mathscr{S}_{T}^{u_0,v}:{\sf M}\mapsto{\sf R}^{\sf M}\big((\widetilde{\bsu_\beta})_T\big)
$$ 
is Lipschitz continuous on the set $\{{\sf M}\in\Models^m\,;\,\|{\sf M}\|\le L\}$ for any $L>0$; as such it has a unique continuous extension to all of ${\sf M}\in\mModels$.
\end{thm}

\medskip

Above the spaces of modelled distributions are all defined with respect to $\sf g$.

\medskip

\begin{Dem}
Within this proof only, we use the shorthand notation $\mcD_{(2T)}^{\gamma,\eta} = \mcD_m^{\gamma,\eta}\big((0,2T),\star\big)$ to lighten notations, whatever space $\star$.
We find a solution by the Picard iteration. Let $\bsw_1^{(0)} = \bsw_2^{(0)}=0$ and set
\begin{equation}\label{Defunvn}
\begin{aligned}
\bsu_\gamma^{(n)}&= {\sf L}_{<\gamma}\big(Q^{v}u_0\big) + {\sf P}_{\gamma,T}^{v,{\sf M}}\big(\bsw_1^{(n)} + \bsw_2^{(n)}\big)  \qquad  (\gamma\in\{\beta,\beta+\alpha_0\}),   \\
\bsw_1^{(n+1)}&= {\sf Q}_{<\beta+\alpha_0-2} \Big\{F(\bsu_\beta^{(n)})\Xi + G(\bsu_\beta^{(n)})(\bsD\bsu_\beta^{(n)})^2 + c\bsu_\beta^{(n)}\Big\},   \\
\bsw_2^{(n+1)}&={\sf Q}_{<\beta+\alpha_0-2} \Big\{\big( A(\bsu_\beta^{(n)}) - A({\sf L}_{<\beta}v)\big) \bsD^2\bsu_{\beta+\alpha_0}^{(n)}\Big\}.
\end{aligned}
\end{equation}
In what follows $C$ means a constant which is independent of $T$, $u_0$ and $(\bsu_\beta^{(n)},\bsw_1^{(n)},\bsw_2^{(n)})$. Also we write $P(r)$ for a polynomial of a variable $r$ whose coefficients are independent of $T$, $u_0$ and $(\bsu_\beta^{(n)},\bsw_1^{(n)},\bsw_2^{(n)})$. These quantities may change from one occurrence to the others. 

By Theorem {\it \ref{thm:multilevelSchauder}} we have for $\gamma\in\{\beta,\beta+\alpha_0\}$
\begin{align}\label{unbound}
\begin{aligned}
|\!|\!| \bsu_\gamma^{(n+1)} |\!|\!|_{\mcD_{(2T)}^{\gamma,\alpha}} &\le \big|\!\big|\!\big| {\sf L}_{<\gamma}\big(Q^{v} u_0\big) \big|\!\big|\!\big|_{\mcD_{(2T)}^{\gamma,\alpha}} + C\Big(T^{\kappa/2} \, |\!|\!| \bsw_1^{(n+1)} |\!|\!|_{\mcD_{(2T)}^{\beta+\alpha_0-2,2\alpha-2}} + |\!|\!| \bsw_2^{(n+1)} |\!|\!|_{\mcD_{(2T)}^{\beta+\alpha_0-2,\alpha-2}}\Big)   \\
&\le C\Big(\|u_0\|_{C^\alpha(\bbT)} + T^{\kappa/2} \, |\!|\!| \bsw_1^{(n+1)} |\!|\!|_{\mcD_{(2T)}^{\beta+\alpha_0-2,2\alpha-2}}+|\!|\!| \bsw_2^{(n+1)} |\!|\!|_{\mcD_{(2T)}^{\beta+\alpha_0-2,\alpha-2}}\Big),
\end{aligned}
\end{align}
where $\kappa=\alpha\wedge(\alpha_0-\alpha)>0$. Next we consider $\bsw_1^{(n+1)}$. Since $\bsu_\beta^{(n)}$ takes values in the sector $\bbU_{\bf0}$, all $F(\bsu_\beta^{(n)}),G(\bsu_\beta^{(n)}),A(\bsu_\beta^{(n)})$ are well-defined elements of $\mcD_{(2T)}^{\beta,\alpha}$. Since $\Xi$ has a homogeneity $\alpha_0-2$,
$$
F(\bsu_\beta^{(n)})\Xi\in\mcD_{(2T)}^{\beta+\alpha_0-2,\alpha+\alpha_0-2}.
$$
Since $\bsD\bsu_\beta^{(n)}\in\mcD_{(2T)}^{\beta-1,\alpha-1}$ takes values in the sector $\bbU_{(0,1)}$ of regularity $\alpha_0-1$
$$
(\bsD\bsu_\beta^{(n)})^2\in\mcD_{(2T)}^{\beta+\alpha_0-2,2\alpha-2},
$$
and thus
$$
G(\bsu_\beta^{(n)})(\bsD\bsu_\beta^{(n)})^2\in\mcD_{(2T)}^{\beta+\alpha_0-2,2\alpha-2};
$$
therefore
\begin{align}\label{vnbound}
\begin{aligned}
&|\!|\!| \bsw_1^{(n+1)} |\!|\!|_{\mcD_{(2T)}^{\beta+\alpha_0-2,2\alpha-2}} \le P\big(|\!|\!|\bsu_\beta^{(n)}|\!|\!|_{\mcD_{(2T)}^{\beta,\alpha}}\big).
\end{aligned}
\end{align}
Finally we consider $\bsw_2^{(n+1)}$. 
Since $\bsD^2\bsu_{\beta+\alpha_0}^{(n)}\in\mcD_{(2T)}^{\beta+\alpha_0-2,\alpha-2}$ takes values in the sector $\bbU_{(0,2)}$ of regularity $\alpha_0-2$, for any $\eta\in(0,\alpha)$ we have
\begin{align*}
&|\!|\!| \bsw_2^{(n+1)} |\!|\!|_{\mcD_{(2T)}^{\beta+\alpha_0-2,\alpha-2}}   \\
&\le C \big|\!\big|\!\big| A(\bsu_\beta^{(n)}) - A({\sf L}_{<\beta}(v)) \big|\!\big|\!\big|_{\mcD_{(2T)}^{\beta,\eta}} \, 
|\!|\!|\bsD^2\bsu_{\beta+\alpha_0}^{(n)}|\!|\!|_{\mcD_{(2T)}^{\beta+\alpha_0-2,\alpha-2}}
\\
&\le C \big|\!\big|\!\big| A(\bsu_\beta^{(n)}) - A({\sf L}_{<\beta}(v)) \big|\!\big|\!\big|_{\mcD_{(2T)}^{\beta,\eta}} \,
\Big(\|u_0\|_{C^\alpha(\bbT)} + T^{\kappa/2} \, |\!|\!| \bsw_1^{(n)} |\!|\!|_{\mcD_{(2T)}^{\beta+\alpha_0-2,2\alpha-2}}+|\!|\!| \bsw_2^{(n)} |\!|\!|_{\mcD_{(2T)}^{\beta+\alpha_0-2,\alpha-2}}\Big).
\end{align*}
To obtain a small factor from $A(\bsu_\beta^{(n)}) - A({\sf L}_{<\beta}(v))$ we decompose it into
\begin{align*}
\Big\{A(\bsu_\beta^{(n)}) - A\big({\sf L}_{<\beta}(Q^{v}u_0)\big) \Big\} &+
\Big\{ A\big({\sf L}_{<\beta}(Q^{v} u_0)\big) - A\big({\sf L}_{<\beta}(e^{(\cdot)\partial_x^2}u_0)\big)\Big\}   \\
&+
\Big\{ A\big({\sf L}_{<\beta}(e^{(\cdot)\partial_x^2} u_0)\big) - A\big({\sf L}_{<\beta}(v)\big)\Big\}.
\end{align*}
For the first term, since $A$ is locally Lipschitz as a mapping from $\mcD_{(2T)}^{\beta,\eta}$ to itself and $|\!|\!|\cdot|\!|\!|_{\mcD_T^{\gamma,\eta}}\le |\!|\!|\cdot|\!|\!|_{\mcD_{(2T)}^{\gamma,\alpha}}$, we have
\begin{align*}
\big|\!\big|\!\big|A(\bsu_\beta^{(n)}) &- A\big({\sf L}_{<\beta}(Q^{v}u_0)\big) \big|\!\big|\!\big|_{\mcD_{(2T)}^{\beta,\eta}}   \\
&\leq P\big(|\!|\!| \bsu_\beta^{(n)} |\!|\!|_{\mcD_{(2T)}^{\beta,\alpha}},\|u_0\|_{C^\alpha(\bbT)}\big)  \big|\!\big|\!\big| \bsu_\beta^{(n)} - {\sf L}_{<\beta}(Q^{v}u_0) \big|\!\big|\!\big|_{\mcD_{(2T)}^{\beta,\eta}}   \\
&\leq P\big(|\!|\!| \bsu_\beta^{(n)} |\!|\!|_{\mcD_{(2T)}^{\beta,\alpha}},\|u_0\|_{C^\alpha(\bbT)}\big) |\!|\!|{\sf P}_{\beta,T}^{v,{\sf M}}(\bsw_1^{(n)}+\bsw_2^{(n)}) |\!|\!|_{\mcD_{(2T)}^{\beta,\eta}}   \\
&\le P\big(|\!|\!| \bsu_\beta^{(n)} |\!|\!|_{\mcD_{(2T)}^{\beta,\alpha}},\|u_0\|_{C^\alpha(\bbT)}\big)T^{(\alpha-\eta)/2} \big(|\!|\!| \bsw_1^{(n)} |\!|\!|_{\mcD_{(2T)}^{\beta+\alpha_0-2,2\alpha-2}} + |\!|\!| \bsw_2^{(n)} |\!|\!|_{\mcD_{(2T)}^{\beta+\alpha_0-2,\alpha-2}}\big).
\end{align*}
For the second and third terms we use Lemma \textit{\ref{Liftofpropagation}} to estimate
\begin{align*}
\big|\!\big|\!\big|A\big({\sf L}_{<\beta}(Q^{v} u_0)\big) - A\big({\sf L}_{<\beta}(e^{(\cdot)\partial_x^2}u_0)\big) \big|\!\big|\!\big|_{\mcD_{(2T)}^{\beta,\eta}}   
&\le P\big(\|u_0\|_{C^\alpha(\bbT)}\big) \big|\!\big|\!\big| {\sf L}_{<\beta}\big(Q^{v}u_0-e^{(\cdot)\partial_x^2}u_0\big)\big|\!\big|\!\big|_{\mcD_{(2T)}^{\beta,\eta}}   \\
&\le P\big(\|u_0\|_{C^\alpha(\bbT)}\big) T^{(\alpha-\eta)/2} \, \|u_0\|_{C^\alpha(\bbT)}
\end{align*}
and get from Assumption \eqref{eq:neighborhoodu_0} that
\begin{align*}
\big|\!\big|\!\big|A\big({\sf L}_{<\beta}(e^{(\cdot)\partial_x^2}u_0)\big) - A\big({\sf L}_{<\beta}(v)\big) \big|\!\big|\!\big|_{\mcD_{(2T)}^{\beta,\eta}} 
&\le P\big(\|u_0\|_{C^\alpha(\bbT)}\big) \big|\!\big|\!\big| {\sf L}_{<\beta}(e^{(\cdot)\partial_x^2}u_0-v)\big|\!\big|\!\big|_{\mcD_{(2T)}^{\beta,\eta}}   \\
&\le P\big(\|u_0\|_{C^\alpha(\bbT)}\big)\big(\delta+T^{(\alpha-\eta)/2}\big) \|u_0\|_{C^\alpha(\bbT)}.
\end{align*}
As a result one has
\begin{align}\label{wnbound}
\begin{aligned}
|\!|\!| \bsw_2^{(n+1)} &|\!|\!|_{\mcD_{(2T)}^{\beta+\alpha_0-2,\alpha-2}}   \\
&\le 
\delta P\big(\|u_0\|_{C^\alpha(\bbT)}\big) 
\Big(\|u_0\|_{C^\alpha(\bbT)} + |\!|\!| \bsw_1^{(n)} |\!|\!|_{\mcD_{(2T)}^{\beta+\alpha_0-2,2\alpha-2}} + |\!|\!| \bsw_2^{(n)} |\!|\!|_{\mcD_{(2T)}^{\beta+\alpha_0-2,\alpha-2}}\Big)   \\
&\quad
+
T^{\rho}P\Big(\|u_0\|_{C^\alpha(\bbT)}, |\!|\!| \bsu_\beta^{(n)} |\!|\!|_{\mcD_{(2T)}^{\beta,\alpha}}, |\!|\!| \bsw_1^{(n)} |\!|\!|_{\mcD_{(2T)}^{\beta+\alpha_0-2,2\alpha-2}}, |\!|\!| \bsw_2^{(n)} |\!|\!|_{\mcD_{(2T)}^{\beta+\alpha_0-2,\alpha-2}}\Big)
\end{aligned}
\end{align}
for a small exponent $\rho>0$.
By \eqref{unbound}, \eqref{vnbound} and \eqref{wnbound}, by choosing some sufficiently small constants $\delta,t_0>0$ we can find some large constants $R_1,R_2,R_3>0$ depending only on $\|u_0\|_{C^\alpha(\bbT)}$ and $L$ and such that
$$
|\!|\!|\bsu_\beta^{(n)}|\!|\!|_{\mcD_{(2T)}^{\beta,\alpha}} \le R_1,\qquad
|\!|\!| \bsw_1^{(n)} |\!|\!|_{\mcD_{(2T)}^{\beta+\alpha_0-2,2\alpha-2}}  \le R_2,\qquad
|\!|\!| \bsw_2^{(n)} |\!|\!|_{\mcD_{(2T)}^{\beta+\alpha_0-2,\alpha-2}} 
\le R_3
$$
for any $n\in\bbN$. Note that $\delta>0$ is chosen to get $\delta P(\|u_0\|_{C^\alpha(\bbT)})\ll1$, so it depends only on $\|u_0\|_{C^\alpha(\bbT)}$.
By the local Lipschitz estimates of the operations in \eqref{Defunvn}: product, composition with smooth functions, differentiation and integration, we have a similar estimate
\begin{align*}
|\!|\!|\bsu_\beta^{(n+1)} &- \bsu_\beta^{(n)}|\!|\!|_{\mcD_{(2T)}^{\beta,\alpha}} + |\!|\!| \bsw_1^{(n+1)}-\bsw_1^{(n)} |\!|\!|_{\mcD_{(2T)}^{\beta+\alpha_0-2,2\alpha-2}} + |\!|\!| \bsw_2^{(n+1)}-\bsw_2^{(n)} |\!|\!|_{\mcD_{(2T)}^{\beta+\alpha_0-2,\alpha-2}}   \\
&\le P(R_1,R_2,R_3) \, T^\rho\, \Big\{|\!|\!|\bsu_\beta^{(n)}-\bsu_\beta^{(n-1)}|\!|\!|_{\mcD_{(2T)}^{\beta,\alpha}} + |\!|\!| \bsw_1^{(n)}-\bsw_1^{(n-1)} |\!|\!|_{\mcD_{(2T)}^{\beta+\alpha_0-2,2\alpha-2}}   \\
&\hspace{186pt}+ |\!|\!| \bsw_2^{(n)}-\bsw_2^{(n-1)} |\!|\!|_{\mcD_{(2T)}^{\beta+\alpha_0-2,2\alpha-2}}
\Big\}
\end{align*}
for a small exponent $\rho>0$. Hence we can choose $T$ smaller such that $(\bsu_\beta^{(n)},\bsw_1^{(n)},\bsw_2^{(n)})$ is a Cauchy sequence. The limit solves Equation \eqref{EqModifiedQgKPZzeta}. Uniqueness also holds because of the local Lipschitz estimates.
\end{Dem}

\medskip

Otto, Sauer, Smith \& Weber \cite{OSSW2} and Linares, Otto \& Tempelmayr \cite{LOT} set up an analytic and algebraic framework to deal with the quasilinear equation \eqref{EqQgKPZLinear} with an additive forcing, that is $f=1$ and $g=0$. They use in particular a greedy index set for their local expansions and prove an a priori bound for the solutions to a renormalized form of their equation driven by a smooth noise. Their result holds in the full sub-critical regime but they do not prove a well-posedness result for their equation. The a priori result entails a compactness statement that ensures the existence of some converging subsequence when the regularizing parameter in the noise is sent to $0$. The analysis of the present section shows that one can run the analysis of the general equation \eqref{Defunvn} within the variant of the usual regularity structure for the generalized (KPZ) equation described in Section {\sf\ref{SectionSetting}}. The present section can also be seen as a simple alternative to the somewhat convoluted approach of Gerencs\'er \& Hairer \cite{GerencserHairer}. The interest of this reformulation of \eqref{EqQgKPZ} will be clear in the next section. The formulation of \cite{GerencserHairer} does not lend itself to an easy formulation of a renormalized equation for \eqref{EqQgKPZ}. At the level of generality of \cite{GerencserHairer} the counterterm in their renormalized equation is a priori a nonlocal functional of the solution. Our main result, Theorem {\it\ref{ThmMainRenormalized}} in Section {\sf\ref{SectionIntro}}, shows that there is a renormalized equation whose counterterm is a local functional of its solution. (Recall there is not a unique renormalized equation.) We prove Theorem {\it\ref{ThmMainRenormalized}} in the next section.


\section{Renormalization matters}
\label{SectionRenormalization}

This section is dedicated to the analysis of the equation satisfied by the solution $u = \mathscr{S}_T^{u_0,v}({\sf M})$ obtained in Theorem \textit{\ref{ThmMainAnalytic}} -- the so-called renormalized equation.  The parameter $\eta$ of $\mcD_m^{\gamma,\eta}$ does not matter in this section, so we work with the shortened version
\begin{equation}\label{EqModifiedQgKPZContracted}
\left\{
\begin{aligned}
\bsu_\gamma&= {\sf L}_{<\gamma}\big(Q^{v}u_0\big) + {\sf P}_{\gamma,T}^{v, \sf M}(\bsw)  \qquad  (\gamma\in\{\beta,\beta+\alpha_0\})   \\
\bsw&= {\sf Q}_{<\beta+\alpha_0-2}\Big\{F(\bsu_\beta)\Xi + G(\bsu_\beta)(\bsD\bsu_\beta)^2 + c\bsu_\beta+\big( A(\bsu_\beta) - A({\sf L}_{<\beta}v)) \bsD^2\bsu_{\beta+\alpha_0}\Big\}
\end{aligned}
\right.
\end{equation}
rather than with \eqref{EqModifiedQgKPZzeta}.

The first systematic treatment of this equation in a semilinear setting was done by Bruned, Chandra, Chevyrev \& Hairer in \cite{BCCH}. They relied on a morphism property satisfied by the coefficients of solutions to semilinear singular SPDEs as modelled distributions, for some multi-pre-Lie structures. A deeper structure on the elements of BHZ regularity structures was unveiled by Bruned \& Manchon in \cite{BrunedManchon} and applied by Bailleul \& Bruned in \cite{NonInvariant} to simplify a lot the analysis of the renormalized equation. This structure is encoded in the $\star$ product introduced in Section {\sf\ref{SubsectionCoherenceStarProduct}}. Its importance in the analysis of Equation \eqref{EqModifiedQgKPZContracted} is emphasized by Proposition \textit{\ref{PropMorphism}}, which provides a basic morphism property -- the counterpart here of the multi-pre-Lie morphism property used in \cite{BCCH}. We introduce in Section {\sf\ref{SubsectionPreparationMaps}} the class of preparation maps -- special operators on the linear space spanned by preparatory trees, and their associated admissible models. A preliminary form of Theorem \textit{\ref{ThmMainRenormalized}} follows from their properties in Proposition \textit{\ref{ThmPreliminaryMain}}. 
We show in Section {\sf\ref{SubsectionRenormalizedEquation}} that working with the preparation map associated with the analogue in our setting of the BHZ character leads to Theorem \textit{\ref{ThmMainRenormalized}}. 


\smallskip

\subsection{Notations{\boldmath $.$} \hspace{0.15cm}}
\label{SectionNotations}

We fix some notations used throughout this section. We first introduce some combinatorial factors. The combinatorial symmetry factor $\overline{S}(\tau)$ of a preparatory tree
\begin{equation*}
\tau=X^{\bf k}\zeta\prod_{i=1}^N\big(\mcI_{{\bf n}_i}(\tau_i)\big)^{m_i}
\end{equation*}
with $({\bf n}_i,\tau_i)\neq({\bf n}_j,\tau_j)$ for any $i\neq j$ is inductively defined by
$$
\overline{S}(\tau)\defeq{\bf k}! \hspace{0.02cm} \prod_{i=1}^n\big(\overline{S}(\tau_i)^{m_i}m_i!\big).
$$
This definition contains the initial case $\overline{S}(X^{\bf k}\zeta)={\bf k}!$ by the convention that $\prod_{i=1}^0=1$. Similarly the factor $S(\tau)$ of a contracted tree
\begin{equation}\label{Eq:FactorizationContractedGrouped}
\tau=X^{\bf k}\zeta\prod_{i=1}^N\big(\mcI_{{\bf n}_i}^{p_i}(\tau_i)\big)^{m_i}
\end{equation}
with $({\bf n}_i,p_i,\tau_i)\neq({\bf n}_j,p_j,\tau_j)$ for any $i\neq j$ is inductively defined by
$$
S(\tau)\defeq{\bf k}! \hspace{0.02cm} \prod_{i=1}^n\big(S(\tau_i)^{m_i}m_i!\big).
$$

\medskip

\begin{prop}
$S(\tau)=\overline{S}(\mcE(\tau))$ for any $\tau\in\bbB$.
\end{prop}

\medskip

\begin{Dem}
The identity holds for the initial case $\tau=X^{\bf k}\zeta$. We consider the contracted tree $\tau$ of the form \eqref{Eq:FactorizationContractedGrouped} with $\tau_i\in\bbC$ and assume that $S(\tau_i)=\overline{S}(\mcE(\tau_i))$ for each $i$. Then we have
$$
\mcE(\tau)=X^{\bf k}\zeta\prod_{i=1}^N\big(\mcI_{{\bf n}_i}(\mcI_{(0,2)})^{p_i}\mcE(\tau_i)\big)^{m_i}
$$
and $({\bf n}_i,(\mcI_{(0,2)})^{p_i}\mcE(\tau_i))\neq({\bf n}_j,(\mcI_{(0,2)})^{p_j}\mcE(\tau_j))$ for any $i\neq j$ by the injectivity of $\mcE$ on $\bbB$. Since $\overline{S}(\mcI_{(0,2)}(\sigma))=\overline{S}(\sigma)$ for any preparatory trees $\sigma$ by definition we have indeed

\begin{align*}
\overline{S}(\mcE(\tau))={\bf k}! \hspace{0.02cm} \prod_{i=1}^N
\big(\overline{S}(\mcE(\tau_i))^{m_i}m_i!\big)
={\bf k}! \hspace{0.02cm} \prod_{i=1}^N
\big(S(\tau_i)^{m_i}m_i!\big)
=S(\tau).
\end{align*}
\end{Dem}

\medskip

For $\tau\in\mcZ(\Trees)$ and $\bsp:E_\tau\to\bbN$ we have 
$$
S(\tau^{\bsp})\le S(\tau)
$$ 
since the $\bsp$ decoration may break the symmetry of $\tau$. We prove an important lemma which will be used in the proof of Theorem {\it \ref{ThmMainRenormalizedSimplified}}. We consider the set $\widetilde{\Trees}$ of all \textbf{planar} trees $\tau$ with node decorations $\frak{t}:N_\tau\to\{\Xi,{\bf1}\}$ and $\frak{n}:N_\tau\to\bbN^2$, and edge decorations $\frak{e}:E_\tau\to\bbN^2$ and $\bsp:E_\tau\to\bbN$. That is, the order of branches of $\tau$ is not commutative unlike the set $\Trees$ of non-planar trees. We denote by 
$$
\mcP:\widetilde{\Trees}\to\Trees
$$ 
the canonical projection and define the {\bf multiplicity factor of $\tau^{\bsp}\in\Trees$} by
$$
M(\tau^{\bsp}) \defeq |\mcP^{-1}(\tau^{\bsp})|.
$$

\medskip

\begin{lem}\label{NonplanarvsPlanar}
For any $\tau\in\mcZ(\Trees)$ and $\bsp:E_\tau\to\bbN$, we have
\begin{equation}\label{SymmetryMultiplicity}
S(\tau)=S(\tau^{\bsp}) \, M(\tau^{\bsp}).
\end{equation}
Consequently, for any $\tau\in\mcZ(\Trees)$, any function $f:\mcZ^{-1}(\tau)\to\bbR$ and any fixed planar tree $\widetilde{\tau}$ such that $\mcP(\widetilde{\tau})=\tau$, we have
$$
\sum_{\tau^{\bsp}\in\mcZ^{-1}(\tau)}\frac{f(\tau^{\bsp})}{S(\tau^{\bsp})}\tau^{\bsp} 
= 
\frac1{S(\tau)}\sum_{\tau^{\bsp}\in\mcZ^{-1}(\tau)}f(\tau^{\bsp})M(\tau^{\bsp})\tau^{\bsp}
=
\frac1{S(\tau)} \, 
\mcP\bigg(\sum_{\bsq:E_{\widetilde{\tau}}\to\bbN}
f\big(\mcP(\widetilde{\tau}^{\bsq})\big)\widetilde{\tau}^{\bsq}\bigg),
$$
if the far left hand side converges in $T_\tau^{(m,1)}$ for some $m>0$.
\end{lem}

\medskip

For example, the contracted tree $\tau^{p,q} \defeq \mcI^p(\Xi)\mcI^q(\Xi)$ is distinguished with $\tau^{q,p}$ if $p\neq q$ as a planar tree but $\tau^{p,q}=\tau^{q,p}$ as non-planar trees. Then we have
\begin{align*}
\sum_{\tau^{p,q}\in\mcZ^{-1}(\tau^{0,0})}\frac{f(\tau^{p,q})}{S(\tau^{p,q})}\tau^{p,q}
=\sum_{p\in\bbN}\frac{f(\tau^{p,p})}{2}\tau^{p,p}+\sum_{p<q\in\bbN}f(\tau^{p,q})\tau^{p,q}
=\frac1{S(\tau^{0,0})}\sum_{(p,q)\in\bbN^2}f(\tau^{p,q})\tau^{p,q}.
\end{align*}

\medskip

\begin{Dem}[of Lemma \ref{NonplanarvsPlanar}]
It is sufficient to show the identity \eqref{SymmetryMultiplicity}. The initial case $\tau=X^{\bf k}\zeta$ is obvious. We consider the contracted tree $\tau\in\mcZ(\Trees)$ of the form 
$$
\tau=X^{\bf k}\zeta\prod_{i=1}^N\big(\mcI_{{\bf n}_i}^0(\tau_i)\big)^{m_i}
$$
with $({\bf n}_i,\tau_i)\neq({\bf n}_j,\tau_j)$ for any $i\neq j$ and the decoration $\bsp:E_\tau\to\bbN$ of the form
$$
\tau^{\bsp}=X^{\bf k}\zeta\prod_{i=1}^N\bigg(\prod_{j=1}^{L_i}\big(\mcI_{{\bf n}_i}^{p_{ij}}(\tau_i^{\bsp_{ij}})\big)^{\ell_{ij}}\bigg),
$$
where $\sum_{j=1}^{L_i}\ell_{ij}=m_i$ and $(p_{ij},\bsp_{ij})\neq(p_{ik},\bsp_{ik})$ for any $j\neq k$. If the identity \eqref{SymmetryMultiplicity} holds for any $\tau_i$, we have
\begin{equation}\label{Eq:SymmetryMultiplicity}
\frac{S(\tau)}{S(\tau^{\bsp})}
=\prod_{i=1}^N
\frac{m_i!}{\ell_{i1}!\cdots\ell_{iL_i}!}
\bigg(\prod_{j=1}^{L_i}M(\tau_i^{\bsp_{ij}})^{\ell_{ij}}\bigg).
\end{equation}
We count the total number of ways of putting the $\bsp$ decoration at the subtree $\big(\mcI_{{\bf n}_i}^0(\tau_i)\big)^{m_i}$ to make it into $\prod_{j=1}^{L_i}\big(\mcI_{{\bf n}_i}^{p_{ij}}(\tau_i^{\bsp_{ij}})\big)^{\ell_{ij}}$. For each $i$, the factor $\frac{m_i!}{\ell_{i1}!\cdots\ell_{iL_i}!}$ represents the total number of ways to assign $\ell_{i1}$ identical pairs $(p_{i1},\bsp_{i1})$, $\ell_{i2}$ identical pairs $(p_{i2},\bsp_{i2})$,..., and $\ell_{iL_i}$ identical pairs $(p_{iL_i},\bsp_{iL_i})$ to the $m_i$ branches of $\big(\mcI_{{\bf n}_i}^0(\tau_i)\big)^{m_i}$. On the other hand, $M(\tau_i^{\bsp_{ij}})$ is the total number of ways of putting the $\bsp_{ij}$ decoration at $\tau_{ij}$. Therefore $\frac{m_i!}{\ell_{i1}!\cdots\ell_{iL_i}!}\prod_{j=1}^{L_i}M(\tau_i^{\bsp_{ij}})^{\ell_{ij}}$ is the total number of planar trees sent to $\prod_{j=1}^{L_i}\big(\mcI_{{\bf n}_i}^{p_{ij}}(\tau_i^{\bsp_{ij}})\big)^{\ell_{ij}}$ by $\Pi$. Since the ways of putting the $\bsp$ decoration at $\big(\mcI_{{\bf n}_i}^0(\tau_i)\big)^{m_i}$ are independent for distinct $i$ and $j$, the right hand side of \eqref{Eq:SymmetryMultiplicity} is equal to $M(\tau^{\bsp})$.
\end{Dem}

\medskip

Next we introduce the inner product on the linear space spanned by preparatory or contracted trees. We define the bilinear map $\langle\cdot,\cdot\rangle:\treesp\times\treesp\to\bbR$ by setting
$$
\langle\tau,\sigma\rangle\defeq \overline{S}(\tau){\bf1}_{\tau=\sigma}
$$
for any $\tau,\sigma\in\trees$ and define the bilinear map $\langle\cdot,\cdot\rangle:\Treesp\times\Treesp\to\bbR$ by 
$$
\langle\tau,\sigma\rangle \defeq \langle\mcE(\tau),\mcE(\sigma)\rangle
$$ 
for any $\tau,\sigma\in\Trees$.

\medskip

\begin{prop}\label{prop:dualityTmT1/m}
For any $\tau\in\mcZ(\Trees)$ and $m>0$, the bilinear map $\langle\cdot,\cdot\rangle:{\sf T}_\tau\times{\sf T}_\tau\to\bbR$ is continuously extended into $\langle\cdot,\cdot\rangle:T_\tau^{(m,2)}\times T_\tau^{(1/m,2)}\to\bbR$.
Furthermore, this map defines an equivalence between the topological dual $(T_\tau^{(m,2)})^*$ of $T_\tau^{(m,2)}$ and $T_\tau^{(1/m,2)}$.
\end{prop}

\medskip

\begin{Dem}
The inequality $|\langle\lambda,\mu\rangle|\le S(\tau)\|\lambda\|_{1/m,2}\|\mu\|_{m,2}$ holds for any $\lambda, \mu\in{\sf T}_\tau$ by the definition of norms and Cauchy--Schwarz inequality.
By this continuity, the injection map
$$
\iota:\lambda\in T_\tau^{(1/m,2)}\mapsto\langle\lambda,\cdot\rangle\in (T_\tau^{(m,2)})^*
$$ 
is continuous. On the other hand the norm of $T_\tau^{(m,2)}$ is induced by another inner product
$$
\Bigg\langle\hspace{-0.2cm}\Bigg\langle
\sum_{\bsp} \lambda_{\bsp} \tau^{\bsp} \,,\, \sum_{\bsp}\mu_{\bsp}\tau^{\bsp}
\Bigg\rangle\hspace{-0.2cm}\Bigg\rangle
\defeq\sum_{\bsp}\lambda_{\bsp} \, \mu_{\bsp} \, m^{2|\bsp|}.
$$
By Riesz representation theorem, for any $f\in(T_\tau^{(m,2)})^*$, there exists $\lambda=\sum_{\bsp}\lambda_{\bsp}\tau^{\bsp}\in T_\tau^{(m,2)}$ such that $f(\mu)=\langle\!\langle\lambda,\mu\rangle\!\rangle$ for any $\mu\in T_\tau^{(m,2)}$.
Setting $\nu_{\bsp}\defeq\lambda_{\bsp}m^{2|\bsp|}/S(\tau^{\bsp})$ and $\nu\defeq\sum_{\bsp}\nu_{\bsp}\tau^{\bsp}$, we have that $\|\nu\|_{1/m,2}\lesssim\|\lambda\|_{m,2}=\|f\|_{(T_\tau^{(m,2)})^*}<\infty$ and $f(\mu)=\langle\nu,\mu\rangle$ for any $\mu\in T_\tau^{(m,2)}$.
Thus the map $\iota$ is bicontinuous.
\end{Dem}

\medskip

\subsection{Coherence and morphism property for the $\star$ product{\boldmath $.$} \hspace{0.15cm}}
\label{SubsectionCoherenceStarProduct}

We recall some algebraic tools from Bruned \& Manchon \cite{BrunedManchon} and Bailleul \& Bruned \cite{NonInvariant}.
We consider mainly preparatory trees in this section.


\subsubsection{\textit{Coherence property}.} 

We consider functions of an abstract variable 
$$
{\sf u}=({\sf u}_{\bf n})_{{\bf n}\in\bbN^2}\in\bbR^{\bbN^2}
$$
and the spacetime variable $z=(z_1,z_2)=(t,x)\in\bbR^2$.
The variable ${\sf u}$ has the role of the placeholder of derivatives $(\partial^{\bf n}u)_{{\bf n}\in\bbN^2}$.
We say that a function $f$ of $({\sf u},z)\in\bbR^{\bbN^2}\times\bbR^2$ is \textbf{cylindrical} if it depends only on a finite number of entries among $({\sf u}_{\bf n})_{{\bf n}\in\bbN^2}$ and $z$.
For any smooth cylindrical function $f$ we denote by $D_{\bf n}f\defeq\frac{\partial}{\partial {\sf u}_{\bf n}}f$ the derivative with respect to ${\sf u}_{\bf n}$ and define the derivative operator $\partial_i$ ($i\in\{1,2\}$) by
$$
\partial_if({\sf u},z) \defeq \sum_{{\bf n}\in\bbN^2}{\sf u}_{{\bf n}+{\bf e}_i}D_{\bf n}f({\sf u},z)+\partial_{z_i}f({\sf u},z).
$$
The operators $\partial_1$ and $\partial_2$ commute, so one can define the derivative operator for generic ${\bf k}=(k_1,k_2)\in\bbN^2$ by
$$
\partial^{\bf k} \defeq \partial_1^{k_1}  \partial_2^{k_2}.
$$

\medskip

\begin{defn}
For any $\tau\in\trees$ we define the smooth cylindrical function $\frak{F}_v[\tau]$ of $({\sf u},z)$ as follows. For $\tau\in\{\Xi,{\bf1}\}$, we set
\begin{equation*} \label{EqDefnFrakFElementarySymbols} 
\begin{split} 
\frak{F}_v[\Xi]({\sf u}_{\bf 0}) &\defeq  f({\sf u}_{\bf 0}),   \\
\frak{F}_v[{\bf 1}]({\sf u}_{\bf 0},{\sf u}_{(0,1)},{\sf u}_{(0,2)},z) &\defeq g({\sf u}_{\bf 0})\, {\sf u}_{(0,1)}^2 +c\, {\sf u}_{\bf 0}
+\big\{a({\sf u}_{\bf 0})-a(v(z))\big\}{\sf u}_{(0,2)}.
\end{split} \end{equation*}
We can see that ${\sf u}$ is a placeholder for $(u,\partial_xu,\partial_x^2u)$ in the right hand side of Equation \eqref{EqQgKPZBis}. For generic $\tau= X^{\bf k}\zeta\prod_{i=1}^N\mcI_{{\bf n}_i}(\tau_i)\in\trees$ we set
\begin{equation} \label{EqDefnFrakFTau}
\frak{F}_v[\tau] \defeq \Big\{\partial^{\bf k} D_{{\bf n}_1}\dots D_{{\bf n}_N} \frak{F}_v[\zeta]\Big\}\,
\prod_{i=1}^N \frak{F}_v[\tau_i].
\end{equation}
\end{defn}

\medskip

Equivalently, by setting $\frak{F}_v[p;\tau]=\partial^{\bf k} D_{{\bf n}_1}\dots D_{{\bf n}_N} \frak{F}_v[\zeta]$ for each $p\in N_\tau$ with node decoration $(\frak{t}(p),\frak{n}(p))=(\zeta,{\bf k})$ such that $N$ edges decorated by ${\bf n}_1,\dots,{\bf n}_N$ leave from $p$, we can write
\begin{equation} \label{EqDefnFrakFTau2}
\frak{F}_v[\tau]=\prod_{p\in N_\tau}\frak{F}_v[p;\tau].
\end{equation}
For instance, with 
$$
\tau_1 \defeq \mcI_{\bf0}(\Xi)\mcI_{(0,1)}(\Xi)^2, \quad 
\tau_2 \defeq X_2\mcI_{(0,2)}(\Xi)
$$ 
one has
\begin{equation*} \begin{split}
\frak{F}_v[\tau_1]({\sf u},z) &= \big\{D_{\bf0}D_{(0,1)}^2 \frak{F}_v[{\bf1}]\big\} \, \frak{F}_v[\Xi]^3({\sf u}, z) = 2g'({\sf u}_{\bf 0})f({\sf u}_{\bf 0})^3,   \\
\frak{F}_v[\tau_2]({\sf u},z) &= \big\{\partial_2D_{(0,2)} \frak{F}_v[{\bf1}]\big\} \, \frak{F}_v[\Xi]({\sf u}, z) = \Big\{a'({\sf u}_{\bf 0}){\sf u}_{(0,1)}-a'(v(z))\partial_xv(z)\Big\} \, f({\sf u}_{\bf 0}).
\end{split} \end{equation*}

\ssk

By an abuse of notation, for any contracted tree $\tau$ we write $\frak{F}_v[\tau]$ instead of $\frak{F}_v[\mcE(\tau)]$.

\medskip

\begin{prop}\label{prop:Fvcontinuous}
For any $\tau\in\mcZ(\Trees)$ and $\bsp:E_\tau\to\bbN$, we have
$$
\frak{F}_v[\tau^{\bsp}]({\sf u},z)=\big\{a({\sf u}_{\bf 0})-a(v(z))\big\}^{|\bsp|}\frak{F}_v[\tau]({\sf u},z).
$$
Consequently for any open subset $\mathcal{O}\subset\bbR^{\bbN^2}\times\bbR^2$ such that $|a({\sf u}_{\bf 0})-a(v(z))|\le m$ if $({\sf u},z)\in\mathcal{O}$, the map $\frak{F}_v$  has a unique extention as a continuous map from $T_\tau^{(m,1)}$ into $C_b(\mathcal{O})$.
\end{prop}

\medskip

\begin{Dem}
The first identity follows from $\frak{F}_v[\mcI_{\bf n}^p(\tau)]=\big\{a({\sf u}_{\bf 0})-a(v(z))\big\}^p \, \frak{F}_v[\mcI_{\bf n}^0(\tau)]$.
The second assertion follows from the definition of $\|\cdot\|_{m,1}$.
\end{Dem}

\medskip

For any solution of Equation \eqref{EqModifiedQgKPZContracted} the coefficients of $\bsw$ are written only in terms of the polynomial components of $\bsu_{\beta+\alpha_0}$. The following property was named `{\it coherence}' in \cite{BCCH}.

\medskip

\begin{prop}\label{prop:Bseries}
Let $(\bsu_{\beta+\alpha_0},\bsw)$ be a pair of modelled distributions satisfying Equation \eqref{EqModifiedQgKPZContracted}.
For any $\tau\in\trees$ with $|\tau|\le0$ and $z\in(0,T]\times\bbT$, one has
$$
\langle\mcE\bsw(z),\tau\rangle=\frak{F}_v[\tau]\big(u_{\bf0}(z),u_{(0,1)}(z),u_{(0,2)}(z),z\big),
$$
where $u_{\bf k}$ denotes the $X^{\bf k}$-component of $\bsu_{\beta+\alpha_0}$.
\end{prop}

\medskip

\begin{Dem}
First note that the function $\frak{F}_v[\tau]$ vanishes if $\tau\notin\overline{\bbB}$. If $\tau$ breaks any of the rules ({\bf P}1), ({\bf P}2) or ({\bf P}3) then it has a node $p\in N_\tau$ such that, among edges leaving from $p$, there are edges $\mcI_{\bf n}$ with ${\bf n}\notin\{{\bf 0}, (0,1), (0,2)\}$, or more than two edges $\mcI_{(0,1)}$, or more than one edge $\mcI_{(0,2)}$. In any case $\frak{F}_v[p;\tau]=0$ by definition. Therefore it is sufficient to consider $\tau={\bf1}$ and $\tau\in\basis_-\defeq\mcE(\Basis_-)$.

We assume that $\bsw\in\mcD^{\beta+\alpha_0-2,\alpha-2}(0,2T)$ has an expansion
$$
\mcE\bsw(z)=w_{\bf 0}(z){\bf1}+\sum_{\tau\in\basis_-}w_\tau(z)\tau.
$$
Since $\widetilde{\bsw}_T\vert_{(0,T]\times\bbT}=\bsw\vert_{(0,T]\times\bbT}$, then we can write
$$
\mcE\bsu_{\beta+\alpha_0}(z)=\sum_{|{\bf k}|_\mfs\le2}u_{\bf k}(z)X^{\bf k}+\sum_{\tau\in\basis_-}w_\tau(z)\mcI_{\bf 0}(\tau)
$$
on $z\in(0,T]\times\bbT$ for some $\{u_{\bf k}\}_{|{\bf k}|_\mfs\le2}$.
From now, we fix $z\in(0,T]\times\bbT$ and consider the right hand side of the second equation of \eqref{EqModifiedQgKPZContracted}. 
First, the coefficients of 
$$
F(\bsu_\beta)={\sf Q}_{<\beta}\Bigg(\sum_{n=0}^\infty\frac{f^{(n)}(u_{\bf0})}{n!}(\bsu_\beta-u_{\bf0}{\bf1})^n\Bigg)
$$
can be written in terms of the coefficients $u_{\bf0},u_{(0,1)}$, and $\{w_\tau\}$. 
Indeed, for any $\tau=\prod_{i=1}^N(\mcI_{\bf0}(\tau_i))^{m_i}\in\basis$ with $|\tau|<\beta$ and $\tau_i\neq\tau_j$ for $i\neq j$, the $\tau$-component $F(\bsu_\beta)_\tau$ of $\mcE F(\bsu_\beta)$ is given by
\begin{align*}
F(\bsu_\beta)_\tau
&= \frac{f^{(m_1+\cdots m_N)}(u_{\bf 0})}{(m_1+\cdots+m_N)!} \, \frac{(m_1+\cdots+m_N)!}{m_1!\cdots m_N!} \, w_{\tau_1}^{m_1}\cdots w_{\tau_N}^{m_N}\\
&= \frac{w_{\tau_1}^{m_1}\cdots w_{\tau_N}^{m_N}}{m_1!\cdots m_N!} \, 
D_{\bf 0}^{m_1+\cdots+m_N}f({\sf u})\vert_{{\sf u}_{\bf 0}=u_{\bf 0}}.
\end{align*}
Similarly, for $\tau=X_2\prod_{i=1}^N(\mcI_{\bf0}(\tau_i))^{m_i}\in\basis$ with $|\tau|<\beta$ and $\tau_i\neq\tau_j$ for $i\neq j$, we have
\begin{align*}
F(\bsu_\beta)_\tau
&= \frac{f^{(m_1+\cdots m_N+1)}(u_{\bf 0})}{(m_1+\cdots+m_N+1)!} \, \frac{(m_1+\cdots+m_N+1)!}{m_1!\cdots m_N!} \, w_{\tau_1}^{m_1}\cdots w_{\tau_N}^{m_N}u_{(0,1)}\\
&= \frac{w_{\tau_1}^{m_1}\cdots w_{\tau_N}^{m_N}}{m_1!\cdots m_N!} \, \partial_2D_{\bf 0}^{m_1+\cdots m_N}f^{(m_1+\cdots m_N)}({\sf u}) \vert_{({\sf u}_{\bf 0},{\sf u}_{(0,1)})=(u_{\bf 0},u_{(0,1)})}.
\end{align*}
Note that, since $\min\{|\tau|\,;\,\tau\in\basis\}=\alpha_0-2>-2$, any $\tau\in\basis_-$ cannot have a node decorated by $X^{\bf k}$ with $|{\bf k}|_\mfs\ge2$.
Similar formulas hold for $G(\bsu_\beta)$ and $A(\bsu_\beta)$.
We also have
\begin{align*}
\mcE(\bsD\bsu_\beta)^2
&=u_{(0,1)}^2{\bf1}+2\sum_{\tau\in\basis_-,\,|\tau|<\beta+\alpha_0-3}u_\tau u_{(0,1)}\mcI_{(0,1)}(\tau)\\
&\quad+2\sum_{\tau\in\basis_-,\,|\tau|<\frac{\beta+\alpha_0-4}2}u_\tau^2\mcI_{(0,1)}(\tau)^2
+\sum_{\tau,\sigma\in\basis_-,\,\tau\neq\sigma,\,|\tau|+|\sigma|<\beta+\alpha_0-4}u_\tau u_\sigma\mcI_{(0,1)}(\tau)\mcI_{(0,1)}(\sigma)\\
\mcE\bsD^2\bsu_{\beta+\alpha_0}&=u_{(0,2)}{\bf1}+\sum_{\tau\in\basis_-}u_\tau\mcI_{(0,2)}(\tau).
\end{align*}
Inserting these expansions into the second equation of \eqref{EqModifiedQgKPZContracted}, we have
$$
w_{\bf 0} = g(u_{\bf 0})u_{(0,1)}^2 + cu_{\bf 0} + \big\{a(u_{\bf 0})-a(v)\big\} u_{(0,2)} = \frak{F}_v[{\bf 1}]\big(u_{\bf 0},u_{(0,1)},u_{(0,2)},z\big)
$$
and, for $\tau\in\basis_-$,
\begin{align*}
w_\tau=
\left\{
\begin{aligned}
&F(\bsu_\beta)_\sigma
&&\textrm{ if } \tau=\sigma\Xi,\, \sigma\in\basis_{\bf0},   \\
&G(\bsu_\beta)_\eta w_\sigma u_{(0,1)}
&&\textrm{ if } \tau=\mcI_{(0,1)}(\sigma)\eta,\, \eta\in\basis_{\bf0},   \\
&2G(\bsu_\beta)_\eta w_\sigma^2
&&\textrm{ if } \tau=\mcI_{(0,1)}(\sigma)^2\eta,\, \eta\in\basis_{\bf0},   \\
&G(\bsu_\beta)_\eta w_{\sigma_1}w_{\sigma_2}
&&\textrm{ if } \tau=\mcI_{(0,1)}(\sigma_1)\mcI_{(0,1)}(\sigma_2)\eta,\, \eta\in\basis_{\bf0},\, \sigma_1\neq \sigma_2,   \\
&\{a(u_{\bf 0})-a(v)\}w_\sigma
&&\textrm{ if } \tau=\mcI_{(0,2)}(\sigma),   \\
&\{a'(u_{\bf 0})u_{(0,1)}-a'(v)\partial_xv\}w_\sigma
&&\textrm{ if } \tau=X_2\mcI_{(0,2)}(\sigma),   \\
&A(\bsu_\beta)_\eta w_\sigma
&&\textrm{ if } \tau=\mcI_{(0,2)}(\sigma)\eta,\, \eta\in\basis_{\bf0}\setminus\{X^{\bf k}\},
\end{aligned}
\right.
\end{align*}
where $\basis_{\bf0}$ denotes the set of all preparatory trees of the form $X^{\bf k}\prod_{i=1}^N(\mcI_{\bf0}(\tau_i))^{m_i}$. Starting from $w_{\Xi}=f(u_{\bf 0})=\frak{F}_v[\Xi](u_{\bf 0})/S(\Xi)$, we can recursively show that $w_\tau=\frak{F}_v[\tau]/S(\tau)$ for any $\tau\in\basis_-\cup\{{\bf 1}\}$ from the above relations.
\end{Dem}


\subsubsection{\textit{Star product}.}
We recall from Section 2 of \cite{NoExtendedDecoration} some bilinear operators on $\overline{\sf T}$. For $\tau\in\trees$ with a decoration $\frak{n}:N_\tau\to\bbN^2$, $v\in N_\tau$, and ${\bf k}\in\bbZ^2$ such that $\frak{n}(v)+{\bf k}\ge{\bf 0}$, we define
$$
\uparrow_v^{\bf k}\tau
$$
as the preparatory tree given by changing the decoration $\frak{n}(v)$ at $v$ into $\frak{n}(v)+{\bf k}$. For $\sigma,\tau\in\trees$ and ${\bf n}\in\bbN^2$, set

\begin{equation}\label{Eq:sigmacurvetau}
\sigma \curvearrowright_{\bf n} \tau
\defeq  \sum_{v\in N_{\tau}}\sum_{{\bf k}\in\bbN^2,\, {\bf k}\le\frak{n}(v)\wedge{\bf n}}{\frak{n}(v) \choose {\bf k}} \,\sigma \curvearrowright^v_{{\bf n}-{\bf k}}(\uparrow_v^{-{\bf k}}\tau),
\end{equation}
where the operator $\curvearrowright^v_{{\bf n}-{\bf k}}$ grafts $\sigma$ onto $\tau$ at the node $ v $ with an edge of type $\mcI_{{\bf n}-{\bf k}}$. 
Since $\mcI_{{\bf n}-{\bf k}}$ increases the homogeneity by $2-|{\bf n}|_\mfs+|{\bf k}|_\mfs$ and $\uparrow_v^{-{\bf k}}$ reduces the homogeneity by $|{\bf k}|_\mfs$, all terms in the right hand side have the homogeneity $|\sigma|+|\tau|+2-|{\bf n}|_\mfs$.
We extend $\curvearrowright_{\bf n}$ into the bilinear map $\curvearrowright_{\bf n}:\overline{\sf T}\times\overline{\sf T}\to\overline{\sf T}$.
The family of operators $\{\curvearrowright_{\bf n}\}_{{\bf n}\in\bbN^2}$ satisfies the \textbf{multiple pre-Lie} property as discussed in \cite{BrunedManchon, BCCH}.

In addition, we consider the bilinear map $\star:\overline{\sf T}_\bullet\times\overline{\sf T}\to\overline{\sf T}$ as follows. 
%
%
%
%
First define for $\tau\in\trees$, $B\subset N_{\tau}$, and ${\bf k}\in\bbN^2$, the derivation map $\uparrow^{\bf k}_{B}$ by
$$
\uparrow^{\bf k}_{B} \tau = \sum_{({\bf k}_v)_{v\in B}\subset \bbN^2,\, \sum_{v \in B} {\bf k}_v  = {\bf k} } \, \bigg(\prod_{v \in B}\uparrow_v^{{\bf k}_v}\bigg) \tau.
$$
Moreover, for $\tau,\sigma\in\trees$, $B\subset N_\tau$, and ${\bf n}\in\bbN^2$, we define $\sigma \curvearrowright_{\bf n}^B \tau$ by the identity \eqref{Eq:sigmacurvetau} where the region of the sum for $v$ is restricted into $B$.
Set finally for $\sigma = X^{\bf k}\prod_{i=1}^m \mcI_{{\bf n}_i}(\sigma_i)\in\trees_\bullet$ and $\tau\in\trees$,
$$
\sigma\star\tau\defeq
 \uparrow^{\bf k}_{N_\tau}\big\{
\sigma_1 \curvearrowright_{{\bf n}_1}^{N_\tau} 
(\sigma_2\curvearrowright_{{\bf n}_2}^{N_\tau}
(\dots
(\sigma_m\curvearrowright_{{\bf n}_m}^{N_\tau}\tau)
\dots)
)\big\}.
$$
Some properties of the $\star$ product are discussed in Section 3.3 of \cite{BrunedManchon}.

\medskip

\begin{prop}\label{prop:deltastardual}
\begin{itemize}
\item[(1)]
(Proposition 3.15 of \cite{BrunedManchon})
The $\star$ product is associative in the sense that
$$
{\tau}\star({\sigma}\star\eta)=({\tau}\star{\sigma})\star\eta
$$
for any ${\tau},{\sigma}\in\trees_\bullet$ and $\eta\in\trees$.
\item[(2)]
(Theorems 3.16 and 4.2 of \cite{BrunedManchon})
The $\star$ product is a dual of the coproduct $\icprod:\treesp\to\treesp\mathbin{\widehat{\otimes}}\treesp_\bullet$ in the sense that
$$
\langle\sigma\star\tau,\eta\rangle=\langle\tau\otimes\sigma,\icprod\eta\rangle
$$
for any $\sigma\in\trees_\bullet$ and $\tau,\eta\in\trees$, where the inner product $\langle\cdot,\cdot\rangle$ of $\treesp$ is linearly extended to $\treesp\mathbin{\widehat{\otimes}}\treesp$ by setting
$$
\langle\tau_1\otimes\tau_2,\sigma_1\otimes\sigma_2\rangle\defeq 
\langle\tau_1,\sigma_1\rangle\langle\tau_2,\sigma_2\rangle
$$
for any $\tau_1,\tau_2,\sigma_1,\sigma_2\in\trees$.
\end{itemize}
\end{prop}

\medskip

The following morphism property of $\frak{F}_v$ with respect to the $\star$ product plays in our argument the role that pre-Lie morphisms played in the original approach of Bruned, Chandra, Chevyrev \& Hairer \cite{BCCH}. The interest of the formula \eqref{EqMorphismFrakFTau} will appear below in Proposition \textit{\ref{PropMorphismPropertyRStar}}. 

\medskip

\begin{prop} \label{PropMorphism}
For any $N\in\bbN$, $\tau,\sigma_1,\dots,\sigma_N\in\trees$, ${\bf k},{\bf n}_1,\dots,{\bf n}_N\in\bbN^2$, one has
\begin{equation} \label{EqMorphismFrakFTau} \begin{split}
\frak{F}_v\bigg[\Big(X^{\bf k} \prod_{i=1}^N \mcI_{{\bf n}_i}(\sigma_i)\Big) \star \tau\bigg]
= \Big\{\partial^{\bf k} D_{{\bf n}_1} \dots D_{{\bf n}_N} \frak{F}_v[\tau]\Big\} \, \prod_{i=1}^N \frak{F}_v[\sigma_i].
\end{split} \end{equation}
\end{prop}

\medskip

\begin{Dem}
The proof is the same as that of Proposition {\it 2} in \cite{NoExtendedDecoration}. The only ingredient to be checked is the following identity used in \cite{NoExtendedDecoration}: for any ${\bf n}_1,\dots,{\bf n}_N\in\bbN^2$ and ${\bf k}\in\bbN^2$ one has
\begin{equation} \label{EqFundamentalRelation}
\sum_{\substack{ {\bf k}_1,\dots,{\bf k}_N \in \bbN^2 \\ {\bf k}_1+\cdots+{\bf k}_N\le{\bf k},\, {\bf k}_i\le{\bf n}_i}} {{\bf k} \choose {{\bf k}_1,\dots,{\bf k}_N}} \, \partial^{{\bf k}-\sum_{i=1}^N{\bf k}_i}
\prod_{i=1}^N D_{{\bf n}_i - {\bf k}_i} = D_{{\bf n}_1}\dots D_{{\bf n}_N}\partial^{{\bf k}},
\end{equation}
where
$$
{{\bf k} \choose {{\bf k}_1,\dots,{\bf k}_N}}
\defeq
\frac{{\bf k}!}{{\bf k}_1!\cdots{\bf k}_N!({\bf k}-\sum_{i=1}^N{\bf k}_i)!}.
$$
For the case that $N=1$ and ${\bf k}={\bf e}_i$, the identity is checked directly.
Indeed, by the Leibniz rule,
\begin{align*}
D_{\bf n}\partial_if
&=\sum_{{\bf m}}D_{\bf n}({\sf u}_{{\bf m}+{\bf e}_i}D_{\bf m}f)+D_{\bf n}\partial_{z_i}f\\
&={\bf1}_{{\bf n}\ge{\bf e}_i}D_{{\bf n}-{\bf e}_i}f+\sum_{{\bf m}}{\sf u}_{{\bf m}+{\bf e}_i}D_{\bf n}D_{\bf m}f+D_{\bf n}\partial_{z_i}f\\
&=\partial_iD_{\bf n}f+{\bf1}_{{\bf n}\ge{\bf e}_i}D_{{\bf n}-{\bf e}_i}f.
\end{align*}
By induction this identity can be extended into 
$$
D_{\bf n}\partial^{\bf k}f=\sum_{{\bf k}_1\le{\bf k}\wedge{\bf n}} {{\bf k} \choose {\bf k}_1}\partial^{{\bf k}-{\bf k}_1}D_{{\bf n}-{\bf k}_1}f
$$
for all ${\bf k}\in\bbN^2$ first, and to \eqref{EqMorphismFrakFTau} next.
\end{Dem}

\medskip

\subsection{Strong preparation maps and their associated models{\boldmath $.$} \hspace{0.15cm}}
\label{SubsectionPreparationMaps}

In this section we recall from \cite{BrunedRecursive} and \cite{NoExtendedDecoration} the building blocks of an inductive construction of a renormalized model and a renormalized equation, and give a detailed statement of Assumption {\it\ref{asmp1}}.

\subsubsection{Preparation maps.}

For $\tau\in \Trees$, denote by $\vert\tau\vert_{\Xi}$ the number of noise symbols $\Xi$ that appear in $\tau$. Recall from Bruned's work \cite{BrunedRecursive} that a \textit{\textbf{preparation map}} is a linear map $R : \Treesp^\Basis\rightarrow \Treesp^\Basis$ with the following properties.
\begin{itemize}
\item[(a)]
One has 
\begin{equation*}
R(\zeta)=\zeta\qquad(\zeta\in\{\Xi,{\bf1}\}),\qquad
R(X^{\bf k}\tau)=X^{\bf k}R(\tau),\qquad
R(\mcI_{\bf n}^p(\tau))=\mcI_{\bf n}^p(\tau).
\end{equation*}
\item[(b)]
For each $\tau\in\Basis$ there exist {\it finitely many} $\tau_i \in \Basis$ and constants $\lambda_i$ such that
\begin{equation*} 
R \tau = \tau + \sum_i \lambda_i \tau_i, \quad\textrm{with}\quad \vert\tau_i\vert> \vert\tau\vert \quad\textrm{and}\quad |\tau_i|_{\Xi} < |\tau|_{\Xi}.
\end{equation*} 
\item[(c)]
One has the `commutation' relation
\begin{equation} \label{Commutation_R}
( R \otimes \id ) \Cprod = \Cprod R.
\end{equation}
\end{itemize}
The role of $R$ is to provide a recursive definition of the product of two trees in terms of other trees that have already been `renormalized'. Its use in Section {\sf\ref{SectionAdmissibleModelsPrepMaps}} in the recursive definition of the renormalized model will make that point clear -- see in particular \eqref{defPiR}. Accordingly the second and third identities of (a) account for the fact that there is no need to `renormalize' elements of the form $X^{\bf k}\tau$ and $\mcI_{\bf n}^p(\tau)$ if the element $\tau$ has already been renormalized. 
We introduce a renormalization character to construct a concrete preparation map.

\medskip

\begin{defn}
A \textbf{(spacetime-dependent) renormalization character of a growth factor $m_0>0$} is a map
$$
\vcst : (\bbR\times\bbT)\times\bbB_-\to\bbR
$$
which is continuous in $\bbR\times\bbT$ and vanishes on the elements of the forms
$$
X^{\bf k}\tau\ \ ({\bf k}\neq0),\qquad
\mcI_{\bf n}^p(\tau),
$$
and such that for any $\tau\in\bbF_-$ there exists a constant $C(\tau)$ such that
\begin{equation}\label{CharacterGrowth}
\big|\vcst(z,\tau^{\bsp})\big|\le C(\tau)\, m_0^{|\bsp|}
\end{equation}
for any $\bsp:E_\tau\to\bbN$ and $z\in\bbR\times\bbT$.
\end{defn}

\medskip

The following lemma is needed to show Proposition {\it\ref{prop:Rellcomplete}}.

\medskip

\begin{lem}\label{lem:completeness}
For $\tau\in\Basis$ and $v\in N_\tau$, we denote by $\tau_{\ge v}$ the unique subtree of $\tau$ consisting of all nodes $w\in N_\tau$ such that $w\ge v$. Then for any $v\in N_\tau\setminus\{\varrho_\tau\}$, one has $|\tau_{\ge v}|\le|\tau|$.
Moreover, for any $u,v\in N_\tau\setminus\{\varrho_\tau\}$ with $u\neq v$, one has $|\tau_{\ge u}|+|\tau_{\ge v}|+2\le|\tau|$.
\end{lem}

\medskip

\begin{Dem}
For the first assertion, it is sufficient to show that $|\tau_{\ge v}|\le|\tau_{\ge u}|$ for the unique node $u$ such that $\{u,v\}\in E_\tau$ and $u\le v$. The possible form of $\tau_{\ge u}\in\bbB$ with the minimal homogeneity is $\tau_{\ge u}=\mcI_{(0,2)}^p(\tau_{\ge v})$. In this case, we have $|\tau_{\ge v}|=|\tau_{\ge u}|$.

For the second assertion, it is sufficient to consider the case that $\{w,u\},\{w,v\}\in E_\tau$ for some $w\in N_\tau$. The possible form of $\tau_{\ge w}\in\bbB$ with the minimal homogeneity is $\tau_{\ge w}=\mcI_{(0,1)}^p(\tau_{\ge u})\mcI_{(0,1)}^q(\tau_{\ge v})$ or $\mcI_{(0,2)}^p(\tau_{\ge u})\mcI_{\bf 0}^q(\tau_{\ge v})$. In either case we have $|\tau_{\ge w}|=|\tau_{\ge u}|+|\tau_{\ge v}|+2$.
\end{Dem}

\medskip

\begin{prop}\label{prop:Rellcomplete}
For any renormalization character $\ell$ of a growth factor $m_0>0$,
we can define the map $R_\ell:(\bbR\times\bbT)\times\Treesp^\Basis\to\Treesp^\Basis$ which is continuous in $\bbR\times\bbT$ and linear in $\Treesp^\Basis$ by
$$
R_\ell(z,\tau)\defeq \tau+\big(\ell(z,\cdot)\otimes\id\big)\iCprod\tau\qquad(\tau\in\Basis),
$$
where we extend $\ell(z,\cdot)$ linearly by setting $\ell(z,\tau)=0$ if $\tau\in\Basis\setminus\Basis_-$.
For any $z\in\bbR\times\bbT$, the linear map $R_\ell(z)\defeq R_\ell(z,\cdot)$ is a preparation map.
Moreover, the map $R_\ell$ is extended into the continuous map $(\bbR\times\bbT)\times T^{(m,1)}\to T^{(m,1)}$ for any $m\ge m_0$ and the continuous map $(\bbR\times\bbT)\times T^{(m,2)}\to T^{(m_0,2)}$ for any $m> m_0$
\end{prop}

\medskip

\begin{Dem}
We assume that the (contracted version of) expansion \eqref{eq:expliciticprod} of $\tau$ is given by
\begin{equation}\label{eq:proof:icprodtau}
\iCprod\tau=\sum_{\sigma,\eta}c_{\sigma\eta}^\tau\, \sigma\otimes\eta,
\end{equation}
where $\sigma$ and $\eta$ are some contracted trees. Since there only finitely many $\sigma\in\Basis_-$ in the above sum, $R_\ell(z)\tau$ becomes a finite linear combination. Let $\sigma\in\Basis_-$ and define $\partial N_\sigma=\{u\in N_\sigma\,;\, \{u,v\}\in\partial\sigma\}$. If all edges leaving from $\partial N_\sigma$ are of the form $\mcI_{\bf 0}^p$, then $\eta\in\Basis$. Indeed, since $\sigma$ cannot have the node decoration greater than $(0,1)$ to have the negative homogeneity, all edges leaving from $\varrho_\eta$ are of the form $\mcI_{\bf 0}^p$ except at most one edge of the form $\mcI_{(0,1)}^q$. Next we assume that there exists $v\in\partial N_\sigma$ from which an edge of the form $\mcI_{(0,2)}^p$ leaves. Since $\frak{t}(v)={\bf1}$, this $v$ cannot be a leaf within $\sigma$ and there is another edge of the form $\mcI_{\bf0}^q$ inside $\sigma$ leaving from $v$. Then by the first assertion of Lemma {\it\ref{lem:completeness}}, we have the contradiction $|\sigma|\ge|\mcI_{\bf0}^q(\Xi)|>0$. We also have a contradiction if we assume that there exists $v\in\partial N_\sigma$ from which two edge of the form $\mcI_{(0,1)}^p$ leave.
Finally, we assume that there exist $u,v\in\partial N_\sigma$ such that $u\neq v$ and there are edges of the form $\mcI_{(0,1)}^p$ leave from $u$ and $v$ respectively. By the rule (C2), there are other edges of the form $\mcI_{\bf0}^q$ or $\mcI_{(0,1)}^q$ inside $\sigma$ leaving from $u$ and $v$ respectively. If $u,v\neq\varrho_\tau$, then by the second assertion of Lemma {\it\ref{lem:completeness}}, we have the contradiction $|\sigma|\ge|\mcI_{(0,1)}^p(\Xi)|+|\mcI_{(0,1)}^q(\Xi)|+2=\alpha_0>0$. If $u=\varrho_\tau$, there exists a unique $w\in N_\sigma$ such that $u\le w\le v$ and $\{w,v\}\in E_\sigma$. Since $\{w,v\}$ is of the form $\mcI_{\bf0}^q$ or $\mcI_{(0,1)}^q$, by the first assertion of Lemma {\it\ref{lem:completeness}}, we have the contradiction $|\sigma|\ge1+|\sigma_{\ge w}|\ge1+|\mcI_{(0,1)}^q(\Xi)|=\alpha_0>0$.

We can show the properties (a) and (b) of preparation maps by the representation \eqref{eq:expliciticprod} and the fact that the equality $|\tau|=|\sigma|+|\eta|$ holds for the expansion \eqref{eq:proof:icprodtau}. We have the property (c) from the associativity of $\iCprod$ as follows.
\begin{align*}
\Cprod\big(R_\ell(z)-\id\big)&=\big(\ell(z,\cdot)\otimes\id\otimes p_+\big)(\id\otimes\iCprod)\iCprod\\
&=\big(\ell(z,\cdot)\otimes\id\otimes p_+\big)(\iCprod\otimes\id)\iCprod
=\big\{\big(R_\ell(z)-\id\big)\otimes\id\big\}\Cprod.
\end{align*}
Since $\ell(z,\cdot)$ is continuous with respect to the norm $\|\cdot\|_{m,1}$ uniformly over $z$ for any $m\ge m_0$ by the condition \eqref{CharacterGrowth}, the continuity of $R_\ell$ follows from Propositions {\it \ref{prop:coproductcontinuous}} and {\it\ref{prop:T1T2embedding}}.
\end{Dem}

\medskip

The dual operator of $R_\ell(z)$ played an important role in \cite{NoExtendedDecoration}.

\medskip

\begin{prop}\label{prop:SRdual}
For any renormalization character $\ell$ of a growth factor $m_0>0$ and any $m>m_0$, we can define the continuous map $S_\ell:(\bbR\times\bbT)\times T^{(1/m_0,2)}\to T^{(1/m,2)}$ which is linear in $T^{(1/m_0,2)}$ by
$$
S_\ell(z,\tau)\defeq \tau+\sum_{\sigma\in\Basis_-}\frac{\ell(z,\sigma)}{S(\sigma)}p_{\Basis}\big(\pi(\mcE\tau\star\mcE\sigma)\big)\qquad(\tau\in\Basis),
$$
where $p_{\Basis}:\Treesp\to\Treesp^{\Basis}$ denotes the canonical projection.
The linear map $S_\ell(z)\defeq S_\ell(z,\cdot)$ is a dual of $R_\ell(z)$ in the sense that
$$
\langle R_\ell(z)\tau,\sigma\rangle=\langle\tau,S_\ell(z)\sigma\rangle
$$
for any $\tau,\sigma\in\Basis$.
\end{prop}

\medskip

We need the projection $p_{\Basis}$ because $\mcE\tau\star\mcE\sigma$ may produce non-conforming trees even if $\mcE\tau$ and $\mcE\sigma$ strongly conform to the rule {\bf P}.

\medskip

\begin{Dem}[of Proposition \ref{prop:SRdual}]
For each $N\in\bbN$, we define the renormalization character $\ell_N$ by setting
$$
\ell_N(z,\tau^{\bsp})={\bf1}_{|\bsp|\le N}\ell(z,\tau^{\bsp}).
$$
Then the duality between $\icprod$ and $\star$ (Proposition {\it\ref{prop:deltastardual}} (2)) implies the duality between $R_{\ell_N}(z)$ and $S_{\ell_N}(z)$ as linear maps on $\Treesp^{\Basis}$.
Since $R_{\ell_N}\tau\to R_\ell\tau$ in the norm $\|\cdot\|_{m_0,2}$ for any $\tau\in T^{(m,2)}$, we have the weak convergence of $S_{\ell_N}\sigma$ in $T^{(1/m,2)}$ for any $\sigma\in T^{(1/m_0,2)}$ by Proposition {\it\ref{prop:dualityTmT1/m}}, and the map 
$$
S_\ell(z)\defeq\lim_{N\to\infty}S_{\ell_N}(z):T^{(1/m_0,2)}\to T^{(1/m,2)}
$$ 
is also a dual of $R_\ell(z)$. The continuity of $S_\ell$ follows from that of $R_\ell$.
\end{Dem}

\medskip

The next proposition follows from Proposition \textit{\ref{PropMorphism}}.

\medskip

\begin{prop} \label{PropMorphismPropertyRStar}
Fix $0<m_0<m$. Let $\mathcal{O}$ be an open set subset of $\bbR^{\bbN^2}\times\bbR^2$ such that 
$$
\big|a({\sf u}_{\bf 0})-a(v(z))\big| \leq m \qquad \textrm{if} \qquad ({\sf u},z)\in\mathcal{O}
$$ 
and let $\ell$ be a renormalization character of growth factor $m_0$. Then for every preparatory tree $\tau=X^{\bf k}\zeta\prod_{i=1}^N\mcI_{{\bf n}_i}(\tau_i)\in\basis$ and $({\sf u},z)\in\mathcal{O}$ one has
\begin{equation*} \begin{split}
\frak{F}_v\big[S_\ell(z)\pi(\tau)\big]({\sf u},z) = \Big\{ \partial^{\bf k} D_{{\bf n}_1} \dots D_{{\bf n}_N} \frak{F}_v\big[S_\ell(z)\zeta\big]\Big\}\,\prod_{i=1}^N \frak{F}_v[\tau_i]({\sf u},z).
\end{split} \end{equation*}
\end{prop}

\medskip

\begin{Dem} 
Since we can write
$$
\tau
=\Big(X^{\bf k} \prod_{i=1}^N \mcI_{{\bf n}_i}(\tau_i)\Big)\star\zeta
\eqdef\eta\star\zeta,
$$ 
by using the fact that $\frak{F}_v[p_{\Basis}\tau]=\frak{F}_v[\tau]$ for any $\tau\in\trees$ (see the proof of Proposition {\it\ref{prop:Bseries}}) and the morphism property (Proposition {\it\ref{PropMorphism}}), one gets
\begin{align*}
\frak{F}_v[(S_\ell(z)-\id)\pi(\tau)]
&= \sum_{\sigma\in\Basis_-}\frac{\ell(z,\sigma)}{S(\sigma)}
\frak{F}_v\big[p_{\Basis}\circ\pi( \eta\star\zeta\star\mcE\sigma) \big]
= \sum_{\sigma\in\Basis_-}\frac{\ell(z,\sigma)}{S(\sigma)}
\frak{F}_v\big[\eta\star\zeta\star\mcE\sigma \big]   \\
&= \sum_{\sigma\in\Basis_-}\frac{\ell(z,\sigma)}{S(\sigma)}\Big\{\partial^{\bf k} D_{{\bf n}_1} \dots D_{{\bf n}_N} \frak{F}_v[\zeta\star\mcE\sigma]\Big\} \, \prod_{i=1}^N \frak{F}_v[\tau_i]\\
&= \Big\{ \partial^{\bf k} D_{{\bf n}_1} \dots D_{{\bf n}_N} \frak{F}_v\big[(S_\ell(z)-\id)\zeta\big]\Big\}\,\prod_{i=1}^N \frak{F}_v[\tau_i].
\end{align*}
This concludes the proof.
\end{Dem}

\smallskip

\subsubsection{Admissible model associated with a preparation map.}
\label{SectionAdmissibleModelsPrepMaps}

Fix a regularization parameter $\varepsilon\in(0,1]$ and denote by 
$$
\xi^\varepsilon \in C^\infty(\textbf{\textsf{R}}\times\textbf{\textsf{T}})
$$ 
a regularized version of the spacetime white noise $\xi$. For any spacetime-dependent renormalization character $\ell$ we define inductively the maps ${\sf \Pi}^{\ell}$ and $\Pi^{\ell,\times}$ as follows.
\begin{equation} \label{defPiR} \begin{split}
&{\sf \Pi}^{\ell}\Xi =  \Pi^{\ell,\times}\Xi \defeq \xi^\varepsilon,
\qquad({\sf \Pi}^{\ell}X^{\bf k})(z) =  (\Pi^{\ell,\times}X^{\bf k})(z)\defeq z^{\bf k},\\
&{\sf \Pi}^{\ell}  = \Pi^{\ell,\times} \circ R_\ell, \qquad \Pi^{\ell,\times}(\tau \sigma) = \big(\Pi^{\ell,\times}\tau\big) \big(\Pi^{\ell,\times} \sigma\big),   \\
&\Pi^{\ell,\times}\big(\mcI^p_{\bf n}(\tau)\big) = (\partial_z^{\bf n} K^{v})\big\{\big(\partial_x^2 K^{v})^{p} ({\sf \Pi}^{\ell} \tau)\big\}.
\end{split} \end{equation}
The operator $\partial_x^2K^{v}$ makes sense here because ${\sf\Pi}^{\ell}\tau$ constructed as above belongs to $\mcC_{\mfs}^{0+}(\bbR\times\bbT)$. (Note that $K^{v}$ maps $\mcC_{\mfs}^{0+}(\bbR\times\bbT)$ into $\mcC_{\mfs}^{2+}(\bbR\times\bbT)$ -- see Theorem {\it \ref{thm:SchauderK+R}}.) As $R$ is spacetime-dependent the third identity in \eqref{defPiR} reads
$$
\big({\sf \Pi}^{\ell}\tau\big)(z)  = \Pi^{\ell,\times}\big(R_\ell(z)\tau\big)(z),
$$
for all $z\in\bbR\times\bbT$ and all $\tau\in\Basis$. From this definition and the properties of the preparation map $R_\ell(z)$, it follows that the map $ {\sf \Pi}^{\ell}$ satisfies the admissibility conditions 
$$
 {\sf \Pi}^{\ell}\big(\mcI^p_{\bf n}(\tau)\big) = (\partial_z^{\bf n} K^{v}) \big\{\big(\partial_x^2 K^{v})^{p} ({\sf \Pi}^{\ell} \tau)\big\},\qquad
 {\sf\Pi}^\ell(X^{\bf k}\tau)(z)=z^{\bf k}({\sf\Pi}^\ell\tau)(z).
$$
As mentioned after Definition {\sf\ref{defn:smoothadmissiblemodels}}, the map ${\sf\Pi}^\ell$ as above determines the smooth admissible model ${\sf M}^\ell=({\sf g}^\ell,{\sf\Pi}^\ell)$.
It also should be note that, the operator ${\sf \Pi}_z^\ell\defeq({\sf\Pi}^\ell)_z^{{\sf g}^\ell}=\big({\sf \Pi}^\ell\otimes({\sf g}_z^\ell)^{-1}\big)\Delta$ admits the factorization
\begin{equation}\label{eq:PixfactorizationPitimesR}
({\sf\Pi}_z^\ell\tau)(w)=\big(\Pi_z^{\ell,\times}R_\ell(w)\tau\big)(w)
\end{equation}
with the algebra morphism $\Pi_z^{\ell,\times}$ recursively defined as follows.
\begin{equation*} 
\begin{split}
&\Pi_z^{\ell,\times}\Xi \defeq \xi^\varepsilon,
\qquad(\Pi_z^{\ell,\times}X^{\bf k})(w)\defeq (w-z)^{\bf k},\\
&\Pi_z^{\ell,\times}\mcI_{\bf n}\tau
= \big(\partial^{\bf n}K^v\big)({\sf\Pi}_z^{\ell}\tau)-\sum_{|{\bf k}|_\mfs<|\tau|+2-|{\bf n}|_\mfs}\frac{(\cdot-z)^{\bf k}}{{\bf k}!} \, \big(\partial^{{\bf n}+{\bf k}}K^v\big)({\sf\Pi}_z^{\ell}\tau)(z)
\end{split} 
\end{equation*}
Among the renormalization characters we are interested in the character $\vcst^\varepsilon_{v}(z,\cdot)$ is defined as in Section 6.3 of Bruned, Hairer \& Zambotti's work \cite{BHZ}. The associated model ${\sf M}^\varepsilon={\sf M}^{\ell_v^\varepsilon}$ is called the \textit{\textbf{BPHZ model}}. Although ${\sf M}^\varepsilon$ depends on $v$ we emphasize only the $\varepsilon$-dependence to lighten the notations.

\medskip

\begin{asmp} \label{asmp1}
There exists a character $\vcst_{v}^\varepsilon$ of growth factor $m_0>0$ for each $\varepsilon\in(0,1]$ (the constants $C(\tau)$ in \eqref{CharacterGrowth} may be $\varepsilon$-dependent) such that, for any $m>m_0$, the BPHZ renormalized model ${\sf M}^\varepsilon$ is convergent in the space $\mModels$ as $\varepsilon$ goes to $0$.
\end{asmp}

\medskip

We conjecture that \textbf{\textit{Assumption \ref{asmp1}}} holds true in the regime $\alpha_0>\frac14$, but we do not discuss it precisely in this paper. Such a convergence result was proved in several works \cite{CH,HS,RandomModel} in semilinear settings, but we need some modifications. For instance we cannot directly use \cite{CH} because the kernel $\partial_x^2K^{v}$ is too singular to be integrable around the origin. We would be able to overcome this difficulty by considering $K^{v}(\partial_x^2K^{v})^{p}$ as a single integrable kernel -- see Proposition {\it \ref{prop:Lpkernelestimate}} below. The inductive proofs in \cite{HS, RandomModel} are also not directly applicable because $K^{v}(z,w)$ is not translation invariant. However we expect that the inductive proof as in \cite{RandomModel} would reduce the effort for the proof significantly when working with a noise satisfying a Poincar\'e-type inequality.
See the discussion at the beginning of Section {\sf\ref{SubsectionAsmp1}}.

\bigskip

\subsection{Renormalized equation{\boldmath $.$} \hspace{0.15cm}}
\label{SubsectionRenormalizedEquation}

Denote by ${\sf R}^\varepsilon$ the reconstruction map associated with ${\sf M}^\varepsilon$. The proof of Theorem \textit{10} in \cite{NonInvariant} works verbatim and gives in our setting the following result.

\medskip

\begin{prop} \label{ThmPreliminaryMain}
Let $(\bsu_\beta^\varepsilon,\bsw_1^\varepsilon,\bsw_2^\varepsilon)$ be a triplet of modelled distributions in the class \eqref{spaceuvw} solving Equation \eqref{EqModifiedQgKPZzeta} with respect to the BPHZ model ${\sf M}^\varepsilon$ of a growth factor $m>0$. Under Assumption {\it\ref{asmp1}} one can choose $0<T'<T$ and $\varepsilon_0>0$ both small enough for
$$
u^\varepsilon\defeq {\sf R}^\varepsilon\big((\widetilde{\bsu_\beta^\varepsilon})_T\big)
$$ 
to satisfy the bound
\begin{equation}\label{eq:au-avsmall}
\sup_{\varepsilon\in(0,\varepsilon_0)}\sup_{t\in(0,2T')} \big\|a(u^\varepsilon)-a(v)\big\|_{L^\infty(\bbT)} < \frac1m
\end{equation}
and solve the `renormalized' equation
\begin{equation} \label{EqRenormalizedEquationGeneral}
\big(\partial_t - a(u^\varepsilon)\partial_x^2\big) u^\varepsilon = f(u^\varepsilon)\xi^\varepsilon + g(u^\varepsilon)(\partial_x u^\varepsilon)^2 + \sum_{\tau^{\bsp}\in\bbB_-\cap\bbC} \frac{\vcst^\varepsilon_v(\cdot,\tau^{\bsp})}{S(\tau^{\bsp})} \, \frak{F}_v[\tau^{\bsp}]\big(u^\varepsilon,\partial_x u^\varepsilon,\cdot\big)
\end{equation}
on $z\in(0,T')\times\bbT$, with initial condition $u_0$. The last term of \eqref{EqRenormalizedEquationGeneral} has at most linear growth with respect to $\partial_xu^\varepsilon$.
\end{prop}

\medskip

\begin{Dem}
We follow the proof of Theorem \textit{9} of \cite{NoExtendedDecoration}.
Note that $u^\varepsilon={\sf Q}_{\bf1}\bsu_\beta^\varepsilon$ converges in $L^\infty((0,T)\times\bbT)$ as $\varepsilon$ goes to $0$ under Assumption {\it\ref{asmp1}} and the definition of the norm of modelled distributions.
Hence we can choose $0<T'<T$ and $\varepsilon_0>0$ such that \eqref{eq:au-avsmall} holds.

Denote by ${\sf M}^\varepsilon=({\sf g}^\varepsilon,{\sf\Pi}^\varepsilon)$ the given BPHZ model and write $R^\varepsilon=R_{\ell_v^\varepsilon}$ and $S^\varepsilon=S_{\ell_v^\varepsilon}$. Moreover we write the factorization \eqref{eq:PixfactorizationPitimesR} as ${\sf\Pi}_z^\varepsilon=\Pi_z^{\varepsilon,\times}R^\varepsilon$. Since ${\sf M}^\varepsilon$ is a smooth model, the reconstruction of $\bsw^\varepsilon=\bsw_1^\varepsilon+\bsw_2^\varepsilon$ satisfies
$$
\big({\sf R}^\varepsilon(\widetilde{\bsw^\varepsilon})_T\big)(z)
= \big({\sf\Pi}_z^\varepsilon\bsw^\varepsilon(z)\big)(z)
= \big({\Pi}_z^{\varepsilon,\times}R^\varepsilon(z)\bsw^\varepsilon(z)\big)(z)
$$
on $z\in(0,T')\times\bbT$.
By the duality between $R^\varepsilon$ and $S^\varepsilon$ and by Proposition {\it\ref{prop:Bseries}}, we have for any $\tau\in\Basis_-\cup\{{\bf1}\}$,
$$
\langle R^\varepsilon(z)\bsw^\varepsilon(z),\tau\rangle
=\langle\bsw^\varepsilon(z),S^\varepsilon(z)\tau\rangle
=\frak{F}_v[S^\varepsilon(z)\tau]\big(u_{\bf0}^\varepsilon(z),u_{(0,1)}^\varepsilon(z),u_{(0,2)}^\varepsilon(z),z\big),
$$
where $u_{\bf k}^\varepsilon$ denotes the $X^{\bf k}$-component of $\bsu_{\beta+\alpha_0}^\varepsilon$.
In other words, we have the expansion
$$
R^\varepsilon(z)\bsw^\varepsilon(z)=\sum_{\tau\in\Basis_-\cup\{{\bf1}\}}
\frac1{S(\tau)}\frak{F}_v[S^\varepsilon(z)\tau]\big(u_{\bf0}^\varepsilon(z),u_{(0,1)}^\varepsilon(z),u_{(0,2)}^\varepsilon(z),z\big)\tau+\cdots,
$$
where the omitted term $(\cdots)$ is spanned by $\tau\in\Basis$ with $|\tau|>0$, so it disappears after acting $\Pi_z^{\varepsilon,\times}(\cdot)(z)$.
By using Proposition {\it\ref{PropMorphismPropertyRStar}} and by tracing back the proof of Proposition {\it\ref{prop:Bseries}}, we can reorganize the first term into
$$
R^\varepsilon(z)\bsw^\varepsilon(z)
=\sum_{\zeta\in\{\Xi,{\bf1}\}}\frak{F}_v[S^\varepsilon(z)\zeta](\bsu_\beta^\varepsilon(z),\bsD\bsu_\beta^\varepsilon(z),\bsD^2\bsu_{\beta+\alpha_0}^\varepsilon(z),z)\zeta+\cdots,
$$
where smooth functions $\frak{F}_v[S^\varepsilon(z)\zeta]({\sf u}_{\bf0},{\sf u}_{(0,1)},{\sf u}_{(0,2)},z)$ are extended into the functionals on modelled distributions in the sense of Proposition {\it\ref{prop:productMD}} (2).
Note that $\bsu_\beta^\varepsilon,\bsD\bsu_\beta^\varepsilon$, and $\bsD^2\bsu_{\beta+\alpha_0}^\varepsilon$ take values in $\bbU_{\bf0},\bbU_{(0,1)}$, and $\bbU_{(0,2)}$, respectively, and $R^\varepsilon$ preserves any element of $\bbU_{\bf n}$.
Therefore, by Theorem {\it\ref{thm:multilevelSchauder}}, we have
\begin{align*}
&\big(\Pi_z^{\varepsilon,\times}\bsu_\beta^\varepsilon(z)\big)(z)
= \big({\sf\Pi}_z^{\varepsilon}\bsu_\beta^\varepsilon(z)\big)(z)
= {\sf R}^\varepsilon(\widetilde{(\bsu_\beta^\varepsilon)}_T)(z) = u^\varepsilon(z),   \\
&\big(\Pi_z^{\varepsilon,\times}\bsD\bsu_\beta^\varepsilon(z)\big)(z) = \partial_xu^\varepsilon(z),   \\
&\big(\Pi_z^{\varepsilon,\times}\bsD^2\bsu_{\beta+\alpha_0}^\varepsilon(z)\big)(z) = \partial_x^2u^\varepsilon(z),
\end{align*}
on $z\in(0,T')\times\bbT$. Since $\Pi_z^{\varepsilon,\times}$ is multiplicative we have
\begin{align*}
(\partial_t-L^v+c)u^\varepsilon=\sum_{\zeta\in\{\Xi,{\bf1}\}}\frak{F}_v[S^\varepsilon(z)\zeta]\big(u^\varepsilon(z),\partial_xu^\varepsilon(z),\partial_x^2u^\varepsilon(z),z\big) \, {\sf\Pi}^\varepsilon\zeta(z),
\end{align*}
or equivalently,
\begin{align*}
\big(\partial_t - a(u^\varepsilon)\partial_x^2\big) u^\varepsilon = f(u^\varepsilon)\xi^\varepsilon + g(u^\varepsilon)(\partial_x u^\varepsilon)^2 
+
\frak{F}_v\big[(S_v^\varepsilon(z)-\id){\bf1}\big]\big(u^\varepsilon(z),\partial_xu^\varepsilon(z),\partial_x^2u^\varepsilon(z),z\big),
\end{align*}
for any $z\in(0,T')\times\bbT$. By the definition of $S^\varepsilon(z)$ and the fact that $\ell^\varepsilon$ vanishes on $\Basis_-\setminus\bbC$, we obtain the right hand side of \eqref{EqRenormalizedEquationGeneral}.

It remains to show that the most right term does not depend on $\partial_x^2u^\varepsilon$ and is at most linear with respect to $\partial_xu^\varepsilon$.
If the function $\frak{F}_v[\tau]$ depends on ${\sf u}_{(0,2)}$, then $\tau$ has one node $v$ with type $\frak{t}(v)={\bf1}$ and only $\mcI_{\bf0}^p$ type edges leave from $v$. But such $\tau$ cannot have negative homogeneity.
If some function $\frak{F}_v[\tau]$ is quadratic with respect to ${\sf u}_{(0,1)}$, then $\tau$ has at least two nodes with type ${\bf1}$ from which exactly one edge $\mcI_{(0,1)}$ leaves. 
Similarly to the proof of Proposition {\it\ref{prop:Rellcomplete}}, such $\tau$ cannot have negative homogeneity.
\end{Dem}

\medskip

Next we reduce the $v$-dependence of the counterterm of \eqref{EqRenormalizedEquationGeneral}. Since we have proved Proposition {\it\ref{prop:Fvcontinuous}}, we consider only the trees of $\bbF_-$. Recall that each tree cannot have the node decoration greater than $(0,1)$ to have the negative homogeneity. We define $\bbF_-'\subset\bbF$ as the subset of all trees $\tau\in\bbF_-$ having one node $p\in N_\tau$ such that $\tau_{\ge p}=X_2\mcI_{(0,2)}(\sigma)$ for some $\sigma\in\bbF_-$.

\medskip

\begin{lem}
For any $\tau\in(\bbF_-\cap\bbC)\setminus\bbF_-'$ the function $\frak{F}_v[\tau]$ depends only on $({\sf u}_{\bf0},{\sf u}_{(0,1)})$. We write $\frak{F}[\tau]$.
\end{lem}

\medskip

\begin{Dem}
In the factorization \eqref{EqDefnFrakFTau2}, a function depending on $v$ can appear only at a node $p\in N_\tau$ from which an edge with type $\mcI_{(0,2)}^0$ leaves. By the rule ({\bf C}3) the function $\frak{F}_v[p;\tau]$ does not depend on $v$ unless $\tau_{\ge p}=X_2\mcI_{(0,2)}(\sigma)$ for some $\sigma\in\bbF_-$.
\end{Dem}

\medskip

For a positive parameter $\lambda$ we denote by
$$
Z_t^\lambda(x) = Z^\lambda(t,x) \defeq {\bf 1}_{t>0}\frac{e^{-ct}}{\sqrt{4\pi \lambda t}}\exp\bigg(-\frac{|x|^2}{4\lambda t}\bigg),
$$
the fundamental solution built from the constant coefficient parabolic operator $\partial_t-\lambda\partial_x^2+c$. The naive admissible model associated with $Z^\lambda$ and the smooth noise $\xi^\varepsilon$ is the unique multiplicative model such that 
$$
{\sf \Pi}^\varepsilon_\lambda\Xi=\xi^\varepsilon,\qquad
\big({\sf \Pi}^\varepsilon_\lambda X^{\bf k}\big)(z)=z^{\bf k},\qquad
{\sf \Pi}^\varepsilon_\lambda\big(\mcI_{\bf n}^p(\tau)\big) = \Big(\partial_z^{\bf n}Z^\lambda* (\partial_x^2Z^\lambda)^{*p}\Big)*{\sf \Pi}^\varepsilon_\lambda\tau,
$$
where $*$ denotes the spacetime convolution.
Note that the convolution with $\partial_x^2Z^\lambda$ preserves the space $\mcC_\mfs^{0+}(\bbR\times\bbT)$ -- see Theorem {\it\ref{thm:InverseParabolic}}.
The \textit{BHZ character} $\fcst^\epsilon_\lambda(\cdot)$ on $\bbB_-$ is defined in that setting as 
\begin{equation} \label{EqDefnCharacter} \begin{split}
&\fcst^\varepsilon_\lambda(\tau) \defeq h^\varepsilon_\lambda(S_-'\tau),   \\
&h^\varepsilon_\lambda(\tau) \defeq \bbE\big[{\sf \Pi}^\varepsilon_\lambda\tau(0)\big],
\end{split} \end{equation}
where $S'_- : \spa(\Basis_-)\rightarrow \bbR[\Basis]$ (the symmetric algebra generated by $\Basis$) is the natural extension to our setting of the negative twisted antipode -- see Proposition {\it 6.6} in \cite{BHZ} or Section 7 of \cite{RSGuide} for its definition in the usual BHZ setting. 

\medskip

\begin{asmp}\label{asmp2}
Let $\ell_v^\varepsilon$ be the renormalization character of a growth factor $m_0>0$ satisfying Assumption \ref{asmp1}.
For any $\tau\in\bbF_-$ there exists an $\varepsilon$-independent constant $C(\tau)$ such that
\begin{equation*}
\big|\vcst_v^\varepsilon(z,\tau^{\bsp})-\fcst^\varepsilon_{a(v(z))}(\tau^{\bsp})\big|\le C(\tau) \, m_0^{|\bsp|}
\end{equation*}
for any $\bsp:E_\tau\to\bbN$ and $z\in\bbR\times\bbT$. Moreover the map $l_\lambda^\varepsilon(\cdot)$ vanishes on $\bbF_-'$.
\end{asmp}

\medskip

We check that \textbf{\textit{Assumption 2}} is satisfied by some examples in Section {\sf\ref{SubsectionAsmp2}}. The next statement is the core fact to get the renormalized equation under the form \eqref{EqRenormalizedEquation} stated in Theorem \textit{\ref{ThmMainRenormalizedSimplified}}. Recall from Lemma \textit{\ref{NonplanarvsPlanar}} the canonical projection map $\mcP$ from planar trees to non-planar trees.

\medskip

\begin{lem} \label{LemCrucial}
For any $\tau\in \bbF_-$ the function 
$$
\lambda\mapsto\fcst^\varepsilon_\lambda(\tau)
$$
is analytic in any given bounded interval of $\textbf{\textsf{R}}$ whose closure does not contain the point $0$, and for a fixed planar tree $\widetilde{\tau}$ such that $\mcP(\widetilde{\tau})=\tau$ one has
$$
\frac1{n!}\partial_\lambda^n \fcst^\varepsilon_\lambda(\tau) = \sum_{\bsp\in \bbN^{E_{\widetilde{\tau}}},\, \vert\bsp\vert=n} \fcst_\lambda^\varepsilon(\widetilde{\tau}^{\bsp}).
$$
\end{lem}

\medskip

\begin{Dem}
By an elementary computation we have
\begin{align}\label{lambdaderivative}
\partial_\lambda Z^\lambda(t,x)=t\partial_x^2 Z^\lambda(t,x)=\int_0^t\int_{\bbR}Z^\lambda(t-s,x-y)\partial_x^2 Z^\lambda(s,y)\,dy\, ds.
\end{align}
By induction we have 
$$
\partial_\lambda^p(\partial_z^{\bf n}Z)=p!(\partial^{\bf n}Z^{\lambda})*(\partial_x^2 Z^\lambda)^{*p}.
$$
Therefore for any fixed $\bsp=(\bsp(e))_{e\in E_{\widetilde\tau}}\in\bbN^{E_{\widetilde\tau}}$, applying $\partial_\lambda^{\bsp(e)}$ simultaneously to the kernel associated with each edge $e\in E_{\widetilde\tau}$ transforms $\fcst^\varepsilon_\lambda(\tau)$ into $\bsp! \, \fcst_\lambda^\varepsilon(\widetilde{\tau}^{\bsp})$, where 
$$
\bsp!\defeq\prod_{e\in E_{\widetilde\tau}}\bsp(e)!.
$$
By the Leibniz rule, the same term appears $n!/\bsp!$ times in the expansion of $\partial_\lambda^n \fcst^\varepsilon_\lambda(\tau)$.
\end{Dem}

\medskip

\noindent \textit{Proof of Theorem \ref{ThmMainRenormalizedSimplified} --} 
For each $\tau\in(\bbF_-\cap\bbC)\setminus\bbF_-'$ we fix a planar tree $\widetilde{\tau}$ such that $\mcP(\widetilde{\tau})=\tau$. It follows from Lemma \textit{\ref{NonplanarvsPlanar}}, Proposition {\it\ref{prop:Fvcontinuous}}, and Lemma \textit{\ref{LemCrucial}} that the counterterm in the renormalized Equation \eqref{EqRenormalizedEquationGeneral} equals to the following simple form up to an $\varepsilon$-uniform remainder term 
\begin{align*} 
\sum_{\tau\in(\bbF_-\cap\bbC)\setminus\bbF_-',\,\bsp:E_\tau\to\bbN} &\frac{\fcst^\varepsilon_{a(v(\cdot))}(\tau^{\bsp})}{S(\tau^{\bsp})} \, \frak{F}_v[\tau^{\bsp}](u^\varepsilon, \partial_xu^\varepsilon)   \\
&\hspace{-0.3cm}= \sum_{\tau\in(\bbF_-\cap\bbC)\setminus\bbF_-'}\frac1{S(\tau)} 
\sum_{\bsp\in\bbN^{E_{\widetilde{\tau}}}} \fcst^\varepsilon_{a(v(\cdot))}(\widetilde{\tau}^{\bsp})\, \frak{F}_v[\widetilde{\tau}^{\bsp}](u^\varepsilon, \partial_xu^\varepsilon)    \\
&\hspace{-0.3cm}= \sum_{\tau\in(\bbF_-\cap\bbC)\setminus\bbF_-'}\frac1{S(\tau)} \, \frak{F}[\tau](u^\varepsilon, \partial_xu^\varepsilon)\sum_{n=0}^\infty\big(a(u^\varepsilon)-a(v)\big)^n
\sum_{\bsp\in\bbN^{E_{\widetilde{\tau}}},\,|\bsp|=n} \fcst^\varepsilon_{a(v(\cdot))}(\widetilde{\tau}^{\bsp})\\
&\hspace{-0.81cm}\overset{\textrm{Lemma {\it\ref{LemCrucial}}}}{=} \sum_{\tau\in(\bbF_-\cap\bbC)\setminus\bbF_-'} \frac1{S(\tau)} \, \frak{F}[\tau](u^\varepsilon, \partial_xu^\varepsilon) \,\fcst^\varepsilon_{a(u^\varepsilon(\cdot))}(\tau).
\end{align*}
After the third line, we write $\frak{F}[\tau](u^\varepsilon, \partial_xu^\varepsilon)=\frak{F}_v[\tau](u^\varepsilon, \partial_xu^\varepsilon)$ for each $\tau\in(\bbF_-\cap\bbC)\setminus\bbF_-'$ because it does not depend on $v$. This completes the proof of Theorem {\it \ref{ThmMainRenormalizedSimplified}}.   \hfill $\rhd$

\medskip

We finish this section by showing that the diverging term $\fcst^\varepsilon_{a(u^\varepsilon(\cdot))}$ in the counterterm takes a particularly nice form under the condition that the noise is Gaussian and regularized only in the spatial variable by symmetric mollifiers. To avoid situations with time regularization we consider here only a spatial noise or a spacetime noise that is white in time with $f=1$. Recall that $|\tau|_{\Xi}$ denotes the number of $\Xi$-type nodes that appear in $\tau$. Recall also from \eqref{EqDefnCharacter} the definition of the function $h^\varepsilon_\lambda$.

\medskip

\begin{prop} \label{lem:scaling}
Assume that $\xi$ is a stationary centered Gaussian noise and define 
$$
\overline{\xi}^\varepsilon(t,x)=\big(\xi(t,\cdot)*\rho_\varepsilon\big)(x)
$$ 
with an even mollifier $\rho_\varepsilon$ and the spatial convolution operator $*$. Then $h^\varepsilon_\lambda(\tau)=0$ if $|\tau|_\Xi$ is odd, otherwise
$$
h^\varepsilon_\lambda(\tau)=
\begin{cases}
\lambda^{-\sharp N_\tau+1}h^\varepsilon_1(\tau), & \text{if $\xi(x)$ depends on only space},   \\
\lambda^{|\tau|_{\Xi}/2-\sharp N_\tau+1}h^\varepsilon_1(\tau), & \text{if $\xi(t,x)$ is white in time}.
\end{cases}
$$
\end{prop}

\medskip

\begin{Dem}
The former statement holds because $\overline{\xi}^\varepsilon$ is centered Gaussian. Let $|\tau^{\bf 0}|_{\Xi}=2a$ be an even number and let $b$ be the number of ${\bf1}$-type nodes. If the root of $\tau$ is not a $\Xi$-type node, then the expectation of $\big({\sf \Pi}^\varepsilon_\lambda\tau\big)(0)$ is given by an integral of the form
$$
\int C^\varepsilon(z^1-z^2)\cdots C^\varepsilon(z^{2a-1}-z^{2a}) (z^1)^{n_1}\cdots(z^{2a})^{n_{2a}}A^\lambda\big(z^1,\dots,z^{2a}\big) \, dz^1\cdots dz^{2a},
$$
where 
$$
C^\varepsilon(z) \defeq \mathbb{E}\big[\,\overline{\xi}^\varepsilon(z) \overline{\xi}^\varepsilon(0)\big]
$$ 
and 
$$
A^\lambda(z^1,\dots,z^{2a})=\int \widetilde{A}^\lambda\big(z^1,\dots,z^{2a},w^1,\dots,w^{b-1}\big) \, dw^1\cdots dw^{b-1}
$$
with a product $\widetilde{A}$ of polynomials $(w^i)^{m_i}$ and kernels $\partial^{n_{ij}}Z^\lambda$. Because of the form of equation and the restriction $\alpha_0\in(0,1)$, the $\frak{n}$-decorations $m_i$ and $\frak{e}$-decorations $n_{ij}$ are $\bf 0$, $(0,1)$, or $(0,2)$. 
So, they are independent of the change of variables $z\mapsto z_\lambda\defeq(\lambda t,x)$. Hence
\begin{align*}
A^\lambda\big(z^1,\dots,z^{2a}\big) = \lambda^{-b+1}A^1\big(z_\lambda^1,\dots,z_\lambda^{2a}\big)
\end{align*}
by a scaling argument. If $\xi(x)$ is $t$-independent then we have
$$
\big({\sf\Pi}^\varepsilon_\lambda\tau\big)(0) = \lambda^{-2a-b+1} \big({\sf\Pi}^\varepsilon_1\tau\big)(0),
$$
since $C^\varepsilon(x)$ does not depend on time. 
Next let $\xi(t,x)$ be white in time. 
Then $C^\varepsilon(z)$ is of the form $\delta(t)\bar{C}^\varepsilon(x)$ for some function $\bar{C}^\varepsilon$. Since $\delta(t)$ reduces the number of time components $t^1,\dots,t^{2a}$ of $z^1,\dots,z^{2a}$ by a half, we have
$$
\big({\sf\Pi}^\varepsilon_\lambda\tau\big)(0) = \lambda^{-a-b+1} \big({\sf\Pi}^\varepsilon_1\tau\big)(0).
$$
We can perform similar computations when the root of $\tau$ is a $\Xi$-type node.
\end{Dem}

\medskip

In the setting of Proposition \textit{\ref{lem:scaling}} the counterterm is of the form
\begin{equation}\label{SimpleFormWhenWhite}
\sum_{\tau\in(\bbF_-\cap\bbC)\setminus\bbF_-'} \frac1{S(\tau)}\bigg(\sum_\sigma s(\tau,\sigma_1,\dots,\sigma_n)\frac{h_1^\varepsilon(\sigma_1)\cdots h_1^\varepsilon(\sigma_n)}{\lambda^{\theta(\sigma_1)+\cdots+\theta(\sigma_n)}}\bigg)  \frak{F}[\tau](u^\varepsilon,\partial_x u^\varepsilon),
\end{equation}
where $S_-'\tau=\sum_{\sigma_1,\dots,\sigma_n\in\bbB_-} s(\tau,\sigma_1,\dots,\sigma_n) \, \sigma_1\cdots\sigma_n$ and the exponent $\theta(\sigma)$ is given in the statement of Proposition \textit{\ref{lem:scaling}}.

\bigskip

\subsection{Examples satisfying Assumption \ref{asmp1}{\boldmath $.$} \hspace{0.15cm}}
\label{SubsectionAsmp1}

In this and the next section we assume that the regularized noise $\xi^\varepsilon$ is given by the spacetime convolution
$$
\xi^\varepsilon(t,x) = \int_{\bbR^2}\rho^\varepsilon(t-s,x-y)\xi(s,y)dsdy
$$
with some smooth mollifier $\rho^\varepsilon$ on $\bbR^2$, where $\xi$ is regarded as a spatially periodic distribution. Then the covariance function 
$$
C^\varepsilon(z,w) \defeq \bbE\big[ \xi^\varepsilon(z) \xi^\varepsilon(w) \big]
$$ 
is homogeneous in the sense that 
$$
C^\varepsilon(z,w) = C^\varepsilon(z-w,0).
$$ 
For simplicity we write 
$$
C^\varepsilon(z) \defeq C^\varepsilon(z,0).
$$

\ssk

In this section we discuss Assumption {\it\ref{asmp1}} for some examples when $\xi$ is the spacetime white noise. As for the semilinear equation of the form
$$
\big(\partial_t - \partial_x^2 \big) u = f(u)\xi + g(u)(\partial_x u)^2,
$$
we can apply the general convergence result of \cite{CH, HS, RandomModel}. Then we can claim that, if $|\tau|>-3/2$ for all decorated trees $\tau$ except $\Xi$ which strongly conform to the rule associated with the above equation, then the associated BPHZ models ${\sf M}^\varepsilon$ converge as the approximation parameter $\varepsilon$ goes to $0$. If we assume that the spacetime regularity of the noise $\xi$ is $\alpha_0-2$, this sufficient condition for the convergence of BPHZ models is equivalent to
$$
\alpha_0>\frac14.
$$
In particular the spacetime white noise case $\alpha_0=1/2-\kappa$ with small $\kappa>0$ is included.

The proofs of \cite{HS, RandomModel} are inductive and based on Poincar\'e-type inequality of noise functionals. We cannot apply these proofs directly to the quasilinear equations because our model space is infinite dimensional and the integral operator $K^{v}$ is inhomogeneous. However we expect that the first problem is solved by the fact that, for each tree $\tau$, the smallest sector containing $\tau$ is finite dimensional. Also, as far as the second problem is concerned, the homogeneity of the integral kernel is used only to obtain the $\varepsilon$-uniform boundedness of the expectations
\begin{equation}\label{eq:proof:asmp1exp}
\big|\bbE\big[\mcQ_\theta({\sf\Pi}_z^\varepsilon\tau^{\bsp})(z)\big]\big|
\lesssim m^{|\bsp|}\theta^{|\tau|/4}
\end{equation}
for any $\tau\in\bbF_-$ with $|\tau|>-3/2$ and their convergences as $\varepsilon$ goes to $0$ -- see Section 3.2 of \cite{RandomModel}. Therefore it is expected that Assumption {\it\ref{asmp1}} is reduced to the convergence of the distributions
\begin{equation}\label{eq:proof:asmp1noise}
{\sf\Pi}^\varepsilon\Xi=\xi^\varepsilon,\qquad
{\sf\Pi}^\varepsilon\mcI_{(0,2)}^p\Xi=
(\partial_x^2 K^{v})^{(p+1)} \xi^\varepsilon
\end{equation}
with the minimum homogeneity $\alpha_0-2$ and the estimate \eqref{eq:proof:asmp1exp}. It should be noted that the expectation \eqref{eq:proof:asmp1exp} vanishes for any tree with an odd number of $\Xi$ symbols.

\medskip

We now prove \eqref{eq:proof:asmp1noise} and give two examples of trees where \eqref{eq:proof:asmp1exp} holds true. For the proof we need some estimate on the integral kernels
$$
\mcQ^{{\bf n},p}(z,w)\defeq\int_0^1\mcQ_\theta^{{\bf n},p}(z,w)d\theta,\qquad
\mcQ_\theta^{{\bf n},p}(z,w)\defeq(\partial_z^{\bf n}\mcQ_\theta^v)*(\partial_x^2K^v)^{*p}(z,w),
$$
for any ${\bf n}\in\bbN^2$ such that $|{\bf n}|_\mfs\le4$, where $*$ denotes the spacetime convolution
$$
(A*B)(z,w) \defeq \int_{\bbR^2}A(z,u)B(u,w)du.
$$ 
Recall the Gaussian type function ${\sf G}_\theta^{(c)}$ defined by \eqref{EqDefnGaussianLike}. We omit the letter $(c)$ if it does not need to be emphasized. As shown in Theorem {\it\ref{thm:GaussGamma}} in Appendix {\sf\ref{SubsectionFundamental}}, there exist some constant $c_0,C_0>0$ and $\delta_0\in(0,\alpha)$ such that
\begin{align}
\label{prop:K"kernelestimate1}
\big|\partial_z^{\bf n}\mcQ_\theta^v(z,w)\big|
&\le C_0 \, \theta^{-|{\bf n}|_\mfs/4} \, {\sf G}_\theta^{(c_0)}(z-w),   \\
\label{prop:K"kernelestimate2}
\big|\partial_z^{\bf n}\mcQ_\theta^v(z+h{\bf e}_i,w)-\partial_z^{\bf n}\mcQ_\theta^v(z,w)\big|
&\le C_0 \, \theta^{-(|{\bf n}|_\mfs+\beta)/4}|h|^{\beta/\mfs_i} \Big\{ {\sf G}_\theta^{(c_0)}(z+h{\bf e}_i-w) + {\sf G}_\theta^{(c_0)}(z-w)\Big\},   \\
\label{prop:K"kernelestimate3}
\big|\partial_z^{\bf n}\mcQ_\theta^v(z,w')-\partial_z^{\bf n}\mcQ_\theta^v(z,w)\big|
&\le C_0 \, \theta^{-(|{\bf n}|_\mfs+\delta_0)/4}\|w'-w\|_\mfs^{\delta_0} \Big\{ {\sf G}_\theta^{(c_0)}(z-w') + {\sf G}_\theta^{(c_0)}(z-w) \Big\}
\end{align}
for any $\theta\in(0,1]$, $z,w,w'\in\bbR^2$, ${\bf n}\in\bbN^2$ such that $|{\bf n}|_\mfs\le4$, $h\in\bbR$, $i\in\{1,2\}$, and $\beta\in[0,\mfs_i]$ is such that $\beta\le4-|{\bf n}|_\mfs+\delta_0$.

\medskip

\begin{prop}\label{prop:Lpkernelestimate}
For any $\delta\in(0,\delta_0)$, there exist constants $c\in(0,c_0)$ and $C\in(0,C_0)$ such that,
one has
\begin{align}
\label{prop:Lpkernelestimate1}
\big|\mcQ_\theta^{{\bf n},p}(z,w)\big|
&\le C^p \, \theta^{-(|{\bf n}|_\mfs+\delta)/4} \, {\sf G}_\theta^{(c)}(z-w),   \\
\label{prop:Lpkernelestimate2}
\big|\mcQ_\theta^{{\bf n},p}(z+h{\bf e}_i,w)-\mcQ_\theta^{{\bf n},p}(z,w)\big|
&\le C^p \, \theta^{-(|{\bf n}|_\mfs+\delta+\beta)/4} \, |h|^{\beta/\mfs_i} \Big\{{\sf G}_\theta^{(c)}(z+h{\bf e}_i-w) + {\sf G}_\theta^{(c)}(z-w) \Big\},   \\
\label{prop:Lpkernelestimate3}
\big|\mcQ_\theta^{{\bf n},p}(z,w')-\mcQ_\theta^{{\bf n},p}(z,w)\big|
&\le C^p \, \theta^{-(|{\bf n}|_\mfs+\delta)/4} \, \|w'-w\|_\mfs^{\delta} \, \Big\{ {\sf G}_\theta^{(c)}(z-w') + {\sf G}_\theta^{(c)}(z-w) \Big\}
\end{align}
for any $p\ge0$, $\theta\in(0,1]$, $z,w,w'\in\bbR^2$, ${\bf n}\in\bbN^2$ such that $|{\bf n}|_\mfs\le4$, $h\in\bbR$, $i\in\{1,2\}$, and $\beta\in[0,\mfs_i]$ is such that $\beta\le4-|{\bf n}|_\mfs+\delta$.
\end{prop}

\medskip

\begin{Dem}
For the initial case $p=0$, all estimates follow from \eqref{prop:K"kernelestimate1}, \eqref{prop:K"kernelestimate2}, and \eqref{prop:K"kernelestimate3}.
Suppose that $\mcQ_\theta^{{\bf n},p}$ satisfy both \eqref{prop:Lpkernelestimate1} and \eqref{prop:Lpkernelestimate3} and consider the integral
$$
\mcQ_{\theta,\theta'}^{{\bf n},p+1}(z,w) \defeq \int_{\bbR^2}\mcQ_\theta^{{\bf n},p}(z,u)\partial_{u_2}^2K_{\theta'}^v(u,w)du.
$$
Since $\int_{\bbR^2}\partial_{u_2}^2K_{\theta'}^v(u,w)du=0$, by using \eqref{prop:Lpkernelestimate3} we have
\begin{align*}
\big|\mcQ_{\theta,\theta'}^{{\bf n},p+1}&(z,w)\big|   \\
&\leq \int_{\bbR^2} \big| \mcQ_{\theta}^{{\bf n},p}(z,u)-\mcQ_{\theta}^{{\bf n},p}(z,w)\big| \, |\partial_{u_2}^2K_{\theta'}^v(u,w)| \, du   \\
&\leq C^p \, C_0 \, \theta^{-(|{\bf n}|_\mfs+\delta)/4}(\theta')^{-1}\int_{\bbR^2}
\Big\{{\sf G}_\theta^{(c)}(z-u) + {\sf G}_\theta^{(c)}(z-w)\Big\} \, \|u-w\|_\mfs^\delta \, {\sf G}_{\theta'}^{(c_0)}(u-w)du.
\end{align*}
By using Lemma {\it\ref{lem:ConvolutionGauss}}-{\it\ref{lem:ConvolutionGauss1}} and {\it\ref{lem:ConvolutionGauss3}}, for some $p$-independent constant $C'>0$ we have
\begin{align}\label{prop:Lpkernelestimateproof}
\begin{aligned}
|\mcQ_{\theta,\theta'}^{{\bf n},p+1}(z,w)| &\le C^p \, C'\theta^{-(|{\bf n}|_\mfs+\delta)/4} \, (\theta')^{-1+\delta/4}
\Big\{(\theta+\theta')^{3/4}\theta^{-3/4}{\sf G}_{\theta+\theta'}^{(c)}(z-w) + {\sf G}_\theta^{(c)}(z-w)\Big\}   \\
&\le 2C^p \, C'\theta^{-(|{\bf n}|_\mfs+\delta)/4} \, (\theta')^{-1+\delta/4} \, {\sf G}_{\theta}^{(c)}(z-w).
\end{aligned}
\end{align}
By integrating this inequality over $\theta'\in(0,1]$ we have
$$
\big|\mcQ_\theta^{{\bf n},p+1}(z,w)\big| \leq \int_0^1|\mcQ_{\theta,\theta'}^{{\bf n},p+1}(z,w)| \, d\theta'
\leq C^{p+1} \, \theta^{-(|{\bf n}|_\mfs+\delta)/4} \, {\sf G}_\theta^{(c)}(z-w)
$$
for some choice of $C$ independent of $p$. To prove \eqref{prop:Lpkernelestimate3}, by interpolating the inequality \eqref{prop:Lpkernelestimateproof} and
\begin{align*}
\big| \mcQ_{\theta,\theta'}^{{\bf n},p+1}&(z,w') - \mcQ_{\theta,\theta'}^{{\bf n},p+1}(z,w) \big|   \\
&\hspace{-0.1cm}\leq \int_{\bbR^2}|\mcQ_\theta^{{\bf n},p}(z,u)| \, \big|\partial_{u_2}^2K_{\theta'}^v(u,w')-\partial_{u_2}^2K_{\theta'}^v(u,w)\big| \, du   \\
&\hspace{-0.1cm}\leq C^pC_0\theta^{-|{\bf n}|_\mfs/4}(\theta')^{-1-\delta_0/4}\|w'-w\|_\mfs^{\delta_0}\int_{\bbR^2}
{\sf G}_\theta^{(c)}(z-u) \Big\{{\sf G}_{\theta'}^{(c_0)}(u-w') + {\sf G}_{\theta'}^{(c_0)}(u-w)\Big\} \, du   \\
&\hspace{-0.1cm}\leq C^pC'\theta^{-|{\bf n}|_\mfs/4}(\theta')^{-1-\delta_0/4}\|w'-w\|_\mfs^{\delta_0}
\Big\{{\sf G}_{\theta}^{(c)}(z-w')+{\sf G}_{\theta}^{(c)}(z-w)\Big\},
\end{align*}
we have
\begin{align*}
\big|\mcQ_{\theta,\theta'}^{{\bf n},p+1}(z,w') &- \mcQ_{\theta,\theta'}^{{\bf n},p+1}(z,w) \big|   \\
&\le C^pC' \, \theta^{-(|{\bf n}|_\mfs+\delta')/4}(\theta')^{-1+\delta''/4} \, \|w'-w\|_\mfs^{\delta} \,
\Big\{{\sf G}_{\theta}^{(c)}(z-w')+{\sf G}_{\theta}^{(c)}(z-w) \Big\}
\end{align*}
for some $\delta'\in(0,\delta)$ and $\delta''>0$. By integrating this inequality over $\theta'\in(0,1]$ we obtain \eqref{prop:Lpkernelestimate3} since $\theta^{-(|{\bf n}|_\mfs+\delta')/4}\le\theta^{-(|{\bf n}|_\mfs+\delta)/4}$.
The proof of \eqref{prop:Lpkernelestimate2} is similar.
\end{Dem}

\ssk

We first prove the convergence of distributions \eqref{eq:proof:asmp1noise} in the space $\mcC_\mfs^{\alpha_0-2,v}(\bbR\times\bbT)$ for any $\alpha_0<1/2$.
By the usual argument based on the hypercontractivity as in Section 10 of \cite{Hai14}, the proof is reduced to the bounds
\begin{equation}\label{eq:proof:asmp1noise:L2}
\bbE\Big[\mcQ_\theta\big((\partial_x^2 K^{v})^{p} \xi\big)(z)^2\Big]
\lesssim m^{2p}\theta^{-(3+\kappa)/4}
\end{equation}
for some constant $m>0$ and any small $\kappa>0$. (We actually need to slightly modify the norm of $\mcC_\mfs^{\alpha_0-2,v}(\bbR\times\bbT)$ to allow the polynomial divergence of distributions as $t$ goes to $\infty$.) The expectation in the left hand side is equal to
$$
\int_{\bbR\times\bbT}\widetilde{\mcQ_\theta^{{\bf 0},p}}(z,w)^2dw,
$$
where we recall that
$$
\widetilde{f}(z,w) = \sum_{k\in\bbZ}f(z+k{\bf e}_2,w)=\sum_{k\in\bbZ}f(z,w+k{\bf e}_2)
$$ 
denotes the spatial periodization of $f$. By the estimate \eqref{prop:Lpkernelestimate1} and Lemma {\it\ref{lem:ConvolutionGauss}}-{\it\ref{lem:ConvolutionGauss3}} the above integral is estimated as
\begin{align*}
\bigg|\int_{\bbR^2}\mcQ_\theta^{{\bf 0},p}(z,w)\widetilde{\mcQ_\theta^{{\bf 0},p}}(z,w)dw\bigg|
&\lesssim\sum_{k\in\bbZ}C^{2p}\theta^{-\kappa/4}\int_{\bbR^2}{\sf G}_\theta(z-w) \, {\sf G}_\theta(z-w+k{\bf e}_2)\,dw   \\
&\lesssim C^{2p}\theta^{-\kappa/4}\sum_{k\in\bbZ}{\sf G}_{2\theta}(k{\bf e}_2)
\lesssim C^{2p}\theta^{-(3+\kappa)/4}
\end{align*}
for any small $\kappa>0$.
Thus we have \eqref{eq:proof:asmp1noise:L2} with $m=C$.

\ssk

As for the proof of \eqref{eq:proof:asmp1exp}, we consider two trees
$$
\tau
=\mcI_{\bf0}^0(\Xi)\mcI_{(0,2)}^0(\Xi)=\,\rtdb{1}{\pol{2}{1}{120}\pol{3}{1}{60}\drbl{1}{2}\drll{1}{3}\drb{1}\drc{2,3}},\qquad
\sigma
=\mcI_{\bf0}^0(\Xi)^3\mcI_{(0,2)}^0(\Xi)=\rtdb{1}{\pol{2}{1}{144}\pol{3}{1}{108}\pol{4}{1}{72}\pol{5}{1}{36}\drbl{1}{2,3,4}\drll{1}{5}\drb{1}\drc{2,3,4,5}}
$$
as examples, where the thick line denotes the operator $\mcI_{\bf0}^0$ and the double line denotes $\mcI_{(0,2)}^0$. The noise symbol $\Xi$ is denoted by a white circle, while ${\bf1}$ is denoted by a black dot.
We does not consider the $\bsp$-decoration because we can obtain the estimates of $\bsp$-decorated versions in the same ways by replacing $K_\theta^v$ below by $K_\theta^v*(\partial_x^2K^v)^{*p}$ and using Proposition {\it\ref{prop:Lpkernelestimate}}.

\ssk

For $\tau$, the renormalization character is
$$
\ell_v^\varepsilon(z,\tau)
=-\int_{(\bbR\times\bbT)^2}\widetilde{K^v}(z,w_1)
\widetilde{\partial_x^2K^v}(z,w_2)C^\varepsilon(w_1-w_2)dw_1dw_2,
$$
and the action of the BPHZ model is
\begin{align*}
{\sf\Pi}_z^{\varepsilon}\tau
&=\int_{(\bbR\times\bbT)^2}\big(\widetilde{K^v}(\cdot,w_1)-\widetilde{K^v}(z,w_1)\big)
\widetilde{\partial_x^2K^v}(\cdot,w_2)\xi^\varepsilon(w_1)\xi^\varepsilon(w_2)dw_1dw_2+\ell_v^\varepsilon(\cdot,\tau).
\end{align*}
Hence its expectation is
\begin{align*}
\bbE[{\sf\Pi}_z^{\varepsilon}\tau]
=-\int_{(\bbR\times\bbT)^2}\widetilde{K^v}(z,w_1)\widetilde{\partial_x^2K^v}(\cdot,w_2)
C^\varepsilon(w_1-w_2)dw_1dw_2.
\end{align*}
For simplicity we consider the limit $\varepsilon\to0$ and replace $C^\varepsilon$ by the Dirac's delta function.
Then the above integral is written as $\sum_{k\in\bbZ}\int_0^1I_\theta(z+k{\bf e}_2,\cdot)d\theta$, where
$$
I_\theta(z,\cdot)\defeq\int_{\bbR^2}K_\theta^v(z,w)\partial_x^2K^v(\cdot,w)dw.
$$
Since $\int_{\bbR^2}\partial_x^2K^v(\cdot,w)dw=0$, we can use an argument similar to the proof of Proposition {\it\ref{prop:Lpkernelestimate}} and have that
$$
|I_\theta(z,\cdot)|\lesssim\theta^{-(1+\kappa)/2}{\sf G}_\theta(z-\cdot)
$$
for any small $\kappa>0$. Consequently we see that
\begin{align*}
\lim_{\varepsilon\to0} \big| \bbE\big[ \mcQ_\theta({\sf\Pi}_z^{\varepsilon}\tau)(z) \big] \big|
&\lesssim\sum_{k\in\bbZ}\int_0^1d\theta_1\int_{\bbR^2} |\mcQ_\theta^v(z,w)| \, \big| I_{\theta_1}(z+k{\bf e}_2,w)\big| \, dw   \\
&\lesssim\sum_{k\in\bbZ}\int_0^1\theta_1^{-1/2}{\sf G}_{\theta+\theta_1}(k{\bf e}_2) \, d\theta_1
\lesssim\int_0^1\theta_1^{-(1+\kappa)/2}(\theta+\theta_1)^{-3/4}d\theta_1.
\end{align*}
Since $(\theta+\theta_1)^{-3/4}\lesssim\theta^{-3/4}\wedge\theta_1^{-3/4}$, by considering the integrals over $\theta_1\in(0,\theta]$ and $\theta_1\in(\theta,1]$ separately we have that
$$
\lim_{\varepsilon\to0} \big|\bbE\big[ \mcQ_\theta({\sf\Pi}_z^{\varepsilon}\tau)(z) \big]\big| \lesssim \theta^{-(1+2\kappa)/4}.
$$
Since $|\tau|=2\alpha_0-2<-1$ the estimate \eqref{eq:proof:asmp1exp} is obtained.

\ssk

Finally we consider $\sigma$. The action of the BPHZ model ${\sf\Pi}_z^\varepsilon$ is
\begin{align*}
{\sf\Pi}_z^{\varepsilon}\sigma
&=\int_{(\bbR\times\bbT)^4}\prod_{i=1}^3\big(\widetilde{K^v}(\cdot,w_i)-\widetilde{K^v}(z,w_i)\big)
\widetilde{\partial_x^2K^v}(\cdot,w_4)\prod_{i=1}^4\xi^\varepsilon(w_i)dw_i\\
&\quad+3\ell_v^\varepsilon(\cdot,\tau)
\bigg(\int_{\bbR\times\bbT}\big(\widetilde{K^v}(\cdot,w)-\widetilde{K^v}(z,w)\big)
\xi^\varepsilon(w)dw\bigg)^2
\end{align*}
Hence its expectation (in the limit) is
\begin{align*}
\lim_{\varepsilon\to0}\bbE[{\sf\Pi}_z^\varepsilon\sigma]
=-3
\bigg(\int_{\bbR\times\bbT}\big(\widetilde{K^v}(\cdot,w)-\widetilde{K^v}(z,w)\big)^2
dw\bigg)
\bigg(\int_{\bbR\times\bbT}\widetilde{K^v}(z,w)\widetilde{\partial_x^2K^v}(\cdot,w)dw\bigg).
\end{align*}
The last integral has been estimated in terms of $I_\theta(z,\cdot)$.
For the integral 
$$
J(\cdot,z)\defeq\int_{\bbR\times\bbT}\big(\widetilde{K^v}(\cdot,w)-\widetilde{K^v}(z,w)\big)^2
dw,
$$
by using the H\"older estimate \eqref{prop:K"kernelestimate2}, for any $\beta\in(0,1/2)$, we have
\begin{align*}
J(\cdot,z)&\lesssim\|\cdot-z\|_\mfs^{2\beta}
\int_{\bbR\times\bbT}\bigg(\int_0^1\theta^{-(2+\beta)/4}\big(\widetilde{{\sf G}_\theta}(\cdot-w)+\widetilde{{\sf G}_\theta}(z-w)\big)d\theta\bigg)^2dw\\
&\lesssim\|\cdot-z\|_\mfs^{2\beta}
\int_0^1d\theta_1\int_0^1\theta_1^{-(2+\beta)/4}\theta_2^{-(2+\beta)/4}
\sum_{u,v\in\{\cdot,z\}}\widetilde{{\sf G}_{\theta_1+\theta_2}}(u-v)d\theta_2\\
&\lesssim\|\cdot-z\|_\mfs^{2\beta}
\int_0^1d\theta_1\int_0^1\theta_1^{-(2+\beta)/4}\theta_2^{-(2+\beta)/4}(\theta_1+\theta_2)^{-3/4}d\theta_2\\
&\lesssim\|\cdot-z\|_\mfs^{2\beta}
\int_0^1d\theta_1\int_0^{\theta_1}\theta_1^{-(5+\beta)/4}\theta_2^{-(2+\beta)/4}d\theta_2\\
&\lesssim\|\cdot-z\|_\mfs^{2\beta}
\int_0^1\theta_1^{-(3+2\beta)/4}d\theta_1\lesssim\|\cdot-z\|_\mfs^{2\beta}.
\end{align*}
We have as a consequence
\begin{align*}
\lim_{\varepsilon\to0} \big|\bbE\big[ \mcQ_\theta({\sf\Pi}_z^{\varepsilon}\sigma)(z) \big]\big| 
&\lesssim \sum_{k\in\bbZ}\int_0^1d\theta_1\int_{\bbR^2}\|w-z\|_\mfs^{2\beta} |\mcQ_\theta^v(z,w)| \, \big|I_{\theta_1}(z+k{\bf e}_2,w)\big| \, dw.
\end{align*}
Since $\|w-z\|_\mfs^{2\beta}|\mcQ_\theta^v(z,w)|\lesssim\theta^{\beta/2}{\sf G}_\theta(z-w)$ we have
$$
\lim_{\varepsilon\to0} \big|\bbE\big[ \mcQ_\theta({\sf\Pi}_z^{\varepsilon}\sigma)(z) \big]\big| \lesssim \theta^{(2\beta-1-2\kappa)/4}
$$
similarly to $\tau$. Since $|\sigma|=4\alpha_0-2<0$ and $2\beta-1<0$, the estimate \eqref{eq:proof:asmp1exp} is obtained.

\bigskip

\subsection{Examples satisfying Assumption \ref{asmp2}{\boldmath $.$} \hspace{0.15cm}}
\label{SubsectionAsmp2}

In this section we give some examples satisfying \textbf{\textit{Assumption \ref{asmp2}}}. Recall from \cite{Hai14} that we can associate to each character $\vcst_v^\varepsilon(\tau)$ a directed graph, called Feynman diagram, whose edges are related with kernels $\partial_z^{\bf n}K^v$ and estimate it by using the singularity of each kernel around the origin. To estimate the difference between $\vcst_v^\varepsilon(z, \tau)$ and $\fcst^\varepsilon_{a(v(z))}(\tau)$, we show that the difference between $K^v$ and $Z^{a(v(\cdot))}$ whose coefficient is frozen at the root is sufficiently regular. First recall from Theorem {\it\ref{thm:SchauderK+R}} that the integration operator $K^v((t,x),(s,y))$ can be replaced with the operator
$$
\big(\partial_t-L^{a(v)}+c\big)^{-1}(\cdot)=\int_{-\infty}^tQ_{ts}^v(\cdot)ds
$$
up to the cost of the good remainder $S^v$ sending $\mcC_{\mfs}^{-2+,v}(\bbR\times\bbT)$ into $\mcC_{\mfs}^{2+}(\bbR\times\bbT)$. Furthermore by Proposition {\it \ref{Qbdecomposition}} below we can replace $Q_{ts}^v(x,y)$ with $Z_{t-s}^{a(v(t,x))}(x-y)$ up to the cost of a less singular kernel.

\medskip

In this section, we consider spacetime two-variable integral kernels $A_{ts}(x,y)$ defined for $x,y\in\bbR$ and $s<t\in\bbR$. 
For such kernels $A_{ts}(x,y)$ and $B_{ts}(x,y)$, we denote by $*$ the spacetime convolution
$$
(A*B)_{ts}(x,y) \defeq \int_{(s,t)\times\bbR^d}A_{tu}(x,z)B_{us}(z,y) \, dudz
$$
as defined in Appendix \textsf{\ref{SubsectionGaussian}}.
By abuse of notation, we use the same symbol to denote the usual convolution $(C*D)_t(x)=\int_{(0,t)\times\bbR}C_{t-s}(x-y)D_s(y)dy$ for spacetime functions $C_t(x)$ and $D_t(x)$ defined for $x\in\bbR$ and $t>0$.
We first prove that the integral kernel
$$
Q_{ts}^{k,p}(x,y) \defeq \big\{ (\partial_x^kQ^v)*(\partial_x^2Q^v)^{*p} \big\}_{ts}(x,y)
$$
has a good estimate and can be replaced with
\begin{align*}
Z_{ts}^{k,p}(x,y) &\defeq \frac1{p!} \, (\partial_\lambda^p\partial_x^kZ_{t-s}^\lambda)(x-y)\vert_{\lambda=a(v(t,x))}   \\
&= \big\{ (\partial_x^kZ^{\lambda})*(\partial_x^2 Z^\lambda)^{*p} \big\}_{t-s}(x-y)\vert_{\lambda=a(v(t,x))},
\end{align*}
up to a better kernel 
$$
D_{ts}^{k,p}(x,y)\defeq Q_{ts}^{k,p}(x,y)-Z_{ts}^{k,p}(x,y).
$$
For simplicity we write 
$$
b=a(v)
$$
and
$$
G_t^{(c_0,c_1)}(x)\defeq \frac{e^{-c_0t}}{t^{1/2}}\exp\bigg(-c_1\frac{|x|^2}t\bigg)
$$
in what follows. We omit the letter $(c_0,c_1)$ if it is not important. Then the functions represented by $G_t(x)$ can be different from one line to another.

\medskip

\begin{prop}\label{Qbdecomposition}
Let $b$ be an element of $\mcC_{\mfs}^\alpha(\bbR\times\bbT)$ with $\alpha\in(0,1]$ and such that $0<b_1\defeq\inf b\le b_2\defeq\sup b<\infty$.
Then for any $\beta\in(0,\alpha)$ and any $k\in\{0,1,2\}$, there exist constants $c_0,c_1>0$ and $C>0$ such that one has
\begin{align}
\label{DecompositionQ1}
\big|Q_{ts}^{k,p}(x,y)\big| + \big| Z_{ts}^{k,p}(x,y) \big|
&\le C^p(t-s)^{-k/2}{G}_{t-s}^{(c_0,c_1)}(x-y),\\
\label{DecompositionQ1'}
\big| D_{ts}^{k,p}(x,y) \big|
&\le C^p(t-s)^{(\beta-k)/2}{G}_{t-s}^{(c_0,c_1)}(x-y)
\end{align}
for any $p\in\bbN$, $x,y\in\bbR$, and $s<t\in\bbR$.
\end{prop}

\medskip

\begin{Dem}
For the estimate \eqref{DecompositionQ1} of $Z_{ts}^{k,p}$, it is sufficient to show that there exist constants $c_0,c_1,C>0$ such that
\begin{equation}\label{proof:Qbdecomposition:Z}
\big| \partial_\lambda^p\partial_x^kZ_t^\lambda(x) \big| \le p! \, C^p \, t^{-k/2} \, G_t^{(c_0,c_1)}(x)
\end{equation}
holds for any $\lambda\in[b_1,b_2]$.
This can be proved by a direct calculation of the heat kernel $e^{t\Delta}(x)=\frac1{\sqrt{4\pi t}}\exp(-\frac{|x|^2}{4t})$.
Since the estimate \eqref{DecompositionQ1} of $Q_{ts}^{k,p}$ immediately follows from \eqref{DecompositionQ1'} and the estimate of $Z_{ts}^{k,p}$, we show \eqref{DecompositionQ1'} in the rest of the proof. 

From the proof of Theorem {\it \ref{thm:GaussGamma}}, we have the decomposition
\begin{equation}\label{ExpansionQ}
\partial_x^kQ_{ts}^v(x,y)
=\partial_x^kL_{ts}(x,y)+R_{ts}^k(x,y),
\end{equation}
where 
$$
L_{ts}(x,y)\defeq Z_{t-s}^{b(s,y)}(x-y)
$$ 
and $|R_{ts}^k(x,y)|\lesssim(t-s)^{(\alpha-k)/2}G_{t-s}(x-y)$ for any $k\in\{0,1,2\}$. By the estimate \eqref{proof:Qbdecomposition:Z} we can replace $\partial_x^kL_{ts}(x,y)$ above with $\partial_x^kZ_{ts}^\lambda(x-y)\vert_{\lambda=b(t,x)}$ up to the cost
$$
\big| b(s,y)-b(t,x) \big| \sup_{\nu\in[b_1,b_2]} \big| \partial_\nu\partial_x^kZ_{t-s}^\nu(x-y) \big|
\lesssim (t-s)^{(\alpha-k)/2} \, G_{t-s}^{(c_0,c_1')}(x-y)
$$
with some $c_1'\in(0,c_0)$. Hence we obtain \eqref{DecompositionQ1'} for $p=0$.

We show the inequality for generic $p$ by induction. By the definition of $D^{k,p}$ and the decomposition \eqref{ExpansionQ} we have
\begin{align*}
Q^{k,p+1}&=Q^{k,p}*\partial_x^2Q^v = Z^{k,p}*\partial_x^2L+Z^{k,p}*R^k+D^{k,p}*\partial_x^2Q^v.
\end{align*}
To show that $Z^{k,p}*\partial_x^2L$ is replaced with $Z^{k,p+1}$, it is sufficient to show that there exist constants $c_0',c_1',C'>0$ such that
\begin{equation}\label{proof:Qbdecomposition:WZ}
\sup_{\mu,\lambda\in[b_1,b_2]} \big| (Z^{\lambda,k,p}*\partial_\mu\partial_x^2Z^\mu)_t(x) \big|
\lesssim (C')^pt^{-k/2}G_{t}^{(c_0',c_1')}(x),
\end{equation}
where $Z^{\lambda,k,p}\defeq(\partial_x^kZ^\lambda)*(\partial_x^2Z^\lambda)^{*p}$.
This can be proved by \eqref{proof:Qbdecomposition:Z} and an argument similar to the proof of Proposition {\it\ref{prop:Lpkernelestimate}}. For $k=2$, since both $Z^{\lambda,2,p}_t(x)$ and $\partial_\mu\partial_x^2Z^\mu_t(x)$ have non-integrable singularity $t^{-1}$, we use the H\"older continuity of $Z^{\lambda,2,p}$ to get
\begin{align*}
\bigg|\int_{\bbR}Z^{\lambda,2,p}_{t-s}(x-y)\partial_\mu\partial_x^2Z^\mu_{s}(y) \, dy\bigg|
\lesssim (C')^p \, (t-s)^{-1-\delta/2} \, s^{-1+\delta/2} \, G_t^{(c_0',c_1')}(x)
\end{align*}
for the integral over $0\le s\le t/2$, and use the H\"older continuity of $\partial_\mu\partial_x^2Z^\mu$ to get
\begin{align*}
\bigg|\int_{\bbR}Z^{\lambda,2,p}_{t-s}(x-y)\partial_\mu\partial_x^2Z^\mu_{s}(y) \, dy\bigg|
\lesssim (C')^p \, (t-s)^{-1+\delta/2} \, s^{-1-\delta/2} \, G_{t}^{(c_0',c_1')}(x)
\end{align*}
for the integral over $t/2\le u\le t$. In the end we can replace $Z^{k,p}*\partial_x^2L$ with $Z^{k,p+1}$ up to the cost
$$
|b(s,y)-b(t,x)|\sup_{\nu\in[b_1,b_2]} \big| (Z^{\lambda,k,p}*\partial_\mu\partial_x^2Z^\mu)_t(x) \big|
\lesssim (C')^p \, (t-s)^{(\alpha-k)/2} \, G_{t-s}^{(c_0',c_1'')}(x-y)
$$
with some $c_1''\in(0,c_0')$. Therefore by writing $D^{k,p+1}$ as
\begin{align*}
D^{k,p+1}=D^{k,p}*\partial_x^2Q^v+Z^{k,p}*R^k+(Z^{k,p}*\partial_x^2L-Z^{k,p+1}),
\end{align*}
we can obtain \eqref{DecompositionQ1'} inductively by an argument similar to the proof of Proposition {\it\ref{prop:Lpkernelestimate}}, as we assume and prove H\"older estimates if necessary.
(The loss $\alpha-\beta$ of regularity is due to the argument based on the H\"older estimates.)
\end{Dem}

\medskip

Equipped with the previous estimates we can now look at three examples where \textbf{\textit{Assumption \ref{asmp2}}} is satisfied. We do not consider the $\bsp$-decoration except in the first half part of Section {\sf\ref{subsubsection:PAM}} because the same argument can be applied to $\bsp$-decorated versions by using Proposition {\it\ref{Qbdecomposition}}.

\subsubsection{Two dimensional parabolic Anderson model.}\label{subsubsection:PAM}

Consider the slightly singular setting of the quasilinear parabolic Anderson model equation 
$$
\partial_tu-a(u)\Delta u=f(u)\xi
$$
on $(t,x)\in(0,\infty)\times\bbT^2$.
We consider the space white noise $\xi$, so one can choose $2/3<\alpha_0<1$. 
We modify the notations on two dimensional spacetime $\bbR\times\bbT$ into the three dimensional spacetime $\bbR\times\bbT^2$, so we write
$$
|{\bf k}|_\mfs = 2k_1+k_2+k_3
$$
for multiindices ${\bf k}=(k_1,k_2,k_3)\in\bbN^3$, and we use the edge decoration $\bsp$ to contract the operator $\Delta\mcI\defeq\mcI_{(0,2,0)}+\mcI_{(0,0,2)}$. Moreover we consider the integral kernels
\begin{align*}
Q_{ts}^{{\bf k},p}(x,y) &\defeq \Big\{ (\partial_{x_2}^{k_2}\partial_{x_3}^{k_3}Q^v)*(\Delta Q^v)^{*p} \Big\}_{ts}(x,y),   \\
Z_{ts}^{{\bf k},p}(x,y) &\defeq \frac1{p!} \, \big(\partial_\lambda^p\partial_{x_2}^{k_2}\partial_{x_3}^{k_3}Z_{t-s}^\lambda\big)(x-y)\vert_{\lambda=a(v(t,x))}.
\end{align*}
instead of $Q_{ts}^{k,p}$ and $Z_{ts}^{k,p}$ defined in the previous section.

The only elements $\tau\in\bbB_-$ with an even number of $\Xi$ type nodes are the trees
\begin{align*}
\tau=\Xi\mcI_{\bf0}^0(\Xi),
\qquad \sigma_i=\mcI_{\bf0}^0(\Xi)\mcI_{2{\bf e}_i}^0(\Xi)\qquad(i\in\{2,3\}),
\end{align*} 
where ${\bf e}_2=(0,1,0)$ and ${\bf e}_3=(0,0,1)$, and their $\bsp$-decorated versions $\tau^{p}=\Xi\mcI_{\bf0}^p(\Xi)$ and $\sigma_i^{p,q}=\mcI_{\bf0}^p(\Xi)\mcI_{2{\bf e}_i}^q(\Xi)$.
The corresponding characters are
\begin{align*}
\vcst_v^\varepsilon((t,x),\tau^p)
&=\int_{(-\infty,t)\times\bbT^2}\widetilde{Q_{ts}^{{\bf0},p}}(x,y)\,C^\varepsilon(x-y)\,dsdy,   \\
\vcst_v^\varepsilon((t,x),\sigma_i^{p,q})
&=\int_{\{(-\infty,t)\times\bbT^2\}^2}\widetilde{Q_{ts}^{{\bf0},p}}(x,y)\,\widetilde{Q_{ts'}^{2{\bf e}_i,q}}(x,y')\,C^\varepsilon(y-y')\,dsds'dydy',
\end{align*}
where $C^\varepsilon(x) = \mathbb{E}\big[\xi^\varepsilon(x) \, \xi^\varepsilon(0)\big]$ and $\widetilde{f}(x,y)=\sum_{{\bf k}\in\bbZ^2}f(x+{\bf k},y)=\sum_{{\bf k}\in\bbZ^2}f(x,y+{\bf k})$ denotes the spatial periodization. By Proposition {\it \ref{Qbdecomposition}}, we can replace $\widetilde{Q_{ts}^{{\bf k},p}}(x,y)$ above by $\widetilde{Z_{ts}^{{\bf k},p}}(x,y)$ up to integrable kernels multiplied by a constant of the form $m^p$. Indeed for $\tau^p$ the difference is estimated by
\begin{align*}
&\int_{(-\infty,t)\times\bbT^2}\widetilde{D_{ts}^{{\bf0},p}}(x,y)\,C^\varepsilon(x-y)\,dsdy
=\int_{(0,\infty)\times\bbT^2}s^{\beta/2}\widetilde{G_s}(y)\,C^\varepsilon(y)\,dsdy \\
&\sim \sum_{{\bf k}\in\bbZ^2}\int_0^{\infty}s^{\beta/2}G_{s}({\bf k})\,ds
\lesssim\int_0^{\infty} e^{-\gamma s}s^{(\beta-1)/2}\,ds < \infty,
\end{align*}
where $a\sim b$ means that $a$ is equal to $b$ up to an $\varepsilon$-uniform remainder term, and $\gamma$ is a positive constant depending on the sufficiently large constant $c>0$ fixed at Section {\sf\ref{SubsectionFunctionSpaces}}.
Similarly for $\sigma_i^{p,q}$ the difference is estimated by
\begin{align*}
&\int_{\{(-\infty,t)\times\bbT^2\}^2} \Big\{ (t-s)^{\beta/2}+(t-s')^{(\beta-2)/2} \Big\}
\, \widetilde{G_{t-s}}(x-y) \, \widetilde{G_{t-s'}}(x-y')\,C^\varepsilon(y-y')\,dsds'dydy'   \\
&\sim \sum_{{\bf k}\in\bbZ^2}\int_{(0,\infty)^2} \Big\{ s^{(\beta-2)/2}+(s')^{(\beta-2)/2}\Big\} \, G_{s+s'}({\bf k})\,dsds'   \\
&\lesssim \int_{(0,\infty)^2} \Big\{ s^{(\beta-2)/2}+(s')^{(\beta-2)/2}\Big\} \, (s+s')^{-1/2}e^{-\gamma(s+s')}\,dsds'
< \infty.
\end{align*}
Thus one can replace $\ell_v^\varepsilon((t,x),\cdot)$ with $l_{a(v(t,x))}^\varepsilon(\cdot)$ up to an $\varepsilon$-uniform remainder term and \textbf{\textit{Assumption \ref{asmp2}}} holds in that case.

To apply Theorem {\it\ref{ThmMainRenormalizedSimplified}}, we see the asymptotics of $\fcst_\lambda^\varepsilon$ for vanishing $\bsp$ decorations.
For $\tau$, we have
\begin{align*}
\fcst^\varepsilon_\lambda(\tau)
&= \int_{(-\infty,t)\times\bbT^2}\widetilde{Z_{t-s}^\lambda}(x-y)\,C^\varepsilon(x-y)\,dsdy
\sim-\frac1{2\pi\lambda}\int_{|y|<1/4}\log|y|\, C^\varepsilon(y)\,dy
\end{align*}
uniformly over $\lambda$ in a compact interval contained in $(0,\infty)$.
Here we have used that $Z_t^\lambda$ is the fundamental solution of the parabolic operator $\partial_t-\lambda\Delta-c$ with $c>0$ and that its integral over $t\in(0,\infty)$ is the Green function of $(c-\lambda\Delta)^{-1}$, which is approximated by $-\frac1{2\pi\lambda}\log|y|$ around the origin.
Similarly, one has
\begin{align*}
\sum_{i\in\{2,3\}}\fcst^\varepsilon_\lambda(\sigma_i)
&= \int_{\{(-\infty,t)\times\bbT^2\}^2}\widetilde{Z_{t-s}^\lambda}(x-y)\,\widetilde{\Delta Z_{t-s'}^\lambda}(x-y')\,C^\varepsilon(y-y')\,dsds'dydy'   \\
&\sim - \frac1{\lambda}\int_{(-\infty,t)\times\bbT^2\times\bbT^2}\widetilde{Z_{t-s}^\lambda}(x-y)\,\delta_0(x-y')\,C^\varepsilon(y-y')\,dsdydy'   \\
&\sim\frac1{2\pi\lambda^2}\int_{|y|<1/4}\log|y|\, C^\varepsilon(y)\,dy.
\end{align*}
Formula \eqref{EqRenormalizedEquation} then takes the form
\begin{align*}
\Bigg(\fcst^\varepsilon_{a(\cdot)}(\tau)f'f+
\bigg(\sum_{i\in\{2,3\}}\fcst^\varepsilon_{a(\cdot)}(\sigma_i)\bigg)
a'f^2\Bigg)(u^\varepsilon) \sim c^\varepsilon\bigg(\frac{f'f}{a}-\frac{a'f^2}{a^2}\bigg)(u^\varepsilon)
\end{align*}
with a constant $c^\varepsilon=-\frac1{2\pi}\int_{|y|<1/4}\log|y|\, C^\varepsilon(y)dy$. This matches the previous works by Bailleul, Debussche \& Hofmanov\'a \cite{BDH}, Furlan \& Gubinelli \cite{FurlanGubinelli} and Otto \& Weber \cite{OttoWeber}.
	
\medskip

\subsubsection{Quasilinear generalized (KPZ) equation with regularized noise.}\label{subsubsection:KPZ}

We return to the two dimensional spacetime $\bbR\times\bbT$. Let $\xi$ be the mildly singular case of a spacetime Gaussian noise of parabolic regularity $\alpha_0-2$ with $1/2<\alpha_0<2/3$ and consider the quasilinear equation
$$
\partial_tu-a(u)\partial_x^2 u=f(u)\xi+g(u)(\partial_xu)^2.
$$
Then the only elements $\tau\in\bbF_-$ with an even number of $\Xi$ type nodes are the trees
\begin{align} \label{EqTreesQgKPZMild}
\tau_1=\Xi\,\mcI_{\bf0}^0(\Xi)=\, \rtdb{1}{\pol{2}{1}{90}\drbl{1}{2}\drc{1,2}}\ ,
\qquad \tau_2=\,\mcI_{(0,1)}^0(\Xi)^2=\,\rtdb{1}{\pol{2}{1}{120}\pol{3}{1}{60}\drl{1}{2,3}\drb{1}\drc{2,3}}\ ,
\qquad \tau_3=\,\mcI_{\bf0}^0(\Xi)\mcI_{(0,2)}^0(\Xi)=\,\rtdb{1}{\pol{2}{1}{120}\pol{3}{1}{60}\drbl{1}{2}\drll{1}{3}\drb{1}\drc{2,3}}\ ,
\end{align}
Here we use the graphical notation similarly to Section {\sf\ref{SubsectionAsmp1}}, and in addition, the thin line denotes the operator $\mcI_{(0,1)}^0$. Since all these trees have homogeneity $2\alpha_0-2>-1$, we can replace the kernel $Q^{k,p}$ by $Z^{k,p}$ up to integrable kernel similarly to Section {\sf\ref{subsubsection:PAM}}.
Thus they satisfy \textbf{\textit{Assumption \ref{asmp2}}} and the counterterm takes the form
\begin{align*}
\Big(\fcst^\varepsilon_{a(\cdot)}(\tau_1)\,f'f + \fcst^\varepsilon_{a(\cdot)}(\tau_2)\,gf^2 + \fcst^\varepsilon_{a(\cdot)}(\tau_3)\,a'f^2\Big)(u^\varepsilon).
\end{align*}
As mentioned in Gerencs\'er \& Hairer's work \cite{GerencserHairer} the renormalization constants are cancelled as follows. 
We write $Z^\lambda(z)=Z_t^\lambda(x)$ for $z=(t,x)$, as a spacetime one variable function.
We assume that the mollifier $\rho^\varepsilon$ is an even function, so is the covariance function $C^\varepsilon$.
By explicit calculations,
\begin{align*}
\fcst^\varepsilon_\lambda(\tau_1) &= -\int_{\bbR\times\bbT}C^\varepsilon(z)\widetilde{Z^\lambda}(z)\,dz,   \\
\fcst^\varepsilon_\lambda(\tau_2) &= -\int_{(\bbR\times\bbT)^2}C^\varepsilon(z-z')\,\widetilde{\partial_xZ^\lambda}(z)\widetilde{\partial_xZ^\lambda}(z')\,dzdz' = -\int_{\bbR\times\bbT}C^\varepsilon(z) \, \widetilde{W_2^\lambda}(z)\,dz,   \\
\fcst^\varepsilon_\lambda(\tau_3) &= -\int_{(\bbR\times\bbT)^2}C^\varepsilon(z-z')\widetilde{Z^\lambda}(z)\,\widetilde{\partial_x^2 Z^\lambda}(z')\,dzdz' = -\int_{\bbR\times\bbT}C^\varepsilon(z) \, \widetilde{W_3^\lambda}(z)\,dz,
\end{align*}
where 
$$
W_2^\lambda(z)\defeq\int_{\bbR^2}\partial_xZ^\lambda(z')\partial_xZ^\lambda(z'-z)dz',\qquad
W_3^\lambda(z)\defeq\int_{\bbR^2}Z^\lambda(z')\partial_x^2Z^\lambda(z'-z)dz'.
$$
As we see from the identity
$$
W_2^\lambda(z) = -W_3^\lambda(z)=\frac1{2\lambda}Z^\lambda(|t|,x)+O(1),
$$ 
that we have 
$$
\lambda\fcst^\varepsilon_\lambda(\tau_2)=-\lambda\fcst^\varepsilon_\lambda(\tau_3)=\fcst^\varepsilon_\lambda(\tau_1)+O(1),
$$ 
our formula matches with Gerencs\'er \& Hairer's formula \cite{GerencserHairer}
$$
\fcst^\varepsilon_{a(u^\varepsilon)}(\tau_1)\bigg(f'f + \frac{gf^2}{a} - \frac{a'f^2}{a}\bigg)(u^\varepsilon).
$$

\subsubsection{Quasilinear generalized (KPZ) equation with space-time white noise.}\label{subsubsection:qgKPZ}

Let $\xi$ be a spacetime Gaussian noise of parabolic regularity $\alpha_0-2$ with $2/5<\alpha_0<1/2$ and consider the stochastic heat equation
$$
\partial_tu-a(u)\partial_x^2 u=\xi.
$$
Then the only elements $\tau\in\bbF_-$ with an even number of $\Xi$ type nodes are the tree $\tau_3$ from \eqref{EqTreesQgKPZMild} together with the trees
\begin{align}\label{treesqheat}
\rtdb{1}{\pol{2}{1}{120}\pol{3}{1}{60}\drbl{1}{2}\drll{1}{3}\drcc{1}\drc{2,3}}\ ,
\qquad
\rtdb{1}{\pol{2}{1}{120}\pol{3}{1}{60}\pol{4}{2}{90}\drbl{1}{2}\drll{1}{3}\drll{2}{4}\drb{1}\drc{3,4}\drcc{2}}\ ,
\qquad
\rtdb{1}{\pol{2}{1}{120}\pol{3}{1}{60}\pol{4}{2}{90}\drbl{1}{3}\drll{1}{2}\drll{2}{4}\drb{1}\drc{3,4}\drcc{2}}\ ,
\qquad
\rtdb{1}{\pol{2}{1}{144}\pol{3}{1}{108}\pol{4}{1}{72}\pol{5}{1}{36}\drbl{1}{2,3,4}\drll{1}{5}\drb{1}\drc{2,3,4,5}}\ ,
\qquad
\rtdb{1}{\pol{2}{1}{135}\pol{3}{1}{45}\pol{4}{2}{135}\pol{5}{2}{90}\pol{6}{2}{45}
\drbl{1}{2}\drbl{2}{4,5}\drll{1}{3}\drll{2}{6}
\drb{1,2}\drc{3,4,5,6}},
\qquad
\rtdb{1}{\pol{2}{1}{120}\pol{3}{1}{60}\pol{4}{2}{120}\pol{5}{2}{60}\pol{6}{4}{120}\pol{7}{4}{60}
\drbl{1}{6}\drll{1}{3}\drll{2}{5}\drll{4}{7}\drb{1,2,4}\drc{3,5,6,7}}\ ,
\end{align}
where the circle dot $\odot$ denotes the node with decorations $\frak{t}={\bf1}$ and $\frak{n}={\bf e}_2$.
Recall again that the renormalization characters vanish on trees with an odd number of $\Xi$ type nodes.
Since $\alpha_0$ is close to $1/2$, trees in \eqref{treesqheat} have homogeneities close to $0$. Indeed, the first three trees have homogeneity $(2\alpha_0-1)$ and the others have homogeneity $(4\alpha_0-2)$.
Hence we can replace the kernel $Q^{k,p}$ by $Z^{k,p}$ in the same way as before, except for the last two trees having edges that do not include the root as an end point.
For such edges $e$ whose lower node is associated with the spacetime variable $w$, the kernel $Q^{k,p}(w,\cdot)$ is replaced with $Z^{\lambda,k,p}(w-\cdot)\vert_{\lambda=b(w)}$, rather than $Z^{\lambda,k,p}(w-\cdot)\vert_{\lambda=b(z)}$, where $z$ is the spacetime variable associated with the root.  
However, by the analyticity of $Z^\lambda$, we can swap them up to the cost of $|b(w)-b(z)|\partial_\lambda Z^{\lambda,k,p}(w-\cdot)\vert_{\lambda=b(z)}$ even for this case.
Since $|b(w)-b(z)|=O(\|w-z\|_\mfs^\alpha)$ smears the singularity of the Feynman diagram by $\alpha$, we can show that \textbf{\textit{Assumption \ref{asmp2}}} is satisfied for the trees in \eqref{treesqheat}.

\ssk

For the remaining tree $\tau_3$ with homogeneity $2\alpha_0-2<-1$, Proposition {\it\ref{Qbdecomposition}} is not sufficient to replace $Q^{k,p}$ by $Z^{k,p}$ with the kernel estimated by $(t-s)^{(\beta-k)/2}G_t(x-y)$.
To overcome this difficulty, we choose a $t$-independent sufficiently regular function $v\in C^2(\bbT)$ -- one choice is $v=e^{\delta\partial_x^2}(u_0)$ for a sufficiently small $\delta>0$ as given after Definition {\it\ref{def:Valpha}}.
In this case, the fundamental solution $Q^v_t$ is $t$-homogeneous, so we denote the convolution $(A*B)_t(x,y)=\int_{(0,t)\times\bbR}A_{t-s}(x,z)B_s(z,y)dy$ for functions $A_t(x,y)$ and $B_t(x,y)$ defined for $x,y\in\bbR$ and $t>0$.

\medskip

\begin{prop}\label{Qbmoredecomposition}
There exists a family $\{Y^{\lambda,k,p}=Y_t^{\lambda,k,p}(x)\}_{\lambda>0,\,k\in\{0,1,2\},\,p\in\bbN}$ of spacetime functions satisfying the following properties.
\begin{enumerate}
\renewcommand{\theenumi}{(\roman{enumi})}
\renewcommand{\labelenumi}{(\roman{enumi})}
\setlength{\itemsep}{5pt}
\item
When $k$ is even, respectively odd, the function $Y_t^{\lambda,k,p}(\cdot)$ is odd, respectively even.
\item
For any $\delta\in(0,1)$ and $k\in\{0,1,2\}$, there exist constants $c_0,c_1>0$ and $C>0$ such that one has
$$
\big| Y_t^{k,p,\lambda}(x,y)\big| \leq C^p \, t^{(\delta-k)/2} \, G_t^{(c_0,c_1)}(x-y)
$$
for any $p\in\bbN$, $x,y\in\bbR$, and $s<t\in\bbR$.
\item
Let $b$ be an element of $C^2(\bbT)$ such that $0<\inf b$.
Then for any $\beta\in(1,2)$ and $k\in\{0,1,2\}$, there exist constants $c_0,c_1>0$ and $C>0$ such that one has
\begin{align}
\label{DecompositionQ2}
\begin{aligned}
\Big| Q_t^{k,p}(x,y) - Z_t^{k,p}(x,y) - b'(x)Y_t^{b(x),k,p}(x-y) \Big| \leq C^p \, t^{(\beta-k)/2} \, G_t^{(c_0,c_1)}(x-y)
\end{aligned}
\end{align}
for any $p\in\bbN$, $x,y\in\bbR$, and $s<t\in\bbR$.
\end{enumerate}
\end{prop}

\medskip

\begin{Dem}
We use an argument similar to the proof of Proposition {\it\ref{Qbdecomposition}}.
From the proof of Theorem {\it \ref{thm:GaussGamma}}, we have the higher decomposition
\begin{equation}\label{ExpansionQmore}
\partial_x^kQ^v
=\partial_x^kL+\partial_x^kL*K+S^k
\end{equation}
than \eqref{ExpansionQ}, where
$$
K_t(x,y)=\big(b(x)-b(y)\big)\partial_x^2L_t(x,y)
$$
and $|S_t^k(x,y)|\lesssim t^{(2-k)/2}G_t(x-y)$ for any $k\in\{0,1,2\}$.
We cannot replace $\partial_x^kL$ in the right hand side of \eqref{ExpansionQmore} with $Z^{k,0}$ as in the proof of Proposition {\it\ref{Qbdecomposition}}, since the remainder does not have a sufficient estimate. Instead, we use a higher decomposition
\begin{align}\label{prooftechnical1}
\begin{aligned}
\partial_x^kL_t(x,y) &= \partial_x^kZ_t^{b(y)}(x-y)   \\
&=\big\{\partial_x^kZ_t^\lambda(x,y)+\big(b(y)-b(x)\big)\partial_\lambda\partial_x^kZ_t^\lambda(x-y)\big\}\big|_{\lambda=b(x)}
+ O\Big(t^{(2-k)/2}G_t(x-y)\Big)   \\
&=Z_t^{k,0}(x,y)+b'(x)(y-x)Z_t^{k,1}(x,y)
+ O\Big(t^{(2-k)/2}G_t(x-y)\Big).
\end{aligned}
\end{align}
On the other hand since $K$ can be decomposed as
\begin{align*}
K_t(x,y) =b'(y)(x-y)\partial_x^2L_t(x,y) + O\big(G_t(x-y)\big),
\end{align*}
by the analyticity of $\lambda\mapsto Z^\lambda$ again, we have
\begin{align}\label{prooftechnical2}
\begin{aligned}
&(\partial_x^kL*K)_t(x,y)\\
&= b'(y) \int_{(0,t)\times\bbR}(\partial_x^kZ_{t-s}^{b(x)})(x-z)(z-y)\partial_x^2Z_s^{b(y)}(z-y)dsdz + O\Big(t^{(2-k)/2}G_t(x-y)\Big)   \\
&= b'(x) \int_{(0,t)\times\bbR}(\partial_x^kZ_{t-s}^{b(x)})(x-z)(z-y)\partial_x^2Z_s^{b(x)}(z-y)dsdz + O\Big(t^{(2-k)/2}G_t(x-y)\Big).
\end{aligned}
\end{align}
By \eqref{prooftechnical1} and \eqref{prooftechnical2} we have the decomposition \eqref{DecompositionQ2} with
\begin{align*}
Y_t^{\lambda,k,0}(x)=-xZ_t^{\lambda,k,1}(x)
+\int_{(0,t)\times\bbR}(\partial_x^kZ_{t-s}^\lambda)(x-y)y\partial_x^2Z_s^\lambda(y)dsdy
\end{align*}
for $p=0$. To show the estimates for generic $p$ by induction, we use another decomposition
\begin{align*}
\partial_x^2Q_t^v(x,y)=\partial_x^2L_t(x,y)+b'(y)\bar{Y}_t^{b(y)}(x-y)+W_t(x,y)
\end{align*}
following from \eqref{ExpansionQmore} and \eqref{prooftechnical2}, where
$$
\bar{Y}_t^\lambda(x)
=\int_{(0,t)\times\bbR}(\partial_x^2Z_{t-s}^\lambda)(x-y)y\partial_x^2Z_s^\lambda(y)dsdz
$$
and $|W_{ts}^k(x,y)|\lesssim G_{t-s}(x-y)$.
By setting 
$$
E_t^{k,p}(x,y)\defeq Q_t^{k,p}(x,y)-Z_t^{k,p}(x,y) - b'(x)Y_t^{b(x),k,p}(x-y),
$$
we have
\begin{align*}
Q_t^{k,p+1}(x,y) &= \big\{ Q^{k,p}*\partial_x^2Q^v \big\}_t(x,y)   \\
&= \big\{ Z^{k,p}*\partial_x^2L\big\}_t(x,y) + b'(y) \big\{ Z^{k,p}*\bar{Y}_t^{b(y)} \big\}_t(x,y) + b'(x) \big\{ Y_t^{b(x),k,p}*\partial_x^2L \big\}_t(x,y)   \\
&\quad+ \big\{ E^{k,p}*\partial_x^2Q^v \big\}_t(x,y) + \cdots,
\end{align*}
where and in the rest of the proof, we omit by ``$\cdots$'' the terms to be proved to have an estimate $O\big(C^pt^{(\beta-k)/2}G_t^{(c_0,c_1)}(x-y)\big)$ for some constants $C,c_0,c_1>0$.
We organize the first three terms in the right hand side.
For the first term, by a calculation similar to \eqref{prooftechnical1}, we have
\begin{align*}
\big\{ Z^{k,p}*\partial_x^2L \big\}_t(x,y)
&= Z_t^{k,p+1}(x,y) + \big(b(y)-b(x)\big) \big\{ Z^{k,p}*\partial_\lambda\partial_x^2Z^\lambda \big\}_t(x,y)|_{\lambda=b(x)}
+\cdots\\
&= Z_t^{k,p+1}(x,y)+b'(x)(y-x)Z^{k,p+2}_t(x,y)
+\cdots.
\end{align*}
For the second term, by the analyticity of $\bar{Y}^\lambda$,
\begin{align*}
b'(y) \big\{ Z^{k,p}*\bar{Y}_t^{b(y)} \big\}_t(x,y)=b'(x)\{Z^{\lambda,k,p}*\bar{Y}^\lambda\}_t(x-y)|_{\lambda=b(x)}
+\cdots.
\end{align*}
For the third term, by the analyticity of $Z^\lambda$, we have
\begin{align*}
b'(x) \big\{ Y_t^{b(x),k,p}*\partial_x^2L \big\}_t(x,y)
= b'(x )\big\{ Y^{\lambda,k,p}*\partial_x^2Z^\lambda \big\}_t(x-y)|_{\lambda=b(x)}
+\cdots.
\end{align*}
Therefore, by defining
\begin{align*}
Y_t^{\lambda,k,p+1}(x)=-xZ_t^{\lambda,k,p+2}(x) + \big\{ Z^{\lambda,k,p}*\bar{Y}^\lambda \big\}_t(x) + \big\{ Y^{\lambda,k,p}*\partial_x^2Z^\lambda \big\}_t(x),
\end{align*}
we can write $E^{k,p+1}$ as
\begin{align*}
E^{k,p+1}=E^{k,p}*\partial_x^2Q^v+\cdots.
\end{align*}
From these identities, we can obtain the desired estimates inductively by an argument similar to the proof of Proposition {\it\ref{prop:Lpkernelestimate}}.
\end{Dem}

\medskip

One sees from \eqref{DecompositionQ2} of Proposition {\it \ref{Qbmoredecomposition}} that 
for the $\bsp$-decorated version $\tau_3^{p,q}=\mcI_{\bf0}^p(\Xi)\mcI_{(0,2)}^q(\Xi)$ we have 
\begin{equation} \label{EqDetailedTrickyCharacter} \begin{split}
\ell_{a(v)}^\varepsilon&\big((t,x),\tau_3^{p,q}\big) 
= \int_{((-\infty,t)\times\bbT)^2} \Big\{\widetilde{Z_{t-s}^{0,p}}(x-y)\widetilde{Z_{t-s'}^{2,q}}(x-y') + (\star)   \\
&+ \sum_{\beta_1+\beta_2>-1}(t-s)^{\beta_1/2}(t-s')^{\beta_2/2}
\widetilde{G_{t-s}}(x-y)\widetilde{G_{t-s'}}(x-y')\Big\} \, C^\varepsilon(s-s',y-y')\,dsds'dydy',
\end{split} \end{equation}
where $(\star)$ is of the form
$$
b'(x) \, \Big\{\widetilde{Y_{t-s}^{b(x),0,p}}(x-y)\widetilde{Z_{t-s'}^{2,q}}(x-y') + \widetilde{Z_{t-s}^{0,p}}(x-y)\widetilde{Y_{t-s'}^{b(x),2,q}}(x-y')\Big\}.
$$
The last term in \eqref{EqDetailedTrickyCharacter} does not matter because one has the $\varepsilon$-uniform estimate
\begin{align*}
\int_{((0,\infty)\times\bbT)^2} s^{\beta_1/2}(s')^{\beta_2/2} \,
&\widetilde{G_s}(y) \, \widetilde{G_{s'}}(y') \, C^\varepsilon(s'-s,y'-y) \, dsds'dydy'   \\
&\lesssim\sum_{k\in\bbZ}\int_0^\infty
s^{(\beta_1+\beta_2)/2} \, G_s(k)\, ds
\lesssim \int_0^{\infty} s^{(\beta_1+\beta_2-1)/2}e^{-\gamma s}ds < \infty.
\end{align*}
The $(\star)$ term in \eqref{EqDetailedTrickyCharacter} is not estimated as above. However if we assume that the mollifier $\rho_\varepsilon$ is spatially even, then since $C^\varepsilon$ is also spatially even, the $(y,y')$-integral $\int(\star)C^\varepsilon(s-s',y-y')dydy'$ vanishes because of the parity of $Y^{\lambda,0,p}$ and $Y^{\lambda,2,q}$. In the end, only the first term of \eqref{EqDetailedTrickyCharacter} survives and \textbf{\textit{Assumption \ref{asmp2}}} is satisfied with
$$
\fcst^\varepsilon_\lambda(\tau_3^{p,q}) = \int_{((-\infty,t)\times\bbT)^2}
\widetilde{Z_{t-s}^{0,p}}(x-y)\widetilde{Z_{t-s'}^{2,q}}(x-y')\,C^\varepsilon(s-s',y-y')\,dsds'dydy'.
$$

\medskip

Finally we show that some of the trees of \eqref{treesqheat} are not involved in Equation \eqref{EqRenormalizedEquation}. The character $l_\lambda^\varepsilon$ vanishes on the first tree since it is of the form $X_2\tau_3$. For the second and third trees, $l_\lambda^\varepsilon$ also vanishes because $Z^\lambda$ is spatially even. Indeed, the value of $l_\lambda^\varepsilon$ on the second tree is given by
\begin{align*}
-\int_{\bbR^2} Z_{-t}^\lambda(x)x\bigg(\int_{((-\infty,t)\times\bbT)^2}\widetilde{Z_{t-s}^\lambda}(x-y)C^\varepsilon(s-s',y-y')\widetilde{\partial_x^2Z_{-s'}^\lambda}(y')dsdyds'dy'\bigg)dtdx,
\end{align*}
where the integral inside the large parentheses is an even function of $x$ and $Z_{-t}^\lambda(x)x$ is an odd function of $x$. Moreover, $l_\lambda^\varepsilon$ vanishes on the fourth tree of \eqref{treesqheat} because $\xi^\varepsilon$ is a centered Gaussian. Indeed by decomposing
\begin{align*}
\fcst_\lambda^\varepsilon\big(\cdot \,,\,\rtdb{1}{\pol{2}{1}{144}\pol{3}{1}{108}\pol{4}{1}{72}\pol{5}{1}{36}\drbl{1}{2,3,4}\drll{1}{5}\drb{1}\drc{2,3,4,5}}\big) = -h_\lambda^\varepsilon\big(\cdot \,,\,\rtdb{1}{\pol{2}{1}{144}\pol{3}{1}{108}\pol{4}{1}{72}\pol{5}{1}{36}\drbl{1}{2,3,4}\drll{1}{5}\drb{1}\drc{2,3,4,5}}\big) + 3h_\lambda^\varepsilon\big(\cdot \,, \,\rtdb{1}{\pol{2}{1}{120}\pol{3}{1}{60}\drbl{1}{2}\drll{1}{3}\drb{1}\drc{2,3}}\big) h_\lambda^\varepsilon\big(\cdot \,, \,\rtdb{1}{\pol{2}{1}{120}\pol{3}{1}{60}\drbl{1}{2,3}\drb{1}\drc{2,3}}\big)
\end{align*}
we see that the right hand side is zero because of Wick theorem for Gaussian random variables. As a consequence, the counterterm takes the form
\begin{align*}
\Bigg\{\fcst^\varepsilon_{a(\cdot)}\big(\, \rtdb{1}{\pol{2}{1}{120}\pol{3}{1}{60}\drbl{1}{2}\drll{1}{3}\drb{1}\drc{2,3}}\, \big) \, a' 
+ \fcst^\varepsilon_{a(\cdot)}\big(\, \rtdb{1}{\pol{2}{1}{135}\pol{3}{1}{45}\pol{4}{2}{135}\pol{5}{2}{90}\pol{6}{2}{45} \drbl{1}{2}\drbl{2}{4,5}\drll{1}{3}\drll{2}{6} \drb{1,2}\drc{3,4,5,6}}\, \big) \,a'a'' 
+ \fcst^\varepsilon_{a(\cdot)}\big(\rtdb{1}{\pol{2}{1}{120}\pol{3}{1}{60}\pol{4}{2}{120}\pol{5}{2}{60}\pol{6}{4}{120}\pol{7}{4}{60}
\drbl{1}{6}\drll{1}{3}\drll{2}{5}\drll{4}{7}\drb{1,2,4}\drc{3,5,6,7}}\, \big) \,(a')^3\Bigg\}(u^\varepsilon),
\end{align*}
which matches Gerencs\'er's formula in Theorem {\it 1.1} of \cite{Gerencser}.

\bigskip

\appendix
\section{Appendix}
\label{SectionAppendix}

In this appendix we prove some technical properties of the fundamental solutions of anisotropic parabolic operators, following the arguments in \cite{Fri64, Eid69}. We believe that the results given here are known but we could not find any suitable references. For the sake of generality for them we work on the space $\bbR^d$ and an arbitrary scaling $\mfs=(\mfs_j)_{j=1}^d\in\bbN^d$. Set
\begin{itemize}
	\item[] $|\mfs|\defeq \sum_{j=1}^d\mfs_j$,
	\item[] $|k|_{\mfs}\defeq \sum_{j=1}^d\mfs_jk_j$,  for $k=(k_j)_{j=1}^d\in\bbN^d$,
	\item[] $\|x\|_{\mfs}\defeq \sum_{j=1}^d|x_j|^{1/\mfs_j}$, for $x=(x_j)_{j=1}^d\in\bbR^d$,
	\item[] $\partial_x^k\defeq \prod_{j=1}^d \partial_{x_j}^{k_j}$, for $k=(k_j)_{j=1}^d\in\bbN^d$.
\end{itemize}
For each $i\in\{1,\dots,d\}$, we denote by $e_i=(0,\dots,\overset{i}{1},\dots,0)$ the $i$-th canonical basis vector. Throughout this appendix we consider the anisotropic parabolic operator
\begin{align}\label{eq:parabolicoperator}
\partial_t-P(t,x,\partial_x)\defeq 
\partial_t-\sum_{|k|_{\mfs}\le N}a_k(t,x)\partial_x^k
\end{align}
with coefficients $a_k(t,x)$ defined in a domain $D=(a,b)\times\bbR^d$, where $-\infty\le a<b\le \infty$. In addition, $N$ is an integer satisfying $N>\max_j\mfs_j$.

\medskip

\begin{defn}\label{defn:FundamentalSolution}
We call a function $Q_{ts}(x,y)$ defined on the domain $D_2\defeq\{(t,s,x,y)\,;\, a<s<t<b,\ x,y\in\bbR^d\}$
\textbf{a fundamental solution} of the operator \eqref{eq:parabolicoperator} if for any $f\in C_b(\bbR^d)$ the integral
$$
F(t,x;s)\defeq\int_{\bbR^d}Q_{ts}(x,y)f(y)dy
$$
exists and satisfies the properties
\begin{align}
\label{defn:FundamentalSolution1}
\big(\partial_t-P(t,x,\partial_x)\big)F(t,x;s)&=0,\qquad t>s,\ x\in\bbR^d,\\
\label{defn:FundamentalSolution2}
\lim_{t\downarrow s}F(t,x;s)&=f(x),\qquad x\in\bbR^d
\end{align}
for any fixed $s\in(a,b)$.
\end{defn}

\medskip

In Appendix \textsf{\ref{SubsectionGaussian}}, we present some preliminary results. We prove the existence of the fundamental solution and Gaussian estimates for it (Theorem {\it \ref{thm:GaussGamma}}) in Appendix \textsf{\ref{SubsectionFundamental}}, and prove uniqueness (Theorem {\it\ref{thm:GammaUnique}}) in Appendix \textsf{\ref{AppendixUniqueness}}. In Appendix \textsf{\ref{AppendixProperties}} we discuss temporally homogeneous operators. The estimates of anisotropic Taylor remainders (Corollary {\it \ref{cor:TaylorFundamentalSolution}}) given in Appendix \textsf{\ref{AppendixTaylor}} are used in the proof of Theorem {\it \ref{thm:multilevelSchauderK}}, although their use is not discussed in detail in this paper -- see \cite{Semi, Singular} for details.

\medskip

\subsection{Gaussian kernels{\boldmath $.$} \hspace{0.15cm}}
\label{SubsectionGaussian}

In this section, we prove some technical properties of exponential functions. For $c>0$ and $\beta\in\bbR$, we define the function
$$
{\sf G}_t^{(c,\beta)}(x)\defeq
t^{(\beta-|\mfs|)/N}\exp\bigg\{-c\sum_{j=1}^d\bigg(\frac{|x_j|^{N/\mfs_j}}{t}\bigg)^{\mfs_j/(N-\mfs_j)}\bigg\},
\qquad t>0,\ x\in\bbR^d.
$$

\medskip

\begin{lem}\label{lem:ConvolutionGauss}
Let $\beta,\beta_1,\beta_2\in\bbR$ and $0<c'<c$. There exists a constant $C$ independent of $t,x,\beta,\beta_1,\beta_2$ such that one has the following.
\begin{enumerate}
\renewcommand{\theenumi}{(\roman{enumi})}
\renewcommand{\labelenumi}{(\roman{enumi})}
\setlength{\itemsep}{5pt}

	\item\label{lem:ConvolutionGauss1}
	For any $\alpha>0$, one has 
$$
\|x\|_{\mfs}^\alpha \, {\sf G}_t^{(c,\beta)}(x)\le C{\sf G}_t^{(c',\alpha+\beta)}(x).
$$ 

	\item\label{lem:ConvolutionGauss2}
	For any $\|h\|_{\mfs}\le t^{1/N}$, one has
$$
{\sf G}_t^{(c,\beta)}(x+h)\le C{\sf G}_t^{(c',\beta)}(x).
$$	
	
	\item\label{lem:ConvolutionGauss3}
	For any $0<s<t$, one has
\begin{align*}
\int_{\bbR^d} {\sf G}_{t-s}^{(c,\beta_1)}(x-y) \, {\sf G}_s^{(c,\beta_2)}(y)\,dy
&\le
C(t-s)^{\beta_1/N}s^{\beta_2/N}\,{\sf G}_t^{(c',0)}(x),\\
\int_{\bbR^d} {\sf G}_{t-s}^{(c',\beta_1)}(x-y) \, {\sf G}_s^{(c,\beta_2)}(y)\,dy
&\le
Ct^{|\mfs|/N}(t-s)^{(\beta_1-|\mfs|)/N}s^{\beta_2/N}\,{\sf G}_t^{(c',0)}(x).
\end{align*}

	\item\label{lem:ConvolutionGauss4}
	If $\beta_1,\beta_2>-N$, one has
$$
\int_0^t\int_{\bbR^d} {\sf G}_{t-s}^{(c,\beta_1)}(x-y) \, {\sf G}_s^{(c,\beta_2)}(y)\,dyds
\le
C\,\frac{\Gamma(\frac{\beta_1+N}N)\Gamma(\frac{\beta_2+N}N)}{\Gamma(\frac{\beta_1+\beta_2+N}N)} \, {\sf G}_t^{(c',\beta_1+\beta_2+N)}(x).
$$

	\item\label{lem:ConvolutionGauss5}
	If $\beta_1>-N+|\mfs|$ and $\beta_2>-N$, one has
$$
\int_0^t\int_{\bbR^d}{\sf G}_{t-s}^{(c',\beta_1)}(x-y) \, {\sf G}_s^{(c,\beta_2)}(y)\,dyds
\le
C\,\frac{\Gamma(\frac{\beta_1-|\mfs|+N}N)\Gamma(\frac{\beta_2+N}N)}{\Gamma(\frac{\beta_1+\beta_2-|\mfs|+N}N)} \, {\sf G}_t^{(c',\beta_1+\beta_2+N)}(x).
$$
\end{enumerate}
\end{lem}

\medskip

\begin{Dem}
The proofs of the statements {\it\ref{lem:ConvolutionGauss1}} and {\it\ref{lem:ConvolutionGauss2}} are elementary and left to the readers. 
For {\it\ref{lem:ConvolutionGauss3}}, note that the elementary inequality
\begin{align}\label{eq:ConvolutionGauss2.5}
{\sf G}_{t-s}^{(c,0)}(x-y){\sf G}_s^{(c,0)}(y)\le t^{|\mfs|/N}(t-s)^{-|\mfs|/N}s^{-|\mfs|/N}{\sf G}_t^{(c,0)}(x)
\end{align}
holds. This inequality reduces to
$$
|x_j|\le t^{\mfs_j/N}
\Bigg\{
\bigg(\frac{|x_j-y_j|}{(t-s)^{\mfs_j/N}}\bigg)^{N/(N-\mfs_j)}+\bigg(\frac{|y_j|}{s^{\mfs_j/N}}\bigg)^{N/(N-\mfs_j)}
\Bigg\}^{(N-\mfs_j)/N},
$$
which follows from the H\"older's inequality. By integration we have
\begin{align*}
&\int_{\bbR^d}{\sf G}_{t-s}^{(c,\beta_1)}(x-y){\sf G}_s^{(c,\beta_2)}(y)dy\\
&=(t-s)^{(\beta_1+|\mfs|)/N}s^{(\beta_2+|\mfs|)/N}\int_{\bbR^d}{\sf G}_{t-s}^{(c',0)}(x-y){\sf G}_{t-s}^{(c-c',0)}(x-y)
{\sf G}_s^{(c',0)}(y){\sf G}_s^{(c-c',0)}(y)dy\\
&\le t^{|\mfs|/N}(t-s)^{\beta_1/N}s^{\beta_2/N}{\sf G}_t^{(c',0)}(x)\int_{\bbR^d}{\sf G}_{t-s}^{(c-c',0)}(x-y){\sf G}_s^{(c-c',0)}(y)dy.
\end{align*}
Since ${\sf G}_t^{(c,0)}(x)\le t^{-|\mfs|/N}$ and
$$
C_c\defeq\int_{\bbR^d}{\sf G}_t^{(c,0)}(x)dx=\int_{\bbR^d}{\sf G}_1^{(c,0)}(x)dx
$$
is a $t$-independent constant, we have
\begin{align*}
&\int_{\bbR^d}{\sf G}_{t-s}^{(c-c',0)}(x-y){\sf G}_s^{(c-c',0)}(y)dy\\
&\le\min\{(t-s)^{-|\mfs|/N},s^{-|\mfs|/N}\}\int_{\bbR^d}{\sf G}_1^{(c-c',0)}(y)dy
\le C_{c-c'}(t/2)^{-|\mfs|/N}.
\end{align*}
Thus we obtain the first inequality of {\it\ref{lem:ConvolutionGauss3}}.
For the second inequality, by similar calculations we have
\begin{align*}
\int_{\bbR^d}&{\sf G}_{t-s}^{(c',\beta_1)}(x-y){\sf G}_s^{(c,\beta_2)}(y)\,dy   \\
&=(t-s)^{\beta_1/N}s^{(\beta_2+|\mfs|)/N}\int_{\bbR^d}{\sf G}_{t-s}^{(c',0)}(x-y){\sf G}_s^{(c',0)}(y){\sf G}_s^{(c-c',0)}(y)dy   \\
&\le t^{|\mfs|/N}(t-s)^{(\beta_1-|\mfs|)/N}s^{\beta_2/N}{\sf G}_t^{(c',0)}(x)\int_{\bbR^d}{\sf G}_s^{(c-c',0)}(y)dy   \\
&\le C_{c-c'}t^{|\mfs|/N}(t-s)^{(\beta_1-|\mfs|)/N}s^{\beta_2/N}{\sf G}_t^{(c',0)}(x).
\end{align*}
We immediately obtain {\it\ref{lem:ConvolutionGauss4}} and {\it\ref{lem:ConvolutionGauss5}} by integrating {\it\ref{lem:ConvolutionGauss3}}.
\end{Dem}

\medskip

\begin{defn} \label{DefnClassGZetaBeta}
For $\beta\in\bbR$, we define ${\sf G}^{\beta}$ as the class of functions $A=A_{ts}(x,y)$ on $D_2$ such that
$$
|A_{ts}(x,y)|\le Ce^{c_0 (t-s)} \, {\sf G}_{t-s}^{(c_1,\beta)}(x-y)
$$
for some positive constants $C,c_0$, and $c_1$.
Moreover, for any $\balpha=(\alpha_i)_{i=1}^d\in\prod_{i=1}^d[0,\mfs_i]$, we define ${\sf G}_{\balpha,0}^\beta$ as the class of functions $A\in{\sf G}^{\beta}$ satisfying
$$
\big|A_{ts}(x+he_i,y)-A_{ts}(x,y)\big|
\le Ce^{c_0 (t-s)} \, |h|^{\alpha_i/\mfs_i}\, {\sf G}_{t-s}^{(c_1,\beta-\alpha_i)}(x-y)
$$
for any $i\in\{1,\dots,d\}$ and $|h|^{1/\mfs_i}\le(t-s)^{1/N}$ -- or, equivalently,
\begin{align}\label{DefnClassGZetaBetaEq}
\big|A_{ts}(x+he_i,y)-A_{ts}(x,y)\big| \le Ce^{c_0 (t-s)} \, |h|^{\alpha_i/\mfs_i} \Big\{ {\sf G}_{t-s}^{(c_1,\beta-\alpha_i)}(x+he_i-y)+{\sf G}_{t-s}^{(c_1,\beta-\alpha_i)}(x-y) \Big\}
\end{align}
for any $h\in\bbR$. We also define ${\sf G}_{0,\balpha}^\beta$ as the set of functions $A\in{\sf G}^\beta$ such that $A_{ts}(y,x)$ is in the class ${\sf G}_{\balpha,0}^\beta$. Finally we define 
$$
{\sf G}_{\balpha,\balpha'}^\beta\defeq{\sf G}_{\balpha,0}^\beta\cap{\sf G}_{0,\balpha'}^\beta.
$$
\end{defn}

\medskip

For any functions $A_{ts}(x,y)$ and $B_{ts}(x,y)$ on $D_2$, we define the spacetime convolution
$$
(A*B)_{ts}(x,y) \defeq \int_{(s,t)\times\bbR^d}A_{tu}(x,z)B_{us}(z,y) \, dudz
$$
for $(t,s,x,y)\in D_2$ if it exists.

\medskip

\begin{lem}\label{lem:ConvolutionGaussClass}
Let $\beta,\beta'\in\bbR$, $\balpha,\balpha'\in\prod_{i=1}^d[0,\mfs_i]$, $A\in{\sf G}_{\balpha,0}^{\beta}$, and $B\in{\sf G}_{0,\balpha'}^{\beta'}$.
We write $\bar{\balpha}<\balpha$ for $\bar{\balpha},\balpha\in\prod_{i=1}^d[0,\mfs_i]$ if $\bar{\alpha}_i<\alpha_i$ for each $i$.
\begin{enumerate}
\renewcommand{\theenumi}{(\roman{enumi})}
\renewcommand{\labelenumi}{(\roman{enumi})}
\setlength{\itemsep}{5pt}

	\item\label{lem:ConvolutionGaussClass1}
Suppose that $\beta,\beta'>-N$, $\max_i\alpha_i<\beta+N$, and $\max_i\alpha_i'<\beta'+N$. Then $A*B\in{\sf G}^{\beta+\beta'+N}_{\balpha,\balpha'}$.

	\item\label{lem:ConvolutionGaussClass2}
Suppose that $\beta\ge-N$, $\beta'>-N$, $\max_i\alpha_i'<\beta'+N$, and $B\in{\sf G}_{(\delta,\dots,\delta),\balpha'}^{\beta'}$ for some $\delta>0$. If
\begin{align}\label{RenormalizeA1}
\bigg|\int_{\bbR^d}A_{ts}(x,y)dy\,\bigg| \lesssim e^{c_0(t-s)}(t-s)^{(\beta+\delta)/N},
\end{align}
then $A*B\in{\sf G}^{\beta+\beta'+N}$.
Moreover, if $\alpha_i'>0$ for each $i$, then $A*B\in{\sf G}_{0,\bar{\balpha}'}^{\beta+\beta'+N}$ for any $\bar{\balpha}'<\balpha'$.
If, in addition to \eqref{RenormalizeA1}, we assume $\max_i\alpha_i<\beta+N+\delta$, $\partial_{x_i}A_{ts}(x,y)\in{\sf G}^{\beta-\mfs_i}$, and
\begin{align}\label{RenormalizeA2}
\bigg|\int_{\bbR^d}\partial_{x_i}A_{ts}(x,y)dy\,\bigg| \lesssim e^{c_0(t-s)}(t-s)^{(\beta-\mfs_i+\delta)/N}
\end{align}
for each $i\in\{1,\dots,d\}$, then $A*B\in{\sf G}_{\balpha,0}^{\beta+\beta'+N}$.

	\item\label{lem:ConvolutionGaussClass3} 
A similar statements to \ref{lem:ConvolutionGaussClass2} holds with the roles of the first and second variables swapped.
\end{enumerate}
\end{lem}

\medskip

\begin{Dem}
Item {\it\ref{lem:ConvolutionGaussClass1}} follows from Lemma {\it \ref{lem:ConvolutionGauss}-\ref{lem:ConvolutionGauss4}}. To show item {\it\ref{lem:ConvolutionGaussClass2}} we decompose
\begin{align*}
(A*B)_{ts}(x,y) &= \int_{(t+s)/2}^tdu\int_{\bbR^d}A_{tu}(x,z)B_{us}(z,y)\,dz \\
&\quad+ \int_s^{(t+s)/2}du\int_{\bbR^d}A_{tu}(x,z)B_{us}(z,y)\,dz
\eqdef I_{ts}(x,y)+J_{ts}(x,y).
\end{align*}
We can prove that $J\in {\sf G}_{\balpha,\balpha'}^{\beta+\beta'+N}$ in the same way as {\it\ref{lem:ConvolutionGaussClass1}}, where we do not have to assume $\max_i\alpha_i<\beta+N$ because $(t-u)^{(\beta-\alpha_i)/N}$ is always integrable over $u\in[s,(t+s)/2]$.
To consider $I_{ts}(x,y)$, we set
$$
I_{ts}(x,y) = \int_s^tC_{tus}(x,y)du,\qquad
C_{tus}(x,y) \defeq \int_{\bbR^d}A_{tu}(x,z)B_{us}(z,y)dz.
$$
If \eqref{RenormalizeA1} holds, then we can estimate $C_{tus}$ by
\begin{align}\label{proof:ConvolutionGaussClass}
\begin{aligned}
&|C_{tus}(x,y)|\\
&\le\bigg|\int_{\bbR^d}A_{tu}(x,z)\,dz\bigg| |B_{us}(x,y)|+ \bigg|\int_{\bbR^d}A_{tu}(x,z)\big(B_{us}(z,y)-B_{us}(x,y)\big)\,dz\bigg|\\
&\lesssim e^{c_0 (t-s)}\bigg\{(t-u)^{(\beta+\delta)/N}{\sf G}_{u-s}^{(c_1,\beta')}(x-y)   \\
&\qquad+\int_{\bbR^d}{\sf G}_{t-u}^{(c_1,\beta)}(x-z)\|z-x\|_{\mfs}^\delta\big({\sf G}_{u-s}^{(c_1,\beta'-\delta)}(z-y) + {\sf G}_{u-s}^{(c_1,\beta'-\delta)}(x-y)\big)\,dz\bigg\}\\
&\lesssim e^{c_0'(t-s)}
(t-u)^{(\beta+\delta)/N}(u-s)^{(\beta'-\delta)/N}{\sf G}_{t-s}^{(c_1',0)}(x-y)
\end{aligned}
\end{align}
for any $c_i'\in(0,c_i)$ for $i\in\{0,1\}$, where we use 
that \eqref{DefnClassGZetaBetaEq} also holds for multidimensional shifts in the second inequality, and that ${\sf G}_{t'}^{(c,0)}(x)\lesssim{\sf G}_t^{(c,0)}(x)$ if $t/2\le t'\le t$ in the third inequality.
Since $\beta+\delta>-N$, the integral in $u$ is finite and we have
\begin{align*}
|I_{ts}(x,y)|
&\lesssim e^{c_0'(t-s)} \, {\sf G}_{t-s}^{(c_1',\beta+\beta'+N)}(x-y).
\end{align*}
Thus we obtain $I\in{\sf G}^{\beta+\beta'+N}$.

Next we show that $I\in{\sf G}_{0,\bar{\balpha}'}^{\beta+\beta'+N}$. Let $i\in\{1,\dots,d\}$ and $|h|^{1/\mfs_i}\le(t-s)^{1/N}$. By \eqref{proof:ConvolutionGaussClass} and Lemma {\it \ref{lem:ConvolutionGauss}-\ref{lem:ConvolutionGauss2}}, we have
$$
\big| C_{tus}(x,y+he_i)-C_{tus}(x,y) \big|
\lesssim
e^{c_0'(t-s)}(t-u)^{(\beta+\delta)/N}(u-s)^{(\beta'-\delta)/N}{\sf G}_{t-s}^{(c_1',0)}(x-y).
$$
On the other hand we also have
\begin{align*}
|C_{tus}(x,y+he_i)-C_{tus}(x,y)|
&\le\int_{\bbR^d}|A_{tu}(x,z)| \, \big| B_{us}(z,y+he_i)-B_{us}(z,y) \big| \, dz   \\
&\lesssim e^{c_0(t-s)}|h|^{\alpha_i'/\mfs_i}\int_{\bbR^d}{\sf G}_{t-u}^{(c_1,\beta)}(x-z) \, {\sf G}_{u-s}^{(c_1,\beta'-\alpha_i')}(z-y) \, dz   \\
&\lesssim e^{c_0(t-s)}|h|^{\alpha_i'/\mfs_i}(t-u)^{\beta/N}(u-s)^{(\beta'-\alpha_i')/N} \, {\sf G}_{t-s}^{(c_1',0)}(x-y).
\end{align*}
By interpolating these two inequalities, for any $\theta\in(0,1)$ we have
\begin{align*}
&\big| I_{ts}(x,y+he_i)-I_{ts}(x,y) \big|
= \int_{(t+s)/2}^t \big| C_{tus}(x,y+he_i)-C_{tus}(x,y)\big| \, du   \\
&\lesssim e^{c_0'(t-s)}|h|^{(1-\theta)\alpha_i'/\mfs_i} \, {\sf G}_{t-s}^{(c_1',0)}(x-y)
\int_{(t+s)/2}^t(t-u)^{(\beta+\theta\delta)/N}(u-s)^{(\beta'-\theta\delta-(1-\theta)\alpha_i')/N}du   \\
&\lesssim e^{c_0'(t-s)}|h|^{(1-\theta)\alpha_i'/\mfs_i} \, {\sf G}_{t-s}^{(c_1',\beta+\beta'+N-(1-\theta)\alpha_i')}(x-y).
\end{align*}

Finally we show that $I\in{\sf G}_{\balpha,0}^{\beta+\beta'+N}$ under the additional assumption of {\it\ref{lem:ConvolutionGaussClass2}}. For each $i\in\{1,\dots,d\}$ and $|h|^{1/\mfs_i}\le(t-s)^{1/N}$, we write $x_h=x+he_i$ and decompose
\begin{align*}
C_{tus}(x_h,y)-C_{tus}(x,y)
&=h\int_0^1\int_{\bbR^d}\partial_{x_i}A_{tu}(x_{\theta h},z)B_{us}(z,y)dz\,d\theta.
\end{align*}
Then by a calculation similar to \eqref{proof:ConvolutionGaussClass}, we have
\begin{align*}
\big| C_{tus}(x_h,y) - C_{tus}(x,y) \big|
&\lesssim e^{c_0'(t-s)}|h|(t-u)^{(\beta-\mfs_i+\delta)/N}(u-s)^{(\beta'-\delta)/N}{\sf G}_{t-s}^{(c_1',0)}(x-y).
\end{align*}
By interpolation between this and \eqref{proof:ConvolutionGaussClass}, we have
\begin{align*}
\big| C_{tus}(x_h,y)-C_{tus}(x,y) \big|
&\lesssim e^{c_0'(t-s)}|h|^{\alpha_i/\mfs_i}(t-u)^{(\beta-\alpha_i+\delta)/N}(u-s)^{(\beta'-\delta)/N}{\sf G}_{t-s}^{(c_1',0)}(x-y)
\end{align*}
for any $\alpha_i\in[0,\mfs_i]$. 
Since $\beta-\alpha_i+\delta>-N$, the integral in $u$ is finite and we have
\begin{align*}
\big| I_{ts}(x_h,y)-I_{ts}(x,y) \big|
&\lesssim e^{c_0'(t-s)}|h|^{\alpha_i/\mfs_i}{\sf G}_{t-s}^{(c_1',\beta_1+\beta_2+N-\alpha_i)}(x-y),
\end{align*}
which completes the proof of {\it\ref{lem:ConvolutionGaussClass2}}.
The proof of {\it\ref{lem:ConvolutionGaussClass3}} is similar.
\end{Dem}

\medskip

\subsection{Existence of a fundamental solution{\boldmath $.$} \hspace{0.15cm}}
\label{SubsectionFundamental}

First we consider the parabolic operator \eqref{eq:parabolicoperator} when the coefficients $a_k$ are constants. Then we write 
\begin{align}\label{eq:parabolicoperatorconst}
\partial_t-P(\partial_x)\defeq\partial_t-\sum_{|k|_{\mfs}\le N}a_k\partial_x^k.
\end{align}

\medskip

\begin{lem}\label{lem:constantcoefficient}
Assume the existence of a constant $\delta>0$ such that the inequality
\begin{align}\label{Petrowskiconst}
\mathop{\text{\rm Re}} P(i\xi)
=\mathop{\text{\rm Re}} \sum_{|k|_{\mfs}\le N}a_k(i\xi)^k\le-\delta\|\xi\|_{\mfs}^N
\end{align}
holds for any $\xi\in\bbR^d$. Then a fundamental solution $Q_{t,s}(x,y)$ of the operator \eqref{eq:parabolicoperatorconst} exists and is of the form $Q_{t,s}(x,y)=Z_{t-s}(x-y)$. (In what follows, we say that $Z_t(x)$ is a fundamental solution of \eqref{eq:parabolicoperatorconst}.) Moreover for any $\varepsilon>0$, $k\in\bbN^d$, and $n\in\bbN$ there exist some positive constants $C$ and $c$ which depend only on $\mfs, N, A\defeq\max_k|a_k|, \delta, \varepsilon, k, n$ such that the inequality
\begin{align}\label{GEforZ}
\big|\partial_t^n\partial_x^kZ_t(x)\big|
\le Ce^{\varepsilon t}{\sf G}_t^{(c,-|k|_{\mfs}-Nn)}(x)
\end{align}
holds for any $t>0$ and $x\in\bbR^d$.
When $(k,n)=(0,0)$, the constant $C$ depends only on $\delta$.
\end{lem}

\medskip

\begin{Dem}
By definition, $Z_t(x)$ is obtained as the Fourier inverse transform of the function $e^{tP(i\xi)}$ of $\xi\in\bbR^d$.
Following the arguments in \cite[Chapter 9, Section 2]{Fri64}, we consider the bound of $e^{tP(i\xi-\eta)}$ for $\eta,\xi\in\bbR^d$.
By the binomial theorem, we can expand 
$$
P(i\xi-\eta)=P(i\xi)+R(\xi,\eta),
$$
where $R(\xi,\eta)$ is a linear combination of monomials $\xi^k\eta^\ell$ with $|k+\ell|_{\mfs}\le N$ and $\ell\neq0$, and with coefficients depending only on $\{a_k\}$.
For any $\varepsilon>0$, by Young's inequality, there exist positive constants $c'$ and $c$ depending only on $A,\varepsilon,\delta$ such that
\begin{align*}
|R(\xi,\eta)|&\le A\sum_{m\ge0, \, n>0,\, m+n\le N}\|\xi\|_{\mfs}^m\|\eta\|_{\mfs}^n\\
&\le \varepsilon+\frac\delta2\|\xi\|_{\mfs}^N+c'\|\eta\|_{\mfs}^N
\le \varepsilon+\frac\delta2\|\xi\|_{\mfs}^N+c\sum_{j=1}^d|\eta_j|^{N/\mfs_j}.
\end{align*}
Therefore, by the condition \eqref{Petrowskiconst}, we have
$$
|e^{tP(i\xi-\eta)}|\le e^{t\mathop{\text{\rm Re}}P(i\xi)}e^{t|R(\xi,\eta)|}\le
\exp\bigg\{t\bigg(\varepsilon-\frac\delta2\|\xi\|_{\mfs}^N+c\sum_{j=1}^d|\eta_j|^{N/\mfs_j}\bigg)\bigg\}.
$$
By using the Cauchy's theorem for each component, for any $\eta\in\bbR^d$ we have
\begin{align*}
|Z_t(x)|&=\bigg|\frac1{(2\pi)^d}\int_{\bbR^d}e^{ix\cdot\xi}e^{tP(i\xi)}d\xi\bigg|
=\bigg|\frac1{(2\pi)^d}\int_{\bbR^d}e^{ix\cdot(\xi+i\eta)}e^{tP(i\xi-\eta)}d\xi\bigg|\\
&\le \frac{e^{\varepsilon t}}{(2\pi)^d}
\exp\bigg(-x\cdot\eta +ct \sum_{j=1}^d|\eta_j|^{N/\mfs_j} \bigg)
\int_{\bbR^d} \exp\bigg(-\frac{\delta t}2\|\xi\|_{\mfs}^N\bigg) d\xi.
\end{align*}
If we choose $\eta_j$ as
\begin{align}\label{choiceeta}
\eta_j=(\mathop{\text{\rm sgn}}x_j)\bigg(\frac{|x_j|}{cp_jt}\bigg)^{1/(p_j-1)},
\end{align}
where $p_j=N/\mfs_j$, then
\begin{align*}
-x_j\eta_j+ct|\eta_j|^{p_j}=-\frac{p_j-1}{p_j}\bigg(\frac{|x_j|^{p_j}}{cp_jt}\bigg)^{1/(p_j-1)},
\end{align*}
which becomes the argument of the exponential function in \eqref{GEforZ}.
The integral with respect to $\xi$ becomes $Ct^{-|\mfs|/N}$ with some constant $C$ depending only on $\delta$.

For the derivatives $\partial_x^kZ_t(x)$, we can derive the desired estimate in a similar way from the identity
\begin{align*}
\partial_x^kZ_t(x)=\frac1{(2\pi)^d}\int_{\bbR^d}e^{ix\cdot\xi}(i\xi)^ke^{tP(i\xi)}d\xi
=\frac1{(2\pi)^d}\int_{\bbR^d}e^{ix\cdot(\xi+i\eta)}(i\xi-\eta)^ke^{tP(i\xi-\eta)}d\xi.
\end{align*}
We decompose $(i\xi-\eta)^k$ into the linear combination of monomials $\xi^\ell\eta^m$ with $\ell+m=k$. The integral of $|\xi^\ell|\exp(-\frac{\delta t}2\|\xi\|_{\mfs}^N)$ over $\xi$ becomes the factor $Ct^{-(|\mfs|+|\ell|_{\mfs})/N}$. For the choice of $\eta$ as in \eqref{choiceeta} we have
\begin{align*}
|\eta^m|=t^{-|m|_{\mfs}/N}\prod_{j=1}^d\bigg(\frac{|x_j|}{cp_jt^{1/p_j}}\bigg)^{m_j/(p_j-1)}.
\end{align*}
Since any powers of $|x_j|/t^{1/p_j}$ are absorbed in the exponential part of \eqref{GEforZ} and the factor $t^{-|m|_{\mfs}/N}$ remains, we have the desired estimate for $\partial_x^kZ_t(x)$. We have similar estimates for the time derivatives because $\partial_t^nZ_t(x)=(P(\partial_x))^nZ_t(x)$.
\end{Dem}

\medskip

Next we consider the operator \eqref{eq:parabolicoperator} with variable coefficients $a_k(t,x)$. 

\medskip

\begin{thm}\label{thm:GaussGamma}
Assume the following conditions for $a_k(t,x)$.
\begin{enumerate}
\renewcommand{\theenumi}{(\alph{enumi})}
\renewcommand{\labelenumi}{(\alph{enumi})}
\setlength{\itemsep}{5pt}
	\item There exists a constant $\delta>0$ such that the inequality
\begin{align}\label{Petrowskiuniform}
\mathop{\text{\rm Re}} P(t,x,i\xi)
=\mathop{\text{\rm Re}} \sum_{|k|_{\mfs}\le N}a_k(t,x)(i\xi)^k
\le-\delta\|\xi\|_{\mfs}^N
\end{align}
holds for any $(t,x)\in D$ and $\xi\in\bbR^d$. 
	\item For some $\alpha\in(0,1]$, one has
\begin{align*}
A&\defeq \max_{|k|_{\mfs}\le N}\sup_{(t,x)\in D}|a_k(t,x)|<\infty,\\
H&\defeq \max_{|k|_{\mfs}\le N}\sup_{(t,x),(s,y)\in D}\frac{|a_k(t,x)-a_k(s,y)|}{\big( |t-s|^{1/N}+\|x-y\|_{\mfs} \big)^\alpha}<\infty.
\end{align*}
\end{enumerate}
Then the operator \eqref{eq:parabolicoperator} has a fundamental solution $Q_{ts}(x,y)$. Moreover for any $k\in\bbN^d$ with $|k|_{\mfs}\le N$, any $\balpha\in\prod_{i=1}^d[0,\mfs_i]$ such that $\alpha_i<N-|k|_{\mfs}+\alpha$, and any $\beta\in(0,\alpha)$, the function $\partial_x^kQ_{ts}(x,y)$ is in the class ${\sf G}_{\balpha,(\beta,\dots,\beta)}^{-|k|_{\mfs}}$ in the sense of Definition \ref{DefnClassGZetaBeta}, where the positive constants $C,c_0$, and $c_1$ involved in this definition depend only on $\mfs, N, \delta,A,H,k,\balpha$, and $\beta$.
\end{thm}

\medskip

We prove this theorem following \cite[Chapter 9]{Fri64}. Let $Z_t^{s,y}(x)$ be the fundamental solution of the operator $\partial_t-P(s,y,\partial_x)$ for fixed $(s,y)$ and set
$$
L_{ts}(x,y)\defeq Z_{t-s}^{s,y}(x-y).
$$
We aim to construct the fundamental solution $Q_{ts}(x,y)$ of the form
\begin{align}\label{GammaZPhi}
Q=L+L*\Phi
\end{align}
with some function $\Phi=\Phi_{ts}(x,y)$. We set
\begin{align*}
K_{ts}(x,y)\defeq \big(P(t,x,\partial_x)-\partial_t\big)L_{ts}(x,y)=\big(P(t,x,\partial_x)-P(s,y,\partial_x)\big)L_{ts}(x,y).
\end{align*}
Then $Q_{ts}(x,y)$ satisfies $\big(\partial_t-P(t,x,\partial_x)\big)Q_{ts}(x,y)=0$ if and only if
\begin{align*}
\Phi=K+K*\Phi.
\end{align*}
This implies that the solution $\Phi$ is formally given by the form
\begin{align}\label{expansionPhi}
\Phi_{ts}(x,y)=\sum_{m=1}^\infty K_{ts}^{(m)}(x,y),\qquad
K^{(m)}\defeq K^{*m}=K^{(m-1)}*K.
\end{align}
It turns out that the series \eqref{expansionPhi} is absolutely convergent and that the function $Q_{ts}(x,y)$ obtained by the formula \eqref{GammaZPhi} is indeed a fundamental solution.

\medskip

\begin{lem}\label{appendixlem:L}
For any $k\in\bbN^d$, $\partial_x^kL_{ts}(x,y)$ is in the class ${\sf G}_{\mfs,(\alpha,\dots,\alpha)}^{-|k|_{\mfs}}$.
\end{lem}

\medskip

\begin{Dem}
The Gaussian estimate and the H\"older estimate for the first variable immediately follow from Lemma {\it \ref{lem:constantcoefficient}}.
The H\"older estimate for the second variable follows from the H\"older continuity of the fundamental solution $Z_t^{\{a_k\}}(x)$ of \eqref{eq:parabolicoperatorconst} with fixed coefficients, the analyticity of $Z_t^{\{a_k\}}(x)$ with respect to $\{a_k\}$, and the H\"older continuity of $a_k(t,x)$. See also Lemma {\it 4} of \cite[Chapter 9]{Fri64}.
\end{Dem}

\medskip

\begin{lem}
The function $K_{ts}(x,y)$ is in the class ${\sf G}_{(\alpha,\dots,\alpha),(\alpha,\dots,\alpha)}^{\alpha-N}$.
\end{lem}

\medskip

\begin{Dem}
Since $K_{ts}(x,y)=\big(P(t,x,\partial_x)-P(s,y,\partial_x)\big)L_{ts}(x,y)$, we have
\begin{align}\label{eq:GaussEstimateK}
\begin{aligned}
|K_{ts}(x,y)|
&\lesssim \big( |t-s|^{1/N}+\|x-y\|_{\mfs} \big)^\alpha e^{\varepsilon (t-s)} \, {\sf G}_{t-s}^{(c,-N)}(x-y)   \\
&\lesssim e^{\varepsilon (t-s)} \, {\sf G}_{t-s}^{(c_1,\alpha-N)}(x-y)
\end{aligned}
\end{align}
for some $c_1<c$ by Lemma {\it \ref{lem:ConvolutionGauss}-\ref{lem:ConvolutionGauss1}}. 
Moreover, for $i\in\{1,\dots,d\}$ and $|h|^{1/\mfs_i}\le(t-s)^{1/N}$,
\begin{align*}
&\big| K_{ts}(x+he_i,y)-K_{ts}(x,y) \big|   \\
&\le\big|\big(P(t,x+he_i,\partial_x)-P(t,x,\partial_x)\big)L_{ts}(x+he_i,y)\big|   \\
&\quad+\big|\big(P(t,x,\partial_x)-P(s,y,\partial_x)\big)\big(L_{ts}(x+he_i,y)-L_{ts}(x,y)\big)\big|   \\
&\lesssim e^{\varepsilon (t-s)}
\Big\{|h|^{\alpha/\mfs_i}{\sf G}_{t-s}^{(c,-N)}(x-y) + (|t-s|^{1/N}+\|x-y\|_{\mfs})^\alpha  |h|^{\alpha/\mfs_i}{\sf G}_{t-s}^{(c,-N-\alpha)}(x-y)\Big\}   \\
&\lesssim e^{\varepsilon (t-s)}|h|^{\alpha/\mfs_i} \, {\sf G}_{t-s}^{(c,-N)}(x-y).
\end{align*}
The H\"older estimate for the second variable is obtained in a similar way.
\end{Dem}

\medskip

\begin{lem}\label{lem:boundPhi}
The function $\Phi_{ts}(x,y)$ is in the class ${\sf G}_{(\beta,\dots,\beta),(\beta,\dots,\beta)}^{\alpha-N}$ for any $\beta<\alpha$.
\end{lem}

\medskip

\begin{Dem}
First we show the estimates
\begin{align}\label{K(m)bound}
\big|K_{ts}^{(m)}(x,y)\big|\le
Ce^{\varepsilon (t-s)} \frac{B^m (t-s)^{m\alpha/N-1}}{\Gamma(\frac{m\alpha-|\mfs|}N)} \, {\sf G}_{t-s}^{(c',0)}(x-y)
\end{align}
for some constants $c',C,B>0$, which depend only on $\mfs,N,\delta,A,H,\varepsilon$.
Let $m_0$ be the smallest integer $m_0$ such that $m_0\alpha>|\mfs|$. Up to $m\le m_0$, \eqref{K(m)bound} is inductively obtained by Lemma {\it \ref{lem:ConvolutionGauss}-\ref{lem:ConvolutionGauss4}}. Indeed, starting from \eqref{eq:GaussEstimateK} we have
\begin{align*}
|K_{ts}^{(m)}(x,y)|&\lesssim e^{\varepsilon (t-s)}\big({\sf G}^{(c_{m-1},(m-1)\alpha-N)}*{\sf G}^{(c_1,\alpha-N)}\big)_{t-s}(x-y)\\
&\lesssim e^{\varepsilon (t-s)} \, {\sf G}_{t-s}^{(c_m,m\alpha-N)}(x-y)
\end{align*}
for some $c_m<c_{m-1}$.
For $m>m_0$, we use Lemma {\it \ref{lem:ConvolutionGauss}-\ref{lem:ConvolutionGauss5}} to obtain
\begin{align*}
|K_{ts}^{(m)}(x,y)|&\le
e^{\varepsilon (t-s)}\frac{B^{m-1}}{\Gamma(\frac{(m-1)\alpha-|\mfs|}N)} \, \big({\sf G}^{(c',(m-1)\alpha-N)}*{\sf G}^{(c_1,\alpha-N)}\big)_{t-s}(x-y)   \\
&\le e^{\varepsilon (t-s)}\frac{B^{m-1}D'\Gamma(\frac\alpha{N})}{\Gamma(\frac{m\alpha-|\mfs|}N)} \, {\sf G}_{t-s}^{(c',m\alpha-N)}(x-y).
\end{align*}
Hence \eqref{K(m)bound} holds with $B=D'\Gamma(\frac{\alpha}N)$.
Summing up \eqref{K(m)bound} over $m\ge1$, we have
\begin{align*}
|\Phi_{ts}(x,y)|\lesssim e^{c_0(t-s)}(t-s)^{\alpha/N-1} \, {\sf G}_{t-s}^{(c',0)}(x-y)=e^{c_0 (t-s)} \, {\sf G}_{t-s}^{(c',\alpha-N)}(x-y)
\end{align*}
for some $c_0$. The H\"older estimates are obtained by applying Lemma {\it \ref{lem:ConvolutionGaussClass}-\ref{lem:ConvolutionGaussClass1}} to the formula
$$
\Phi = K+K*\Phi = K+\Phi*K.
$$
We let the reader complete the associated details.
\end{Dem}

\medskip

\begin{Dem}[of Theorem \ref{thm:GaussGamma}]
We have the Gaussian and H\"older estimates of $\partial_x^kQ$ by applying Lemma {\it \ref{lem:ConvolutionGaussClass}-\ref{lem:ConvolutionGaussClass2}} to the formula
$$
\partial_x^kQ=\partial_x^kL+\partial_x^kL*\Phi.
$$
By Lemmas {\it \ref{appendixlem:L}} and {\it \ref{lem:boundPhi}} we have
\begin{equation}\label{eq:proofofthm:GaussGammaR}
\partial_x^kL*\Phi\in
{\sf G}_{\balpha,(\beta,\dots,\beta)}^{\alpha-|k|_{\mfs}}
\end{equation}
for any $k\in\bbN^d$ with $|k|_{\mfs}\le N$, any $\balpha\in\prod_{i=1}^d[0,\mfs_i]$ such that $\alpha_i<N-|k|_{\mfs}+\alpha$, and any $\beta\in(0,\alpha)$. Note that
\begin{align}\label{eq:proofofthm:GaussGamma}
\begin{aligned}
\bigg|\int_{\bbR^d}\partial_x^\ell L_{ts}(x,y)dy\bigg|
&=\bigg|\int_{\bbR^d}\big(\partial_x^\ell Z_{ts}^{s,y}(x-y)-\partial_x^\ell Z_{ts}^{s,x}(x-y)\big)dy\bigg|\\
&\lesssim e^{\varepsilon (t-s)}\int_{\bbR^d}\|x-y\|_{\mfs}^\alpha {\sf G}_{t-s}^{(c,-|\ell|_{\mfs})}(x-y)\, dz\\
&\lesssim e^{\varepsilon (t-s)} (t-s)^{(\alpha-|\ell|_{\mfs})/N}
\end{aligned}
\end{align}
for any $\ell\in\bbN^d\setminus\{0\}$.

We can check that $Q_{ts}(x,y)$ is indeed a fundamental solution by an argument similar to the proof of Theorem {\it11} of \cite[Section 6, Chapter 1]{Fri64}.
\end{Dem}

\medskip

For $\gamma>0$, define $\mcC_{\mfs}^\gamma(D)$ as the parabolic $\gamma$-H\"older space on the domain $D=(a,b)\times\bbR^d$, that is, $f\in\mcC_{\mfs}^\gamma(D)$ if $\partial_t^n\partial_x^kf$ exists and is bounded for any $Nn+|k|_{\mfs}<\gamma$, and one has
$$
\big| \partial_t^n\partial_x^kf(t+s,x)-\partial_t^n\partial_x^kf(t,x) \big| \lesssim |s|^{(\gamma-Nn-|k|_\mfs)/N}
$$
for each $t,t+s\in(a,b)$ if $Nn+|k|_\mfs<\gamma\le N(n+1)+|k|_\mfs$, and
$$
\big| \partial_t^n\partial_x^kf(t,x+he_i)-\partial_t^n\partial_x^kf(t,x) \big| \lesssim |h|^{(\gamma-Nn-|k|_\mfs)/\mfs_i}
$$
for each $i\in\{1,\dots,d\}$ and $h\in\bbR^d$ if $Nn+|k|_\mfs<\gamma\le Nn+|k+e_i|_\mfs$.

\medskip

\begin{thm}\label{thm:InverseParabolic}
Let $Q_{ts}(x,y)$ be a fundamental solution of the operator \eqref{eq:parabolicoperator} satisfying the assumptions of Theorem \ref{thm:GaussGamma}, and let $\beta\in(0,\alpha)$.
When $a>-\infty$, for any $f\in \mcC_\mfs^\beta(D)$ the function
$$
F(t,x)=\int_a^t\int_{\bbR^d}Q_{ts}(x,y)f(s,y)dyds
$$
belongs to $\mcC_{\mfs}^{\beta+N}(D)$ and satisfies
\begin{equation}\label{thm:InverseParabolic:eq}
\big(\partial_t-P(t,x,\partial_x)\big)F(t,x)=f(t,x).
\end{equation}
Moreover $a=-\infty$ is allowed if both the constant term $a_0$ of \eqref{eq:parabolicoperator} and the constant $c_0$ in Definition \ref{DefnClassGZetaBeta} are strictly negative numbers.
\end{thm}

\medskip

\begin{Dem}
\eqref{thm:InverseParabolic:eq} is obtained by an argument similar to the proof of Theorem {\it9} of \cite[Section 5, Chapter 1]{Fri64}.
We focus on the proof of $F\in\mcC_\mfs^{\beta+N}(D)$ when $a>-\infty$. The proof for the case $a=-\infty$ is similar.
Moreover, by \eqref{thm:InverseParabolic:eq}, we have only to show the estimates of $\partial_x^kF(t,x)$.

\ssk

We decompose
$$
F(t,x)=\int_a^tF(t,x;s)ds,\qquad
F(t,x;s)\defeq\int_{\bbR^d}Q_{ts}(x,y)f(s,y)dy
$$
and consider the estimates of $Q_{ts}f$.
First, note that the identity $\int_{\bbR^d}Q_{ts}(x,y)dy=e^{a_0(t-s)}$ holds by the uniqueness of the fundamental solution.
Therefore, for any $|k|_\mfs\le N$, 
\begin{align*}
|\partial_x^kF(t,x;s)|
&\le\int_{\bbR^d}|\partial_x^kQ_{ts}(x,y)| \, \big| f(s,y)-f(s,x)\big| \, dy
+\bigg|\int_{\bbR^d}\partial_x^kQ_{ts}(x,y)dy\bigg||f(s,x)|\\
&\lesssim e^{c_0(t-s)}\int_{\bbR^d}{\sf G}_{t-s}^{(c_1,-|k|_\mfs)}(x-y)\|y-x\|_\mfs^\beta dy
+e^{a_0(t-s)}{\bf 1}_{k=0}\\
&\lesssim e^{c_0(t-s)}(t-s)^{(\beta-|k|_\mfs)/N}+e^{a_0(t-s)}{\bf 1}_{k=0}.
\end{align*}
Thus $\partial_x^kF(t,x)=\int_a^t\partial_x^kF(t,x;s)ds$ converges for any $|k|_\mfs\le N$.

Next, for any $k\in\bbN^d$ and $i\in\{1,\dots,d\}$ such that $|k|_\mfs<\beta+N\le |k+e_i|_\mfs$, we consider the H\"older estimate of $\partial_x^kF(t,x+he_i;s)-\partial_x^kF(t,x;s)$.
In the region $|h|^{1/\mfs_i}\ge(t-s)^{1/N}$, by estimating each term separately, we have
\begin{align*}
\int_{(t-|h|^{N/\mfs_i})\vee a}^t &\big| \partial_x^kF(t,x+he_i;s)-\partial_x^kF(t,x;s) \big| \, ds   \\
&\lesssim \int_{(t-|h|^{N/\mfs_i})\vee a}^t\{e^{c_0(t-s)}(t-s)^{(\beta-|k|_\mfs)/N}+e^{a_0(t-s)}{\bf 1}_{k=0}\}ds
\lesssim |h|^{(\beta+N-|k|_\mfs)/\mfs_i}.
\end{align*}
In the region $|h|^{1/\mfs_i}<(t-s)^{1/N}$, by using the H\"older estimate of $Q_{ts}$, we have
\begin{align*}
\big| \partial_x^kF(t,x+he_i;s)-\partial_x^kF(t,x;s) \big| &\leq \int_{\bbR^d}|\partial _x^kQ_{ts}(x+he_i,y)-\partial_x^kQ_{ts}(x,y)| \, \big| f(s,y)-f(s,x)\big| \, dy   \\
&\lesssim e^{c_0(t-s)}|h|{\alpha_i/\mfs_i}
\int_{\bbR^d}{\sf G}_{t-s}^{(c_1,-|k|_\mfs-\alpha_i)}(x-y)\|y-x\|_\mfs^\beta dy\\
&\lesssim e^{c_0(t-s)}|h|^{\alpha_i/\mfs_i}(t-s)^{(\beta-|k|_\mfs-\alpha_i)/N}
\end{align*}
for any $\alpha_i<(\alpha+N-|k|_\mfs)\wedge\mfs_i$. By choosing $\alpha_i=\gamma+N-|k|_\mfs$ with $\gamma\in(0,\beta)$, we have
\begin{align*}
\int_a^{t-|h|^{N/\mfs_i}} &\big| \partial_x^kF(t,x+he_i;s)-\partial_x^kF(t,x;s) \big| \, ds   \\
&\lesssim |h|^{(\gamma+N-|k|_\mfs)/\mfs_i}\int_a^{t-|h|^{N/\mfs_i}}
e^{c_0(t-s)}(t-s)^{(\beta-\gamma-N)/N}ds
\lesssim |h|^{(\beta+N-|k|_\mfs)/\mfs_i}.
\end{align*}
Thus we obtain the desired estimate for $\partial_x^kF(t,x+he_i)-\partial_x^kF(t,x)$.
\end{Dem}

\vfill \pagebreak

\subsection{Uniqueness of the fundamental solution{\boldmath $.$} \hspace{0.15cm}}
\label{AppendixUniqueness}

We prove the uniqueness of the fundamental solution $Q_t(x,y)$ of the operator \eqref{eq:parabolicoperator} by the same argument as the proof of Theorem {\it4.3} of \cite[Section I\!I\!I.2]{Eid69}. See also Lemma {\it6.1.2} of \cite{Lun95}.
For any time interval $I\subset(a,b)$, we define $C_{\mfs}^{1,N}(I\times\bbR^d)$ as the collection of bounded continuous functions $f$ on $I\times \bbR^d$ such that $\partial_tf$ and $\partial_x^kf$ for any $|k|_{\mfs}\le N$ are bounded and continuous.

\medskip

\begin{thm}\label{thm:GammaUnique}
Suppose that the coefficients $a_k(t,x)$ of the operator \eqref{eq:parabolicoperator} satisfies the assumptions of Theorem \ref{thm:GaussGamma}.
For any fixed $s\in(a,b)$ and $f\in C_b(\bbR^d)$, the Cauchy problem \eqref{defn:FundamentalSolution1}-\eqref{defn:FundamentalSolution2} has a unique solution $F\in C_{\mfs}^{1,N}([s,b)\times\bbR^d)$.
Consequently, the fundamental solution $Q_{ts}(x,y)$ of the operator \eqref{eq:parabolicoperator} is unique (up to Lebesgue null sets in $y$).
\end{thm}

\medskip

\begin{Dem}
It is sufficient to show that the solution $F$ of \eqref{defn:FundamentalSolution1}-\eqref{defn:FundamentalSolution2} with $f=0$ is equal to zero. Set
$$
W(t)\defeq\sum_{|k|_{\mfs}\le N}\int_s^t\|\partial_x^kF(r,\cdot)\|_{C_b(\bbR^d)}dr.
$$
We fix a point $y\in \bbR^d$ and write the equation in the form
\begin{align*}
\big(\partial_t-P(t,y,\partial_x)\big)F(t,x)=\big(P(t,x,\partial_x)-P(t,y,\partial_x)\big)F(t,x)\eqdef f^y(t,x).
\end{align*}
Note that the fundamental solution of the operator \eqref{eq:parabolicoperator} is unique if the coefficients $a_k(t)$ are $x$-independent continuous functions of $t$. Thus we can write
\begin{align*}
F(t,x)=\int_{s}^t\int_{\bbR^d}Q_{tr}^y(x-z)f^y(r,z)dzdr,
\end{align*}
where $Q_{ts}^y$ is a fundamental solution of $\partial_t-P(t,y,\partial_x)$.
Then $\partial_xQ_{ts}^y\in{\sf G}^{-|k|_\mfs}$ for any $k\in\bbN^d$.
 We estimate the derivatives of $F$ by setting $y=x$. Since
$$
|f^x(r,z)|\lesssim \|z-x\|_{\mfs}^\alpha\sum_{|k|_{\mfs}\le N}\|\partial_x^kF(r,\cdot)\|_{C_b(\bbR^d)},
$$
using the Gaussian estimates of $\partial^kQ$ we have
\begin{align*}
\|\partial_x^kF(t,\cdot)\|_{C_b(\bbR^d)}
\lesssim \int_s^t(t-r)^{(\alpha-|k|_{\mfs})/N}\sum_{|\ell|_{\mfs}\le N}\|\partial_x^\ell F(r,\cdot)\|_{C_b(\bbR^d)}dr.
\end{align*}
By integration we can conclude that there exists a constant $C>0$ such that
\begin{align*}
W(t)\le C(t-s)^{\alpha/N}W(t).
\end{align*}
Hence it follows that $W(t)=0$ for any $s<t<s+t_0$, where $t_0\defeq (1/C)^{N/\alpha}>0$. 
Since the coefficients $a_k(t,x)$ are uniformly bounded and H\"older continuous, we can repeat the same argument as above by replacing $s$ with some $s'\in(t_0/2,t_0)$ and obtain that $W(t)=0$ for any $s'<t<s'+t_0$.
In the end, we can establish that $W(t)=0$ for any $t>s$.
\end{Dem}

\medskip

\subsection{Temporally homogeneous operator{\boldmath $.$} \hspace{0.15cm}}
\label{AppendixProperties}

\medskip

Next we consider the operator
\begin{align}\label{Appendix:timehomogeneousoperator}
\partial_t-P(x,\partial_x)
=\partial_t-\sum_{|k|_{\mfs}\le N}a_k(x)\partial_x^k
\end{align}
with $t$-independent coefficients $a_k(x)$.
Let $P$ satisfy the assumptions of Theorem \textit{\ref{thm:GaussGamma}}, and denote by $Q_t(x,y)$ be its fundamental solution defined for $t\in(0,\infty)$ and $x,y\in\bbR^d$. 
For any $f\in C_b(\bbR^d)$, we define the integral operator
$$
\big(Q_tf\big)(x)\defeq\int_{\bbR^d}Q_t(x,y)f(y)dy.
$$
It should be noted that $Q_t$ satisfies the semigroup property $Q_tQ_sf=Q_{t+s}f$. 

For any $\beta>0$, we denote by $C_{\mfs}^\beta(\bbR^d)$ the collection of $f\in C_b(\bbR^d)$ such that $\partial_x^kf$ is bounded and continuous for any $|k|_{\mfs}<\beta$, and if $|k|_{\mfs}<\beta\le|k+e_i|_{\mfs}$, we have
$$
\big| \partial^kf(x+he_i) - \partial^kf(x) \big| \lesssim |h|^{(\beta-|k|_\mfs)/\mfs_i}
$$
for any $x\in\bbR^d$ and $h\in\bbR$.

\medskip

\begin{prop}\label{smoothingofQt1}
Let $\beta\in(0,1]$, $f\in C^\beta(\bbR^d)$, and $|k|_\mfs\le N$.
\begin{enumerate}
\renewcommand{\theenumi}{(\roman{enumi})}
\renewcommand{\labelenumi}{(\roman{enumi})}
\setlength{\itemsep}{5pt}

\item\label{smoothingofQt1-1}
For any $t>0$, we have
$$
\|\partial_x^kQ_tf\|_{L^\infty(\bbR^d)} \lesssim \Big(e^{c_0t}t^{(\beta-|k|_\mfs)/N}+e^{a_0t}{\bf 1}_{k=0} \Big) \, \|f\|_{C_\mfs^\beta(\bbR^d)}.
$$
This estimate also holds for $\beta=0$ if we regard $C_\mfs^0(\bbR^d)$ as $C_b(\bbR^d)$.

\item
Let $i\in\{1,\dots,d\}$ and $\alpha_i\in[0,\mfs_i]\cap[0,\alpha+N-|k|_\mfs)$.
For any $t>0$, $x\in\bbR^d$, and $h\in\bbR$, we have
$$
\big| \partial_x^kQ_tf(x+he_i)-\partial_x^kQ_tf(x) \big|
\lesssim e^{c_0t}|h|^{\alpha_i/\mfs_i}t^{(\beta-|k|_\mfs-\alpha_i)/N}\|f\|_{C_\mfs^\beta(\bbR^d)}.
$$

\item\label{smoothingofQt1-3}
Let $\theta\in[\beta,N]$. For any $0<s<t\le1$, we have
$$
\|Q_tf-Q_sf\|_{L^\infty(\bbR^d)}
\lesssim |t-s|^{\theta/N}s^{(\beta-\theta)/N}\|f\|_{C_\mfs^\beta(\bbR^d)}.
$$

\item\label{smoothingofQt1-4}
Let $k\neq 0$ and $\theta\in[0,N]$. For any $0<s<t\le1$ we have
$$
\big\| \partial_x^kQ_tf - \partial_x^kQ_sf \big\|_{L^\infty(\bbR^d)}
\lesssim |t-s|^{\theta/N} s^{(\beta-|k|_\mfs-\theta)/N}\|f\|_{C_\mfs^\beta(\bbR^d)}.
$$
\end{enumerate}
\end{prop}

\medskip

\begin{Dem}
The first two estimates are essentially shown in the proof of Theorem {\it\ref{thm:InverseParabolic}}, so we show only \ref{smoothingofQt1-3} and \ref{smoothingofQt1-4}.

\ref{smoothingofQt1-3}:
For the case $\theta=N$, by using the inequality $\|\partial_tQ_tf\|_{L^\infty}\lesssim t^{(\beta-N)/N}\|f\|_{C_\mfs^\beta(\bbR^d)}$ following from $\partial_tQ_tf=P(x,\partial_x)f$, we have the desired inequality as follows.
\begin{align*}
\|Q_tf-Q_sf\|_{L^\infty(\bbR^d)}
&=\int_s^t\|\partial_rQ_rf\|_{L^\infty(\bbR^d)}dr\\
&\lesssim\|f\|_{C_\mfs^\beta(\bbR^d)}\int_s^tr^{(\beta-N)/N}dr
\lesssim\|f\|_{C_\mfs^\beta(\bbR^d)}|t-s|s^{(\beta-N)/N}.
\end{align*}
For the case $\theta=\beta$, we use the inequality 
\begin{equation}\label{smoothingofQt1-3-1}
\|(Q_t-\id)f\|_{L^\infty(\bbR^d)}\lesssim t^{\beta/N}\|f\|_{C_\mfs^\beta(\bbR^d)}
\end{equation}
obtained by performing a calculation similar to the proof of Theorem {\it\ref{thm:InverseParabolic}} for the decomposition
\begin{align*}
Q_tf(x)-f(x)
=\int_{\bbR^d}Q_t(x,y)\big(f(y)-f(x)\big)dy
+(e^{a_0t}-1)f(x).
\end{align*}
By using this inequality, we have
\begin{align*}
\|Q_tf-Q_sf\|_{L^\infty(\bbR^d)}
&= \big\| Q_s(Q_{t-s}-\id)f \big\|_{L^\infty(\bbR^d)}   \\
&\lesssim\|(Q_{t-s}-\id)f\|_{L^\infty(\bbR^d)}
\lesssim |t-s|^{\beta/N}\|f\|_{C_\mfs^\beta(\bbR^d)}.
\end{align*}
The inequality for $\theta\in(\beta,N)$ follows from an interpolation.

\ref{smoothingofQt1-4}:
For the case $\theta=N$, by using the semigroup property of $Q_t$ and \ref{smoothingofQt1-1} with $\beta=0$, we have
\begin{align*}
\big\| \partial_x^kQ_tf - \partial_x^kQ_sf \big\|_{L^\infty(\bbR^d)}
&= \big\|\partial_x^kQ_{s/2}(Q_{t-s/2}f - Q_{s/2}f) \big\|_{L^\infty(\bbR^d)}   \\
&= s^{-|k|_\mfs/N} \big\| Q_{t-s/2}f - Q_{s/2}f \big\|_{L^\infty(\bbR^d)}   \\
&\lesssim\|f\|_{C_\mfs^\beta(\bbR^d)} \, |t-s| \, s^{(\beta-|k|_\mfs-N)/N}.
\end{align*}
The inequality for the case $\theta=0$ immediately follows from \ref{smoothingofQt1-1}.
The generic case is obtained by an interpolation.
\end{Dem}

\medskip

For $\beta<0$, we define the space $C_{\mfs}^\beta(P)$ as the completion of the set of $f\in C_b(\bbR^d)$ under the norm
$$
\|f\|_{C_{\mfs}^\beta(P)}\defeq \sup_{0<t\le1}t^{-\beta/N} \|Q_tf\|_{L^\infty(\bbR^d)}.
$$

\medskip

\begin{thm}\label{thm:SchauderGeneral}
Let $Q_t$ satisfy the estimates in Definition \ref{DefnClassGZetaBeta} with some constant $c_0$.
For any $c>c_0$ and any $f\in C_b(\bbR^d)$, we define
$$
\big(c-P(x,\partial_x)\big)^{-1}f(x)\defeq
\int_0^\infty e^{-ct}Q_tf(x)dt.
$$
Then the map $(c-P(x,\partial_x))^{-1}$ is continuously extended to the map from $C_{\mfs}^\beta(P)$ into $C_{\mfs}^{\beta+N}(\bbR^d)$ 
for any $\beta\in(-N,0)$ such that $\beta+N$ is not an integer.
\end{thm}

\medskip

\begin{Dem}
We write $P=P(x,\partial_x)$ and $Q_t^c=e^{-ct}Q_t$ for simplicity. By  the semigroup property and the Gaussian estimate of $Q_t$, we have
\begin{align*}
\|\partial_x^kQ_t^cf\|_{L^\infty(\bbR^d)} = \big\| \partial_x^kQ_{t/2}^cQ_{t/2}^cf \big\|_{L^\infty(\bbR^d)}
\lesssim t^{-|k|_{\mfs}/N}\|Q_{t/2}^cf\|_{L^\infty(\bbR^d)}
\lesssim t^{(\beta-|k|_{\mfs})/N}\|f\|_{C_{\mfs}^\beta(P)}
\end{align*}
for any $t\in(0,2]$, and
\begin{align*}
\|\partial_x^kQ_t^cf\|_{L^\infty(\bbR^d)} = \big\| \partial_x^kQ_{t-1}^cQ_1^cf \big\|_{L^\infty(\bbR^d)}
\lesssim e^{-(c-c_0)(t-1)}\|f\|_{C_{\mfs}^\beta(P)}
\end{align*}
for any $t\ge2$. By integration, for any $|k|_{\mfs}<\beta+N$, we have
\begin{align*}
\big\| \partial_x^k(c-P)^{-1}f \big\|_{L^\infty(\bbR^d)} \leq \int_0^2\|\partial_x^kQ_t^c f\|_{L^\infty(\bbR^d)}\,dt + \int_2^\infty \|\partial_x^kQ_t^cf\|_{L^\infty(\bbR^d)} \, dt
\lesssim \|f\|_{C_{\mfs}^\beta(P)}.
\end{align*}
To show the H\"older estimates of $\partial_x^k(c-P)^{-1}f$ with $|k|_{\mfs}<\beta+N<|k+e_i|_\mfs$, it is sufficient to consider $x'=x+he_i$ with $|h|<2$. We decompose
\begin{align*}
\partial_x^k(c-P)^{-1}f(x')-\partial_x^k(c-P)^{-1}f(x)
=\int_0^\infty \int_{\bbR^d}\Big\{\partial_x^kQ_{t_1}^c(x',y)-\partial_x^kQ_{t_1}^c(x,y)\Big\} \,Q_{t_0}^cf(y) \, dydt
\end{align*}
as before, where $t=t_0+t_1$ and $t_0\defeq\min\{t/2,1\}$. Setting $\gamma=\beta+N-|k|_\mfs\in(0,\mfs_i)$ and choosing sufficiently small $\varepsilon>0$ such that $\gamma+\varepsilon<N-|k|_\mfs+\alpha$, we have
\begin{align*}
&\big| \partial_x^k(c-P)^{-1}f(x')-\partial_x^k(c-P)^{-1}f(x) \big|   \\
&\le
|h|^{(\gamma-\varepsilon)/\mfs_i}
\int_0^{|h|}t^{(-|k|_{\mfs}-(\gamma-\varepsilon)+\beta)/N}dt
+
|h|^{(\gamma+\varepsilon)/\mfs_i}
\int_{|h|}^2t^{(-|k|_{\mfs}-(\gamma+\varepsilon)+\beta)/N}dt\\
&\quad+
|h|^{\gamma/\mfs_i}
\int_2^\infty t^{(-|k|_{\mfs}-\gamma)/N}e^{-(c-c_0)(t-1)}dt
\lesssim |h|^{\gamma/\mfs_i}
\end{align*}
by the H\"older estimates of $Q_t$.
\end{Dem}

\medskip

\subsection{Anisotropic Taylor formula{\boldmath $.$} \hspace{0.15cm}}
\label{AppendixTaylor}

Continuing the previous section, we consider the temporally homogeneous operator \eqref{Appendix:timehomogeneousoperator}. In what follows, we consider the parabolic scaling $\mfs=(2,1,1,\dots,1)$. We denote by
$$
x=(x_1,x_2,\dots,x_d)\eqdef(x_1,\overline{x})
$$
a generic element of $\bbR^d$.
The following anisotropic Taylor formula is an analogue of Proposition {\it A.1} of \cite{Hai14}, but it should be noted that the fundamental solutions we discuss here have only restricted differentiability.

\medskip

\begin{prop}\label{PropAnisotropicTaylor}
Let $n\in\bbN$.
For any function $f$ on $\bbR^d$ which is $k$-times differentiable for any $|k|_{\mfs}\le n$, we have
\begin{align*}
\bigg|f(y)-\sum_{|k|_{\mfs}\le n}\frac{(y-x)^k}{k!}\partial_x^kf(x)\bigg|
&\lesssim \|y-x\|_{\mfs}^{n-1}\sup_{|k|_{\mfs}=n-1}\sup_{(z_1,\bar{z})}\big|\partial^kf(z_1,\bar{z})-\partial^kf(x_1,\bar{z})\big|\\
&\quad+\|y-x\|_{\mfs}^n\sup_{|k|_{\mfs}=n}\sup_{(z_1,\bar{z})}\big|\partial^kf(z_1,\bar{z})-\partial^kf(x)\big|,
\end{align*}
where $z_1$ (resp.  $\bar{z}$) runs over the interval $(x_1,y_1)$ (resp. $(\overline{x},\bar{y})$).
\end{prop}

\medskip

\begin{proof}
We define $A\defeq\{k\in\bbN^d\, ;\, |k|_{\mfs}\le n\}$ and
$$
A^\circ \defeq \Big\{k\in A\, ;\, k+e_i\in A\ \text{for all}\ i\in\{1,\dots,d\}\Big\} = \big\{k\in\bbN^d\, ;\, |k|_{\mfs}\le n-2\big\}.
$$
For $\theta\in[0,1]$, we set
$$
x(\theta)= (x_1(\theta),x'(\theta)) \defeq x+\theta(y-x).
$$ 
We write $\partial_i\defeq\partial_{x_i}$. By repeating the Taylor expansion of first order, we have
\begin{align*}
f(y)-f(x) &= \sum_{e_i\in A^\circ}(y_i-x_i)\int_0^1\partial_if(x(\theta))\,d\theta+\sum_{e_i\notin A^\circ}(y_i-x_i)\int_0^1\partial_if(x(\theta))\,d\theta   \\
&=\sum_{e_i\in A^\circ}(y_i-x_i)\partial_if(x)+\sum_{e_i+e_j\in A^\circ}(y-x)^{e_i+e_j}\int_0^1(1-\theta)\partial_{ij}f(x(\theta))\,d\theta   \\
&\qquad+\sum_{\substack{e_i\in A^\circ \\ e_i+e_j\notin A^\circ}}(y-x)^{e_i+e_j}\int_0^1(1-\theta)\partial_{ij}f(x(\theta))\,d\theta+\sum_{e_i\notin A^\circ}(y_i-x_i)\int_0^1\partial_if(x(\theta))\,d\theta   \\
&=\cdots\\
&=\sum_{k\in A^\circ\setminus\{0\}}\frac{(y-x)^k}{k!}\partial^kf(x)
+\sum_{\substack{k\in A^\circ \\ k+e_i\notin A^\circ}}\frac{(y-x)^{k+e_i}}{k!}\int_0^1(1-\theta)^{|k|}\partial^{k+e_i}f(x(\theta))\,d\theta
\end{align*}
in the end, where $|k|\defeq k_1+\cdots+k_d$.
In other words
\begin{align*}
f(y)
&=\sum_{|k|_{\mfs}\le n-2}\frac{(y-x)^k}{k!}\partial^kf(x)
+\sum_{|k|_{\mfs}=n-1}\frac{(y-x)^{k}}{k!}\int_0^1|k|(1-\theta)^{|k|-1}\partial^kf(x(\theta))\,d\theta\\
&\qquad+\sum_{|k|_{\mfs}=n-2}\frac{(y-x)^{k+e_1}}{k!}\int_0^1(1-\theta)^{|k|}\partial^{k+e_1}f(x(\theta))\,d\theta.
\end{align*}
We treat the last two terms by taking into account the restriction on differentiability. For $|k|_{\mfs}=n-1$, we decompose
\begin{align*}
\partial^kf(x(\theta))
&= \partial^kf(x)+\big(\partial^kf(x_1(\theta),\overline{x}(\theta))-\partial^kf(x_1,\overline{x}(\theta))\big)
+\big(\partial^kf(x_1,\overline{x}(\theta))-\partial^kf(x_1,\overline{x})\big)
\\
&=\partial^kf(x)+\hspace{-0.08cm}\big(\partial^kf(x_1(\theta),\overline{x}(\theta))-\partial^kf(x_1,\overline{x}(\theta))\big)
\hspace{-0.1cm} + \hspace{-0.1cm}\sum_{i=2}^d(y_i-x_i) \hspace{-0.08cm}\int_0^\theta\partial^{k+e_i}f(x_1,\overline{x}(\theta'))\,d\theta'.
\end{align*}
Thus we have
\begin{align*}
f(y)-&\sum_{|k|_{\mfs}\le n-1}\frac{(y-x)^k}{k!}\partial^kf(x)   \\
&=\sum_{|k|_{\mfs}=n-1}\frac{(y-x)^{k}}{k!}\int_0^1|k|(1-\theta)^{|k|-1} \Big\{\partial^kf(x_1(\theta),\overline{x}(\theta))-\partial^kf(x_1,\overline{x}(\theta)) \Big\} \, d\theta   \\
&\qquad+\sum_{|k|_{\mfs}=n-1}\sum_{i=2}^d\frac{(y-x)^{k+e_i}}{k!}\int_0^1(1-\theta)^{|k|}\partial^{k+e_i}f(x_1,\overline{x}(\theta))d\theta\\
&\qquad+\sum_{|k|_{\mfs}=n-2}\frac{(y-x)^{k+e_1}}{k!}\int_0^1(1-\theta)^{|k|}\partial^{k+e_1}f(x(\theta))d\theta,
\end{align*}
and
\begin{align*}
f(y)-&\sum_{|k|_{\mfs}\le n}\frac{(y-x)^k}{k!}\partial^kf(x)   \\
&=\sum_{|k|_{\mfs}=n-1}\frac{(y-x)^{k}}{k!}\int_0^1|k|(1-\theta)^{|k|-1} \Big\{ \partial^kf(x_1(\theta),\overline{x}(\theta))-\partial^kf(x_1,\overline{x}(\theta)) \Big\} \, d\theta   \\
&\qquad+\sum_{|k|_{\mfs}=n-1}\sum_{i=2}^d\frac{(y-x)^{k+e_i}}{k!}\int_0^1(1-\theta)^{|k|} \Big\{\partial^{k+e_i}f(x_1,\overline{x}(\theta))-\partial^{k+e_i}f(x)\Big\} \, d\theta   \\
&\qquad+\sum_{|k|_{\mfs}=n-2}\frac{(y-x)^{k+e_1}}{k!}\int_0^1(1-\theta)^{|k|} \Big\{ \partial^{k+e_1}f(x(\theta))-\partial^{k+e_1}f(x) \Big\} \, d\theta.
\end{align*}
This provides the desired inequality.
\end{proof}

\medskip

\begin{cor}\label{cor:TaylorFundamentalSolution}
Let $Q_t(x,y)$ be the fundamental solution of the operator \eqref{Appendix:timehomogeneousoperator} satisfying the assumptions of Theorem \textit{\ref{thm:GaussGamma}}.
Then for any $k\in\bbN^d$ with $|k|_\mfs\le N$ and $\beta\in(0,\alpha)$,
\begin{align*}
&\bigg|\partial_{x'}^kQ_t(x',y)-\sum_{|k+\ell|_{\mfs}\le N}\frac{(x'-x)^{\ell}}{\ell!}\partial_x^{k+\ell}Q_t(x,y)\bigg|\\
&\lesssim e^{c_0t}\|x'-x\|_{\mfs}^{N+\beta-|k|_{\mfs}}
\Big\{{\sf G}_t^{(c_1,-N-\beta)}(x'-y)+{\sf G}_t^{(c_1,-N-\beta)}(x-y)\Big\}.
\end{align*}
\end{cor}

\medskip

\begin{Dem}
We apply Proposition {\it\ref{PropAnisotropicTaylor}} to $f(x)=\partial_x^kQ_t(x,y)$ and $n=N-|k|_{\mfs}$.
By Theorem {\it\ref{thm:GaussGamma}},
\begin{align*}
\big|\partial_x^mQ_t((x_1',\overline{x}),y)-\partial_x^mQ_t((x_1,\overline{x}),y)\big|
&\lesssim e^{c_0t}|x_1'-x_1|^{(1+\beta)/2} \, {\sf G}_t^{(c_1,-N-\beta)}(x-y)
\end{align*}
for any $|m|_{\mfs}=N-1$, and
\begin{align*}
\big|\partial_x^mQ_t(x',y)-\partial_x^mQ_t(x,y)\big|
&\lesssim e^{c_0t}\|x'-x\|_{\mfs}^{\beta} \, {\sf G}_t^{(c_1,-N-\beta)}(x-y)
\end{align*}
for any $|m|_{\mfs}=N$.
\end{Dem}

\bigskip
\bigskip

\noindent \textcolor{gray}{$\bullet$} {\sf I. Bailleul} -- Univ Brest, CNRS, LMBA - UMR 6205, F- 29238 Brest, France.   \\
{\it E-mail}: ismael.bailleul@univ-brest.fr

\medskip

\noindent \textcolor{gray}{$\bullet$} {\sf M. Hoshino} --  Graduate School of Engineering Science, Osaka University, Japan.   \\
{\it E-mail}: hoshino@sigmath.es.osaka-u.ac.jp

\medskip

\noindent \textcolor{gray}{$\bullet$} {\sf S. Kusuoka} --  Graduate School of Science, Kyoto University, Japan.   \\
{\it E-mail}: kusuoka@math.kyoto-u.ac.jp

\end{document}